\documentclass[11pt]{article}
\usepackage{amssymb, amsmath, euscript, verbatim, theorem, xr, sw,default,title,mathrsfs,array,xypic,graphicx,url,tabularx}
\usepackage[margin=32pt,font=footnotesize]{caption}

\begin{document}




\title{The Seiberg-Witten Equations on Manifolds with Boundary I: \\ The Space of Monopoles and Their Boundary Values}
\author{Timothy Nguyen\footnote{The author was supported by NSF grants DMS-0706967 and DMS-0805841.}}

\abstract{In this paper, we study the Seiberg-Witten equations on a compact $3$-manifold with boundary.  Solutions to these equations are called monopoles.  Under some simple topological assumptions, we show that the solution space of all monopoles is a Banach manifold in suitable function space topologies.  We then prove that the restriction of the space of monopoles to the boundary is a submersion onto a Lagrangian submanifold of the space of connections and spinors on the boundary.  Both these spaces are infinite dimensional, even modulo gauge, since no boundary conditions are specified for the Seiberg-Witten equations on the $3$-manifold.  We study the analytic properties of these monopole spaces with an eye towards developing a monopole Floer theory for $3$-manifolds with boundary, which we pursue in \cite{N2}.}

\tableofcontents

\section{Introduction}

The Seiberg-Witten equations, introduced by Witten in \cite{Wit}, yield interesting topological invariants of closed three and four-dimensional manifolds and have led to many important developments in low dimensional topology during the last two decades.  On a closed $4$-manifold $X$, the Seiberg-Witten equations are a system of nonlinear partial differential equations for a connection and spinor on $X$.  When $X$ is of the form $\R \times Y$ with $Y$ a closed $3$-manifold, a dimensional reduction leads to the $3$-dimensional Seiberg-Witten equations on $Y$.  These latter equations are referred to as the \textit{monopole equations}.  Solutions to these equations are called \textit{monopoles}. For both three and four-dimensional manifolds, the topological invariants one obtains require an understanding of the moduli space of solutions to the Seiberg-Witten equations.  On a closed $4$-manifold $X$, the Seiberg-Witten invariant for $X$ is computed by integrating a cohomology class over the moduli space of solutions.  On a closed $3$-manifold $Y$, the monopole invariants one obtains for $Y$ come from studying the monopole Floer homology of $Y$.  This involves taking the homology of a chain complex whose differential counts solutions of the Seiberg-Witten equations on $\R \times Y$ that connect two monopoles on $Y$.  For further background and applications, see e.g. \cite{KM}, \cite{Mo}, \cite{Ni}.

In this paper and its sequel \cite{N2}, we study the Seiberg-Witten equations on manifolds with boundary.  In this first paper, we study the monopole equations on $3$-manifolds with boundary, where no boundary conditions are specified for the equations. Specifically, as is done in the case of a closed $3$-manifold, we study the geometry of the space of solutions to the monopole equations.  However, unlike in the closed case, where one hopes to achieve a finite dimensional (in fact zero-dimensional) space of monopoles modulo gauge, the space of monopoles on a $3$-manifold with boundary, even modulo gauge, is infinite dimensional, since no boundary conditions are imposed.  Moreover, we study the space obtained by restricting the space of monopoles to the boundary.  Under the appropriate assumptions (see the main theorem), we show that the space of monopoles and their boundary values are each Banach manifolds in suitable function space topologies.  We should emphasize that studying the monopole equations on a $3$-manifold with boundary poses some rather unusual problems.  This is because the linearization of the $3$-dimensional Seiberg-Witten equations are not elliptic, even modulo gauge.  This is in contrast to the $4$-dimensional Seiberg-Witten equations, whose moduli space of solutions on $4$-manifolds with boundary has been studied in \cite{KM}.  What we therefore have in our situation is a nonelliptic, nonlinear system of equations with unspecified boundary conditions.  We will address the nonellipticity of these equations and other issues in the outline at the end of this introduction.

The primary motivation for studying the space of boundary values of monopoles is that the resulting space, which is a smooth Banach manifold under the appropriate hypotheses, provides natural boundary conditions for the Seiberg-Witten equations on $4$-manifolds with boundary.  More precisely, consider the Seiberg-Witten equations on a cylindrical $4$-manifold $\R \times Y$, where $\partial Y = \Sigma$, and let $Y'$ be any manifold such that $\partial Y' = -\Sigma$ and $Y' \cup_\Sigma Y$ is a smooth closed oriented Riemannian $3$-manifold.   We impose as boundary condition for the Seiberg-Witten equations on $\R \times Y$ the following: at every time $t \in \R$, the configuration restricted to the boundary slice $\{t\} \times \Sigma$ lies in the space of restrictions of monopoles on $Y'$, i.e., the configuration extends to a monopole on $Y'$.   This boundary condition has its geometric origins in the construction of a monopole Floer theory for the $3$-manifold with boundary $Y$.  We discuss these issues and the analysis behind the associated boundary value problem in the sequel \cite{N2}.

In order to state our main results, let us introduce some notation (see Section 2 for a more detailed setup). So that we may work within the framework of Banach spaces, we need to consider the completions of smooth configuration spaces in the appropriate function space topologies.  The function spaces one usually considers are the standard Sobolev spaces $H^{k,p}$ of functions with $k$ derivatives lying in $L^p$.  However, working with these spaces alone is inadequate because the space of boundary values of a Sobolev space is not a Sobolev space (unless $p = 2$).  Instead, the space of boundary values of a Sobolev space is a Besov space, and so working with these spaces will be inevitable when we consider the space of boundary values of monopoles.  Thus, while we may work with Sobolev spaces on $Y$, we are forced to work with Besov spaces on $\Sigma$.  However, to keep the analysis and notation more uniform, we will mainly work with Besov spaces on $Y$ instead of Sobolev spaces (though nearly all of our results adapt to Sobolev spaces on $Y$), which we are free to do since the space of boundary values of a Besov space is again a Besov space.  Moreover, since Besov spaces on $3$-manifolds will be necessary for the analysis in \cite{N2}, as $3$-manifolds will arise as boundaries of $4$-manifolds, it is essential that we state results here for Besov spaces and not just for Sobolev spaces.  On the other hand, there will be places where we want to explicitly restate\footnote{Besov spaces and Sobolev spaces are ``nearly identical" in the sense of Remark \ref{RemHB}, so that having proven results for one of these types of spaces, one automatically obtains them for the other type.} our results on Besov spaces in terms of Sobolev spaces (we will need both the Besov and Sobolev space versions of the analysis done in this paper in \cite{N2}), so that the separation of Besov spaces from Sobolev spaces on $Y$ is not completely rigid, see Remark \ref{RemFun}. With these considerations then, if $Y$ is a $3$-manifold with boundary $\Sigma$, consider the Besov spaces $B^{s,p}(Y)$ and $B^{s,p}(\Sigma)$, for $s \in \R$ and $p \geq 2$.  The definition of these spaces along with their basic properties are summarized in Appendix \ref{AppFun}.  When $p = 2$, we have $B^{s,2} = H^{s,2}$, the usual fractional order Sobolev space of functions with $s$ derivatives belonging to $L^2$ (usually denoted just $H^s$).  For $p \neq 2$, the Besov spaces are never Sobolev spaces of functions with $s$ derivatives in $L^p$.  The reader unfamiliar with Besov spaces can comfortably set $p = 2$ on a first reading of this paper.  Moreover, since we will be working with low fractional regularity, the reader may also set $s$ equal to a sufficiently large integer to make a first reading simpler.

Let $Y$ be endowed with a $\spinc$ structure $\s$.  The $\spinc$ structure yields for us the configuration space
$$\fC^{s,p}(Y) = \fC^{s,p}(Y,\s) := B^{s,p}(\A(Y) \times \Gamma(\S))$$
on which the monopole equations are defined.  Here $\Gamma(\S)$ is the space of smooth sections of the spinor bundle $\S = \S(\s)$ on $Y$ determined by $\s$, $\A(Y) = \A(Y,\s)$ is the space of smooth $\spinc$ connections on $\S$, and the prefix $B^{s,p}$ denotes that we have taken the $B^{s,p}(Y)$ completions of these spaces.  The monopole equations are defined by the equations
\begin{equation}
  SW_3(B,\Psi) = 0, \label{SW3-intro}
\end{equation}
where $SW_3$ is the Seiberg-Witten map given by (\ref{SW3map}).  Here, $s$ and $p$ are chosen sufficiently large so that these equations are well-defined (in the sense of distributions). Define
\begin{equation}
  \fM^{s,p}(Y,\s) = \{(B,\Psi): \in \fC^{s,p}(Y, \s): SW_3(B,\Psi) = 0\} \label{sec1:fM}
\end{equation}
to be space of all solutions to the monopole equations in $\fC^{s,p}(Y)$.  Fixing a smooth reference connection $B_\rf \in \A(Y)$, let
\begin{equation}
  \M^{s,p}(Y,\s) = \{(B,\Psi): \in \fC^{s,p}(Y, \s): SW_3(B,\Psi) = 0,\, d^*(B-B_\rf) = 0\}
\end{equation}
denote the space of $B^{s,p}(Y)$ monopoles in Coulomb gauge with respect to $B_\rf$.

On the boundary $\Sigma$, we can define the boundary configuration space in the $B^{s,p}(\Sigma)$ topology,
$$\fC^{s,p}(\Sigma) = \fC^{s,p}(\Sigma,\s) := B^{s,p}(\A(\Sigma) \times \Gamma(\S_\Sigma)),$$
where $\S_\Sigma$ is the bundle $\S$ restricted to $\Sigma$, and $\A(\Sigma)$ is the space of $\spinc$ connections on $\S_\Sigma$. For $s > 1/p$, we have a restriction map \label{p:rSigma}
\begin{align}
  r_\Sigma: \fC^{s,p}(Y) & \to \fC^{s-1/p,p}(\Sigma) \nonumber\\
  (B,\Psi) & \mapsto (B|_\Sigma,\Psi|_\Sigma) \label{sec1:r}
\end{align}
which restricts a connection $B \in \A(Y)$ and spinor $\Psi \in \Gamma(\S)$ to $\Sigma$.  Observe that when $p = 2$, this is the usual trace theorem on $H^s$ spaces, whereby the trace of an element of $H^s(Y)$ belongs to $H^{s-1/2}(\Sigma)$. Thus, we can define the space of boundary values of the space of monopoles
\begin{equation}
  \L^{s-1/p,p}(Y,\s) := r_\Sigma(\fM^{s,p}(Y,\s)) \subset \fC^{s-1/p,p}(\Sigma).
\end{equation}
We will refer to all the spaces $\fM^{s,p}$, $\M^{s,p}$, and $\L^{s-1/p,p}$ as \textit{monopole spaces}.

The boundary configuration space $\fC(\Sigma)$ carries a natural symplectic structure.  Indeed, the space $\fC(\Sigma)$ is an affine space modeled on $\Omega^1(\Sigma; i\R) \oplus \Gamma(\S_\Sigma)$, and we can endow $\fC(\Sigma)$ with the constant symplectic form
\begin{align}
\omega((a,\phi), (b,\psi)) = \int_\Sigma a \wedge b + \int_\Sigma\Re(\phi,\rho(\nu)\psi), \qquad  (a,\phi), (b,\psi) \in \Omega^1(\Sigma; i\R) \oplus \Gamma(\S_\Sigma).\label{omega0}
\end{align}
Here, $\rho(\nu)$ is Clifford multiplication by the outward unit normal $\nu$ to $\Sigma$ and the inner product on spinors is induced from the Hermitian metric on $\S_\Sigma$.  The symplectic form (\ref{omega0}) extends to a symplectic form on $\fC^{0,2}(\Sigma)$, the $L^2$ configuration space on the boundary.  Since $\fC^{s,p}(\Sigma) \subset \fC^{0,2}(\Sigma)$ when $s > 0$ and $p > 2$, these latter spaces are also symplectic Banach configuration spaces (in the sense of Appendix \ref{SecSymp}).

Let $\det(\s) = \Lambda^2\S(\s)$ denote the determinant line bundle of the spinor bundle and let $c_1(\s) = c_1(\det(\s))$ denote its first Chern class.  Then under suitable restrictions on $\s$ and $Y$, our main theorem gives us the following relations among our Besov monopole spaces\footnote{The main theorem also holds with the Besov space $B^{s,p}(Y)$ on $Y$ replaced with the function space $H^{s,p}(Y)$.  The space $H^{s,p}(Y)$ is a known as a \textit{Bessel-potential space} and for $s$ a nonnegative integer, $H^{s,p}(Y) = W^{s,p}(Y)$ is the usual Sobolev space of functions having $s$ derivatives in $L^p(Y)$, $1<p<\infty$.  Thus, the $H^{s,p}(Y)$ can be regarded as fractional Sobolev spaces for $s$ not an integer. See Appendix \ref{AppFun} and Remark \ref{RemFun}.}:\\

\noindent \textbf{Main Theorem. }\textit{Let $Y$ be a smooth compact oriented Riemannian $3$-manifold with boundary\footnote{We always assume the boundary $\Sigma$ to be nonempty. It need not be connected however.} $\Sigma$ and let $\s$ be a $\mathit{spin}^c$ structure on $Y$. Suppose either $c_1(\s)$ is non-torsion or else $H^1(Y,\Sigma)=0$.  Let $p \geq 2$ and $s > \max(3/p,1/2)$.  Then we have the following:
  \begin{enumerate}
    \item The spaces $\fM^{s,p}(Y,\s)$ and $\M^{s,p}(Y,\s)$ are closed\footnote{For Banach submanifolds modeled on an infinite dimensional Banach space (see Definition \ref{DefSub}), we use the adjective closed only to denote that the submanifold is closed as a topological subspace.  For finite dimensional manifolds, closed in addition means that the manifold is compact and has no boundary. As an infinite dimensional Banach manifold is never even locally compact, this distinction in terminology should cause no confusion.} Banach submanifolds of $\fC^{s,p}(Y)$.
    \item If furthermore, $s > 1/2+1/p$, then $\L^{s-1/p,p}(Y,\s)$ is a closed Lagrangian submanifold of $\fC^{s-1/p,p}(\Sigma)$.  The restriction maps
    \begin{align}
      r_\Sigma: \fM^{s,p}(Y,\s) & \to \L^{s-1/p,p}(Y,\s) \label{r1}\\
      r_\Sigma: \M^{s,p}(Y,\s) & \to \L^{s-1/p,p}(Y,\s) \label{r2},
    \end{align}
    are a submersion and covering map, respectively.  The fiber of (\ref{r2}) is isomorphic to the lattice $H^1(Y,\Sigma)$.  In particular, if $H^1(Y,\Sigma) = 0$, then (\ref{r2}) is a diffeomorphism.
    \item Smooth configurations are dense in $\fM^{s,p}(Y,\s)$, $\M^{s,p}(Y,\s)$, and $\L^{s-1/p,p}(Y,\s)$.\\
  \end{enumerate}}

Thus, in particular, our main theorem tells us that our monopole spaces are smooth Banach manifolds for a certain range of $s$ and $p$.  Let us make some remarks on the condition $s > \max(3/p,1/2)$. We need $s > 3/p$ because then $B^{s,p}(Y)$ embeds into the space $C^0(Y)$ of continuous functions on $Y$.  This allows us to use the unique continuation results stated in the appendix.  Unfortunately, for $p \leq 3$, this means we need $s > 1$, which does not seem optimal since the monopole equations only involve one derivative.  For $p > 3$, we can take $s < 1$, in which case, the monopole equations are defined only in a weak sense (in the sense of distributions).  We consider this low regularity case because it arises in the boundary value problem studied in \cite{N2}. Specifically, we will use the Lagrangian submanifold $\L^{s-1/p,p}$ as a boundary condition for the $4$-dimensional Seiberg-Witten equations.  Here, Lagrangian means that every tangent space to $\L^{s-1/p,p}$ is a Lagrangian subspace of the tangent space to $\fC^{s-1/p,p}(\Sigma)$, i.e., the tangent space to $\L^{s-1/p,p}$ is isotropic and has an isotropic complement with respect to the symplectic form (\ref{omega0}).  The Lagrangian property is important because it arises in the context of self-adjoint boundary conditions.  These issues will be further pursued in \cite{N2}.  We should note that the analysis of the monopole equations needs to be done rather carefully at low regularity, since managing the function space arithmetic that arises from multiplying low regularity configurations becomes an important issue. In fact, the low regularity analysis is unavoidable if one wishes to prove the Lagrangian property for $\L^{s-1/p,p}$, since we need to understand the family of symplectic configuration spaces $\fC^{s-1/p,p}(\Sigma)$ as lying inside the strongly symplectic configuration space $\fC^{0,2}(\Sigma)$, the space of $L^2$ configurations on $\Sigma$ (see Appendix \ref{SecSymp} and also Remark \ref{RemReg}).  If one does not care about the Lagrangian property, then the main theorem with $s$ large can be proven without having to deal with low regularity issues. At low regularity, the requirements $s > 1/2$ and $s > 1/2+1/p$ in the theorem are other technicalities that have to do with achieving transversality and obtaining suitable a priori estimates for monopoles (see Section 4). Let us also note that statement (iii) in the main theorem, which establishes the density of smooth monopoles in the monopole spaces, is not at all obvious.  Indeed, our monopole spaces are not defined to be Besov closures of smooth monopoles, but as seen in (\ref{sec1:fM}), they arise from the zero set of the map $SW_3$ defined on a Banach space of configurations.  This way of defining our monopole spaces is absolutely necessary if we are to use the essential techniques from Banach space theory, such as the inverse function theorem.  However, since our monopole spaces are not linear Banach spaces, and since they are infinite dimensional modulo gauge, some work must be done to show that a Besov monopole can be approximated by a smooth monopole.

Let us make the simple remark that our theorem is nonvacuous due to the following example:\\

\noindent \textbf{Example. }Suppose $c_1(\s)$ is torsion.  Then every flat connection on $\det(\s)$ yields a solution of the monopole equations (where the spinor component is identically zero).  If $H^1(Y,\Sigma) = 0$, the main theorem implies that the monopole spaces are smooth nonempty Banach manifolds.  In fact, using Theorem \ref{ThmMcharts}, one can describe a neighborhood of any configuration in the space of monopoles on $Y$, in particular, a neighborhood of a flat connection.\\

Our main theorem will be proven in Theorems \ref{ThmMMan} and \ref{ThmMtoL}.  In addition to these, we have Theorems \ref{ThmMcharts} and \ref{ThmLcharts}, which describe for us certain analytic properties of the local chart maps of our monopole spaces.  These properties are not only of interest in their own right, since our monopole spaces are infinite dimensional Banach manifolds, but they will play an essential role in \cite{N2}.

Finally, let us also remark that our methods, and hence our theorems, carry over straightforwardly if we perturb the Seiberg-Witten equations by a smooth coclosed $1$-form $\eta$. That is, we consider the equations
\begin{equation}
  SW_3(B,\Psi) = (\eta,0). \label{SW3-pert}
\end{equation}
We have the following result:\\

\noindent \textbf{Corollary. }\textit{Suppose either $c_1(\s) \neq \frac{i}{\pi}[*\eta]$ or else $H^1(Y,\Sigma) = 0$.  Then all the conclusions of the main theorem remain true for the monopole spaces associated to the perturbed monopole equations (\ref{SW3-pert}).\\}

\noindent Thus, for any $\s$, the corresponding perturbed spaces of monopoles will be smooth for generic coclosed perturbations. Moreover, these monopole spaces will be nonempty for many choices of $\eta$, since given any smooth configuration $(B,\Psi)$ such that $\Psi$ lies in the kernel of $D_B$, the Dirac operator determined by $B$ (see Section 2), we can simply define $\eta$ to be the value of $SW_3(B,\Psi)$, in which case $(B,\Psi)$ automatically solves (\ref{SW3-pert}).\\

\noindent\textit{Outline of Paper:} This paper is organized as follows.  In Section 2, we define the basic setup for the monopole equations on $Y$.  In Section 3, we establish the foundational analysis to handle the linearization of the monopole equations.  This primarily involves understanding the various gauge fixing issues involved as well as understanding how elliptic operators behave on manifolds with boundary.  The presence of a boundary makes this latter issue much more difficult than the case when there is no boundary.  Indeed, on a closed manifold, elliptic operators are automatically Fredholm when acting between standard function spaces (e.g. Sobolev spaces and Besov spaces).  On the other hand, on a manifold with boundary, the kernel of an elliptic operator is always infinite dimensional.  To fully understand the situation, we need to use the pseudodifferential tools summarized in Appendix \ref{AppEBP}, which allows us to handle elliptic boundary value problems on a variety of function spaces, in particular, Besov spaces of low regularity.  From this, what we will find is that the tangent spaces to our monopole spaces are given essentially by the range of pseudodifferential projections. Having established the linear theory, we use it in Section 4 to study the nonlinear monopole equations and prove our main results concerning the monopole spaces.  The appendix summarizes and synthesizes many of the analytic results needed in the main body of the paper, including the basic definitions and properties of the function spaces we use.

As we pointed out earlier, the linearization of the $3$-dimensional Seiberg-Witten equations is unfortunately not elliptic, even modulo gauge.  To work around this, we embed these equations into an elliptic system and use tools from elliptic theory to derive results for the original equations from the enlarged system.  This procedure is described in Section 3.3, where issues regarding ellipticity and gauge-fixing intertwine.  Furthermore, when we restrict to the boundary, passing from the enlarged elliptic system back to the original non-elliptic system involves a symplectic reduction, and so there is also an important interplay of symplectic functional analysis in what we do.\\

\noindent\textit{Acknowledgements: }This paper and its sequel are outputs of the author's PhD thesis.  The author would like to thank his thesis advisors Tom Mrowka and Katrin Wehrheim for their guidance in this project.  He would also like to thank David Jerison and Richard Melrose for helpful discussions.\\

\section{The Basic Setup}

We give a quick overview of the setup for the Seiberg-Witten equations on a $3$-manifold.  For a more detailed setup, see \cite{KM}.  Let $Y$ be a smooth compact oriented Riemannian $3$-manifold with boundary $\Sigma$.  A $\spinc$ structure $\s$ on $Y$ is a choice of $U(2)$ principal bundle over $Y$ that lifts the $SO(3)$ frame bundle of $Y$.  The space of all $\spinc$ structures on $Y$ is a torsor over $H^1(Y;\mathbb{Z})$.  Any given $\spinc$ structure $\s$ determines for us a spinor bundle \label{p:s} $\S = \S(\s)$ over $Y$, which is the two-dimensional complex vector bundle over $Y$ associated to the $U(2)$ bundle corresponding to $\s$. Endow $\S$ with a Hermitian metric. From this, we obtain Clifford multiplication bundle maps $\rho: TY \to \mathrm{End}(\S)$ and $\rho: T^*Y \to \mathrm{End}(\S)$, where the two are intertwined by the fact that the Riemannian metric gives a canonical isomorphism $TY \cong T^*Y$.  The map $\rho$ extends complex linearly to a map on the complexified exterior algebra of $T^*Y$ and we choose $\rho$ so that $\rho$ maps the volume form on $Y$ to the identity automorphism on $\S$.  This determines the spinor bundle $\S = (\S,\rho)$ uniquely up to isomorphism.

Fix a $\spinc$ structure $\s$ for the time being on $Y$. Only later in Section \ref{SecMonopoles} will be impose restrictions on $\s$.  A $\spinc$ connection on $\S$ is a Hermitian connection $\nabla$ on $\S$ for which Clifford multiplication is parallel, i.e., for all $\Psi \in \Gamma(\S)$ and $e \in \Gamma(TY)$, we have $\nabla (\rho(e)\Psi) = \rho(\nabla^{LC}e)\Psi + \rho(e)\nabla\Psi$, where $\nabla^{LC}$ denotes the Levi-Civita connection. Let \pageref{p:cA} $\A(Y)$ denote the space of $\spinc$ connections $\A(Y)$ on $Y$.  The difference of any two $\spinc$ connections acts on a spinor via Clifford multiplication by an imaginary-valued $1$-form.  Thus, given any fixed $\spinc$ connection $B_0 \in \A(Y)$, we can identify \label{p:B}
$$\A(Y) = \{B_0 + b : b \in \Omega^1(Y; i\R)\},$$
so that $\A(Y)$ is an affine space over $\Omega^1(Y; i\R)$.

Let \label{p:fC}
$$\fC(Y) = \fC(Y,\s) = \A(Y) \times \Gamma(\S)$$
denote the configuration space of all smooth $\spinc$ connections and smooth sections of the spinor bundle $\S$.  It is an affine space modeled on $\Omega^1(Y;i\R) \oplus \Gamma(\S)$.  By abuse of notation we let the inner product $(\cdot,\cdot)$ denote the following items: the Hermitian inner product on $\S$, linear in the first factor, the Hermitian inner product on complex differential forms induced from the Riemannian metric on $Y$, and finally the real inner product on $\Omega^1(Y;i\R) \oplus \Gamma(\S)$ induced from the real part of the inner products on each factor.

The Seiberg-Witten equations on $Y$ are given by the pair of equations
\begin{equation}
  \begin{split}
    \frac{1}{2}*F_{B^t} + \rho^{-1}(\Psi\Psi^*)_0 &= 0\\
      D_B\Psi &= 0,
  \end{split}\label{monopoles}
\end{equation}
where $(B,\Psi) \in \fC(Y)$.  Here $B^t$ is the connection induced from $B$ on the determinant line bundle $\det(\s) = \Lambda^2(\S)$ of $\S$, the element $F_{B^t} \in \Omega^2(Y, i\R)$ is its curvature, and $*$ is the Hodge star operator on $Y$.  For any spinor $\Psi$, the term $(\Psi\Psi^*)_0 \in \mathrm{End}(\S)$ is the trace-free Hermitian endomorphism of $\S$ given by the trace-free part of the map $\varphi \mapsto (\varphi,\Psi)\Psi$.  Since $\rho$ maps $\Omega^1(Y; i\R)$ isomorphically onto the space of trace-free Hermitian endomorphisms of $\S$, then $\rho^{-1}(\Psi\Psi^*)_0 \in \Omega^1(Y; i\R)$ is well-defined.  Finally, $D_B: \Gamma(\S) \to \Gamma(\S)$ is the $\spinc$ Dirac operator associated to the $\spinc$ connection $B$, i.e., in local coordinates, we have $D_B = \sum_{i=1}^3\rho(e_i)\nabla_{B,e_i}$ where $\nabla_B$ is the $\spinc$ covariant derivative associated to $B$ and the $e_i$ form a local orthonormal frame of tangent vectors.

Altogether, the left-hand side of (\ref{monopoles}) defines for us a Seiberg-Witten map \label{p:SW3}
\begin{align}
  SW_3: \fC(Y) & \to \Omega^1(Y;i\R) \times \Gamma(\S) \nonumber\\
  (B,\Psi) & \mapsto \left(\frac{1}{2}*F_{B^t} + \rho^{-1}(\Psi\Psi^*)_0, D_B\Psi\right). \label{SW3map}
\end{align}
Thus, solutions to the Seiberg-Witten equations are precisely the zero set of the map $SW_3$.  We will refer to a solution of the Seiberg-Witten equations as a \textit{monopole}. Let \label{p:fM}
\begin{equation}
  \fM(Y,\s) = \{(B,\Psi) \in \fC(Y) : SW_3(B,\Psi) = 0\} \label{smoothfM}
\end{equation}
denote the solution space of all monopoles on $Y$.  Fixing a smooth reference connection $B_\rf \in \A(Y)$ once and for all, let \label{p:M}
\begin{equation}
  \M(Y,\s) = \{(B,\Psi) \in \fC(Y) : SW_3(B,\Psi) = 0, d^*(B-B_\rf) = 0\} \label{smoothM}
\end{equation}
denote the space of all smooth monopoles that are in Coulomb gauge with respect to $B_\rf$.   Without any assumptions, the spaces $\fM(Y,\s)$ and $\M(Y,\s)$ are just sets, but we will see later, by transversality arguments, that these spaces of monopoles are indeed manifolds under suitable assumptions on $Y$ and $\s$. Since $\partial Y = \Sigma$ is nonempty and no boundary conditions have been specified for the equations defining $\fM(Y,\s)$ and $\M(Y,\s)$, these spaces will be infinite dimensional, even modulo the full gauge group.  Note that the space $\M(Y,\s)$ is obtained from $\fM(Y,\s)$ through a partial gauge-fixing, see Section \ref{SecGG}.

Let the boundary $\Sigma$ be given the usual orientation induced from that of $Y$, i.e., if $\nu$ is the outward normal vector field along $\Sigma$ and $dV$ is the oriented volume form on $Y$, then $\nu \llcorner dV$ yields the oriented volume form on $\Sigma$. On the boundary $\Sigma$, we have the configuration space
$$\fC(\Sigma) = \A(\Sigma) \times \Gamma(\S_\Sigma),$$
where $S_\Sigma$ is the bundle $\S$ restricted to $\Sigma$, and $\A(\Sigma)$ is the space of $\spinc$ connections on $S_\Sigma$. We have a restriction map
\begin{align}
  r_\Sigma: \fC(Y) & \to \fC(\Sigma) \nonumber \\
  (B,\Psi) & \mapsto (B|_\Sigma,\Psi|_\Sigma) \label{rSigma}
\end{align}
From this, we can define the space of (tangential) boundary values of the space of monopoles \label{p:L}
\begin{equation}
  \L(Y,\s) = r_\Sigma(\fM(Y,\s)). \label{sec2:L=rM}
\end{equation}

Observe that the space $\L(Y,\s)$ is nonlocal in the sense that its elements, which belong to $\fC(\Sigma)$, are not defined by equations on $\Sigma$.  Indeed, $\L(Y,\s)$ is determined by the full Seiberg-Witten equations in the interior of the manifold.  This makes the analysis concerning the manifold $\L(Y,\s)$ rather delicate, since one has to control both the space $\fM(Y,\s)$ and the behavior of the map $r_\Sigma$.

Ultimately, we want our manifolds to be Banach manifolds, and so we must complete our smooth configuration spaces in the appropriate function space topologies.  As explained in the introduction, the topologies most suitable for us are the Besov spaces $B^{s,p}(Y)$ and $B^{s,p}(\Sigma)$ on $Y$ and $\Sigma$, respectively, where $s \in \R$ and $p \geq 2$. These are the familiar $H^s$ spaces when $p = 2$  and for $p \neq 2$, the Besov spaces are never Sobolev spaces, i.e., spaces of functions with a specified number of derivatives lying in $L^p$.  Nevertheless, much of the analysis we will do applies to Sobolev spaces as well, since the analysis of elliptic boundary value problems is flexible and applies to a wide variety of function spaces.  To keep the notation minimal, we work mainly with Besov spaces and make a general remark at the end about how statements generalize to Sobolev spaces and other spaces (see Remark \ref{RemFun}). The Besov spaces, other relevant function spaces, and their properties are summarized in the appendix. On a first reading, one may set $p = 2$ and $s$ a large number, say a large integer, wherever applicable, so that the function spaces are as familiar as desired.

Thus, for $p \geq 2$ and $s \in \R$, we consider the Besov spaces $B^{s,p}(Y)$ and $B^{s,p}(\Sigma)$ of scalar-valued functions on $Y$ and $\Sigma$, respectively.  These topologies induce topologies on vector bundles over $Y$ and $\Sigma$ in the natural way, and so we may define the Besov completions of the configuration spaces \label{p:Csp}
\begin{align}
  \fC^{s,p}(Y) &= B^{s,p}(Y) \textrm{ closure of } \fC(Y)\\
  \fC^{s,p}(\Sigma) &= B^{s,p}(\Sigma) \textrm{ closure of } \fC(\Sigma).
\end{align}
Of course, when defining Besov norms on the space of connections in the above, we have to first choose a (smooth) reference connection, which then identifies the Besov space of connections with the Besov space of $1$-forms.

For $s,p$ such that the Seiberg-Witten equations make sense on $\fC^{s,p}(Y)$ (in the sense of distributions), we have the monopole spaces
\begin{align}
  \fM^{s,p}(Y,\s) & = \{(B,\Psi) \in \fC^{s,p}(Y) : SW_3(B,\Psi) = 0\}\\
  \M^{s,p}(Y,\s) & = \{(B,\Psi) \in \fC^{s,p}(Y) : SW_3(B,\Psi) = 0, d^*(B-B_\rf)=0.\}
\end{align}
in $\fC^{s,p}(Y)$.  Observe that for the range of $s$ and $p$ that are relevant for us, namely $p \geq 2$ and $s > \max(3/p,1/2)$, the Seiberg-Witten equations are well-defined on $\fC^{s,p}(Y)$.  This follows from Corollary \ref{CorDiffMap} and Theorem \ref{ThmMult}.

For $s > 1/p$, the restriction map (\ref{rSigma}) extends to a map
\begin{equation}
  r_\Sigma: \fC^{s,p}(Y) \to \fC^{s-1/p,p}(\Sigma), \label{BrSigma}
\end{equation}
and so we can define
$$\L^{s-1/p,p}(Y,\s) := r_\Sigma(\fM^{s,p}(Y,\s)).$$
Having defined our monopole spaces in the relevant topologies, we now begin the study of their properties as Banach manifolds.  With $\s$ and $Y$ fixed, we will often write $\fM(Y)$ or simply $\fM$ instead of $\fM(Y,\s)$.  Likewise for the other monopole spaces.

\begin{Rem}\label{RemReg}
  In the remainder of this paper, we will be stating results for various values of $s$ and $p$.  Unless stated otherwise, we will always assume
  \begin{equation}
    2 \leq p < \infty.
  \end{equation}
  Many of the statements of this paper are phrased in such a way that the range of permissible $s$ and $p$ is quite large and moreover, several topologies are often simultaneously  involved (e.g. Lemma \ref{LemmaTdecomp}).  This is not merely an exercise in function space arithmetic and there are several important reasons for stating our results this generally.

  First, we will need to work in the low regularity regime with $s < 1$ for applications in \cite{N2}.  In particular, when a first order operator acts on a configuration with regularity $s < 1$, we obtain a configuration with negative regularity and hence our results must be stated in enough generality to account for this.  Second, as mentioned in the introduction, the Lagrangian property of $\L^{s-1/p,p}$, even at high regularity (i.e. large $s$), requires an analysis of the $\L^{s-1/p,p}$ at low regularity.  Indeed, among all the spaces $\fC^{s,p}(\Sigma)$, only $\fC^{0,2}(\Sigma)$ is modeled on a strongly symplectic Hilbert space (see Appendix \ref{SecSymp}), and we will need to study all the symplectic spaces $\fC^{s,p}(\Sigma)$, $s > 0$, as subspaces of the space $\fC^{0,2}(\Sigma)$.  Thus, in a fundamental way, we will generally be considering multiple topologies simultaneously.  Observe that from these considerations, it is necessary to have the pseudodifferential tools summarized in Appendix \ref{AppEBP}.  Indeed, we need to understand elliptic boundary value problems at low (even negative) regularity, and furthermore, we have to deal with the fact that there is no trace map $\fC^{1/2,2}(Y) \to \fC^{0,2}(\Sigma)$.

  Hence, it is natural to state our results for a range of $s$ and $p$ that are as flexible as possible.  In fact, based on the function space arithmetic alone, many of the proofs involved are natural for the range $s > 3/p$ say (since then $B^{s,p}(Y)$ is an algebra), and it would be unnatural to restrict the range of $s$ based on the particular applications we have in mind.  Finally, it may be desirable to sharpen the range of $s$ and $p$ considered in this paper and so we try to state our results in a sufficiently general way at the outset.
\end{Rem}

\noindent \textbf{Notation.} Given any space $\frak{X}$ of configurations over a manifold $X = Y$ or $\Sigma$, we write $B^{s,p}\frak{X}$ to denote the closure of $\frak{X}$ with respect to the $B^{s,p}(X)$ topology.  We define $L^p\frak{X}$, $C^0\frak{X}$, and $H^{s,p}\frak{X}$ similarly.  For brevity, we may refer to just the function space which defines the topology of a configuration, e.g.,  if $\frak{X}$ is a space of configurations on $Y$, we may say an element $u \in B^{s,p}\frak{X}$ belongs to $B^{s,p}(Y)$ or just $B^{s,p}$ for short. If $E$ is a vector bundle over a space $X$, we write $B^{s,p}(E)$ as shorthand for $B^{s,p}\Gamma(E)$, the closure of the space $\Gamma(E)$ of smooth sections of $E$ in the topology $B^{s,p}(X)$.  If $X$ has boundary, we write $E_{\partial X}$ to denote $E|_{\partial X}$, the restriction of the bundle $E$ to the boundary $\partial X$.

From now on, we will make free use of the basic properties of the function spaces employed in this paper (multiplication and embedding theorems in particular), all of which can be found in the appendix.\\

\section{Linear Theory}\label{SecLinearTheory}

To study our monopole spaces, we first study their linearization, that is, their formal tangent spaces.  This involves studying the linearization of the Seiberg-Witten map.  Furthermore, since we have an action of a gauge group, we must take account of this action in our framework.  This section therefore splits into three subsections.  In the first section, we study the gauge group and how it acts on the space of configurations.  Next, we study how this action decomposes the tangent space to the configuration space into natural subspaces.  Finally, we apply these decompositions to the study of the linearized Seiberg-Witten equations, where modulo gauge and other modifications, we can place ourselves in an elliptic situation.

\subsection{The Gauge Group}\label{SecGG}

The gauge group $\G = \G(Y) = \mathrm{Maps}(Y, S^1)$ \label{p:G} is the space of smooth maps $g: Y \to S^1$, where we regard $S^1 = \{e^{i\theta} \in \mathbb{C} : 0 \leq \theta < 2\pi\}$. Elements of the gauge group act on $\fC(Y)$ via
\begin{equation}
  (B,\Psi) \mapsto g^*(B,\Psi) = (B - g^{-1}dg, g\Psi). \label{gaction}
\end{equation}
It is straightforward to check that the Seiberg-Witten map $SW_3$ is gauge equivariant (where gauge transformations act trivially on $\Omega^1(Y;i\R)$).  In particular, the space of solutions to the Seiberg-Witten equations is gauge-invariant.

The gauge group decomposes into a variety of important subgroups, which will be important for the various kinds of gauge fixing we will be doing.  First, observe that $\pi_0(\G)$, the number of connected components of $\G$, satisfies
\begin{align}
  \pi_0(\G) & \cong H^1(Y; 2\pi i\Z). \label{pi0}
\end{align}
The correspondence (\ref{pi0}) is given by
\begin{equation}
  g \mapsto [g^{-1}dg], \label{du}
\end{equation}
where the latter denotes the cohomology class of the closed $1$-form $g^{-1}dg$.  Among subgroups of the gauge group, one usually considers the group of harmonic gauge transformations, i.e., gauge transformations such that $g^{-1}dg \in \ker d^*$.  However, on a manifold with boundary, $\ker (d + d^*)$ is infinite dimensional and we need to impose some boundary conditions.

On a manifold with boundary, Hodge theory tells us that we can make the following identifications between cohomology classes and harmonic forms with the appropriate boundary conditions\footnote{\label{footnotedf} For a differential form $a$ over a manifold $X$ with boundary, $a|_{\partial X}$ always denotes the differential form on $\partial X$ obtained via the restriction of those components tangential to $\partial X$.  Otherwise, given a section $u$ of a general vector bundle over $X$, $u|_X$ denotes the restriction of $u$ to the boundary, which therefore has values in the bundle restricted to the boundary. This clash of notation should not cause confusion since it will always be clear which restriction map we are using based on the context.}
\begin{align}
  H^1(Y; \R) &\cong \{\alpha \in \Omega^1(Y): da = d^*a = 0, *a|_{\Sigma} = 0\} \label{H1Y} \\
  H^1(Y,\Sigma; \R) &\cong \{\alpha \in \Omega^1(Y): da = d^*a = 0, a|_{\Sigma} = 0.\}. \label{H1YSigmaR}
\end{align}
In fact, we have two different Hodge decompositions, given by
\begin{align}
  \Omega^1(Y) &= \im d \oplus \im *d_n \oplus H^1(Y;\R)\\
  &= \im d_t \oplus \im *d \oplus H^1(Y,\Sigma;\R).
\end{align}
where
\begin{align}
  d_n: \{a \in \Omega^1(Y) : *a|_{\Sigma} = 0\} & \to \Omega^2(Y) \label{eq3.1:dn}\\
  d_t: \{\alpha \in \Omega^0(Y): \alpha|_\Sigma = 0\} & \to \Omega^1(Y). \label{eq3.1:dt}
\end{align}
Any gauge transformation $g$ in the identity component of the gauge group $\G_{\id}(Y)$ lifts to the universal cover of $S^1$ and so it can be expressed as $g = e^\xi$ for some $\xi \in \Omega^0(Y;i\R)$.  For such $g$, we have $g^{-1}dg = d\xi$, and thus we see that $\G/\G_{\id}$ is isomorphic to the integer lattice inside $\ker d/\im d$, which establishes the correspondence (\ref{pi0}).  Corresponding to the two cohomology groups (\ref{H1Y}) and (\ref{H1YSigmaR}), we can consider the following two subgroups of the harmonic gauge transformations
\begin{align}
  \G_{h,n}(Y) &= \{g \in \G : g^{-1}dg \in \ker d^*, *dg|_\Sigma = 0\} \label{Ghn}\\
  \G_{h,\partial}(Y) &= \{g \in \G : g^{-1}dg \in \ker d^*, g|_\Sigma = 1\}. \label{Ghp}
\end{align}
The group (\ref{Ghn}) is isomorphic to $S^1 \times H^1(Y;\Z)$, where the $S^1$ factor accounts for constant gauge transformations, and the group (\ref{Ghp}) is isomorphic to the integer lattice $H^1(Y,\Sigma;\Z)$ inside (\ref{H1YSigmaR}).

Next, we have the subgroup
$$\G_{\bot}(Y) = \{e^\xi \in \G_{\id} : \int_Y\xi = 0\}.$$
Thus, identifying  constant gauge transformations with $S^1$, we have the decompositions \label{p:Gpartial}
\begin{align*}
  \G_{\id}(Y) &= S^1 \times \G_\bot(Y)\\
  \G(Y) &= \G_{h,n}(Y) \times \G_\bot(Y).
\end{align*}
We have the following additional subgroups of the gauge group consisting of gauge transformations whose restriction to the boundary is the identity:
\begin{align}
  \G_{\partial}(Y) &= \{g \in \G(Y) : g|_\Sigma = 1\}\\
  \G_{\id,\partial}(Y) &= \G_{\id}(Y) \cap \G_{\partial}(Y)
\end{align}
Thus, we have
\begin{equation}
  \G_{\partial}(Y) = \G_{h,\partial}(Y)\times\G_{\id,\partial}(Y)
\end{equation}
and
\begin{equation}
  T_{\id}\G_{\id,\partial}(Y) = \{\xi \in \Omega^0(Y; i\R) : \xi|_\Sigma = 0\}.
\end{equation}

Since we consider the completion of our configuration spaces in Besov topologies, we must do so for the gauge groups as well.  Thus, let $\G^{s,p}(Y)$ denote the completion of $\G(Y)$ in $B^{s,p}(Y)$ and similarly for the other gauge groups.

\begin{Lemma}\label{Ggroup}
  For $s > 3/p$, the $B^{s,p}(Y)$ completions of $\G(Y)$ and its subgroups are Banach Lie groups.  If in addition $s > 1/2$, these groups act smoothly on $\fC^{s-1,p}(Y)$.
\end{Lemma}

\Proof For $s > 3/p$, the multiplication theorem, Theorem \ref{ThmMult}, implies $B^{s,p}(Y)$ is a Banach algebra.  Thus, $\G^{s,p}(Y)$ is closed under multiplication and has a smooth exponential map.  The second statement follows from (\ref{gaction}), Theorem \ref{ThmMult}, and the fact that $d: B^{s,p}(Y) \to B^{s-1,p}\Omega^1(Y)$ for all $s \in \R$ by Corollary \ref{CorDiffMap}.  Here the requirement $s > 1/2$ comes from the fact that we need $s + (s-1) > 0$ in Theorem \ref{ThmMult}.\End

Fix a smooth reference connection $B_\rf$.  From this, we obtain the Coulomb slice and Coulomb-Neumann slice through $B^{\rf}$, given by
\begin{align}
  \fC_C^{s,p}(Y) = &\{(B,\Psi) \in \fC^{s,p}(Y) : d^*(B-B^\rf)=0\}\\
  \fC_{C_n}^{s,p}(Y) = &\{(B,\Psi) \in \fC^{s,p} : d^*(B-B^\rf)=0, *(B-B^\rf)|_\Sigma=0\}, \qquad (s > 1/p)
\end{align}
respectively. The next lemma tells us that we can find gauge transformations which place any configuration into either of the above slices.

\begin{Lemma}\label{PiCdecomp}
  Let $s+1>\max(3/p,1/2)$. The action of the gauge group gives us the following decompositions of the configuration space:
  \begin{enumerate}
    \item We have\footnote{The direct products appearing in (\ref{Cslice}) and (\ref{Cnslice}) mean that the gauge group factor acts freely on the subspace appearing in the second factor so that the space on the left is equal to the resulting orbit space obtained from the right-hand side.}
  \begin{align}
    \fC^{s,p}(Y) &= \G_{\id,\partial}^{s+1,p}(Y) \times \fC_C^{s,p}(Y). \label{Cslice}
  \end{align}
  \item Suppose in addition $s > 1/p$.  Then we have
  \begin{align}
    \fC^{s,p}(Y) &= \G_{\bot}^{s+1,p}(Y) \times \fC_{C_n}^{s,p}(Y). \label{Cnslice}
  \end{align}
  \end{enumerate}
\end{Lemma}

\Proof (i) Since $s + 1 > \max(3/p,1/2)$, the previous lemma implies $\G^{s+1}_{\id,\partial}(Y)$ is a Banach Lie group and it acts on $\fC^{s,p}(Y)$. If $u = e^\xi \in G^{s+1,p}_{\id,\partial}$ puts a configuration $(B_\rf + b,\Psi)$ into the Coulomb slice through $B_\rf$, then $\xi$ satisfies
\begin{equation}
\left\{\begin{array}{rl}
  \Delta \xi &= d^*b \in B^{s-1,p}(Y;i\R),\\
  \xi|_\Sigma &= 0.
\end{array}\right.  \label{DL}
\end{equation}
The Dirichlet Laplacian is an elliptic boundary value problem and since $s+1 > 1/p$, we may apply Corollary \ref{CorEBP}, which shows that we have an elliptic estimate
$$\|\xi\|_{B^{s+1,p}} \leq C(\|\Delta \xi\|_{B^{s,p}} + \|\xi\|_{B^{s,p}})$$
for $\xi$ satisfying (\ref{DL}).  A standard computation shows that the kernel and cokernel of the Dirichlet Laplacian is zero, and so we have existence and uniqueness for the Dirichlet problem.  This implies the decomposition.  

(ii) The analysis is the same, only now we have a homogeneous Neumann Laplacian problem for $\xi$:
\begin{equation}
\left\{\begin{array}{rl}
  \Delta \xi &= d^*b \in B^{s-1,p}(Y; i\R)\\
  *d\xi|_\Sigma &= *b|_\Sigma \in B^{s-1/p,p}\Omega^2(\Sigma; i\R).
\end{array}\right.  \label{NL}
\end{equation}
Since the Neumann Laplacian is an elliptic boundary value problem, we can apply Corollary \ref{CorEBP} again.  The inhomogeneous Neumann problem $\Delta \xi = f$ and $\partial_\nu \xi = g$ has a solution if and only if $\int_Y f + \int_\Sigma g = 0$, and this solution is unique up to constant functions.  Since we always have $\int_Y d^*b + \int_\Sigma *b = 0$, then
(\ref{NL}) has a unique solution $\xi \in B^{s+1,p}(Y)$ subject to $\int_Y \xi = 0$.  The decomposition now follows.\End

In light of Lemma \ref{PiCdecomp}, we can regard the quotient of $\fC^{s,p}(Y)$ by the gauge groups $\G^{s+1,p}_{\id,\partial}$ and $\G^{s+1,p}_{\bot}$ as subspaces of $\fC^{s,p}(Y)$, namely, those configurations in Coulomb and Coulomb-Neumann gauge with respect to $B_\rf$.

\begin{Rem}
  In gauge theory, one usually also considers the quotient of the configuration space by the entire gauge group.  In our case (which is typical) the quotient space is singular since different elements of the configuration space have different stabilizers.  Namely, if $(B,\Psi) \in \fC(Y)$ is such that $\Psi \not\equiv 0$, then it has trivial stabilizer, whereas if $\Psi \equiv 0$, then it has stabilizer $S^1$, the constant gauge transformations.  In the former case, such a configuration is said to be \textit{irreducible}, otherwise it is \textit{reducible}.  We will not need to consider the quotient space by the entire gauge group in this paper, and we will only need to consider the decompositions in Lemma \ref{PiCdecomp}.
\end{Rem}

\subsection{Decompositions of the Tangent Space}\label{SecSpaces}

The action of the gauge group on the configuration space induces a decomposition of the tangent space to a configuration $\c$ into the subspace tangent to the gauge orbit through $\c$ and its orthogonal complement.  More precisely, let \label{p:Tc}
\begin{equation}
  \T_\c := T_\c\fC(Y) = \Omega^1(Y; i\R) \oplus \Gamma(\S)
\end{equation}
be the smooth tangent space to a smooth configuration $\c$.  Define the operator \label{p:bd}
\begin{align}
\dd_\c: \Omega^0(Y; i\R) &\to \T_\c \nonumber\\
\xi &\mapsto (-d\xi, \xi\Psi), \label{bd}
\end{align}
and let \label{p:J}
\begin{equation}
\J_\c := \im \bd_\c \subset \T_\c
\end{equation}
be its image.  Then observe that $\J_\c$ is the tangent space to the gauge orbit at $\c$.  Indeed, this follows from differentiating the action (\ref{gaction}) at the identity. We also have the adjoint operator \label{p:bd^*}
\begin{align}
  \dd_\c^*: \T_\c  &\to \Omega^0(Y; i\R) \nonumber\\
(b,\psi)  &\mapsto -d^*b + i\Re(i\Psi,\psi), \label{bd*}
\end{align}
and we define the subspace \label{p:K}
\begin{equation}
  \K_\c := \ker \bd^*_\c \subset \T_\c.
\end{equation}
On a closed manifold, $\K_\c$ is the $L^2$ orthogonal complement of $\im \bd_\c$.  In this case, the orthogonal decomposition of $\T_\c$ into the spaces $\J_\c$ and $\K_\c$ plays a fundamental role in the analysis of \cite{KM}.  In our case, since we have a boundary, we will impose various boundary conditions on these spaces, and the resulting spaces will play a very important role for us too.  Moreover, we will take the appropriate Besov completions of these spaces.

Thus, let $\c \in \fC^{t,q}(Y)$ be any configuration of regularity $B^{t,q}(Y)$, where $t \in \R$ and $q \geq 2$.  For $s \in \R$ and $p \geq 2$, let
\begin{equation}
  \T_\c^{s,p} := B^{s,p}(\Omega^1(Y; i\R) \oplus \Gamma(\S))
\end{equation}
be the Besov closure of $\T_\c$.  It is independent of $\c$ and is equal to the tangent space $T_\c\fC^{s,p}(Y)$ when $(s,p)=(t,q)$.    So long as we have bounded multiplication maps
\begin{align}
  B^{t,q}(Y) \times B^{s+1,p}(Y) \to B^{s,p}(Y) \label{sec3:Bmult1} \\
  B^{t,q}(Y) \times B^{s,p}(Y) \to B^{s-1,p}(Y),\label{sec3:Bmult2}
\end{align}
then we can define maps
\begin{align*}
\dd_\c: B^{s+1,p}\Omega^0(Y; i\R) &\to \T^{s,p}_\c\\
\xi &\mapsto (-d\xi, \xi\Psi),\\
\dd_\c^*: \T_\c^{s,p}  &\to B^{s-1,p}\Omega^0(Y; i\R) \\
(b,\psi)  &\mapsto -d^*b + i\Re(i\Psi,\psi),
\end{align*}
respectively.  In particular, if $(t,q) = (s,p)$ and $s > 3/p$, then by Theorem \ref{ThmMult}, the multiplications (\ref{sec3:Bmult1}) and (\ref{sec3:Bmult2}) are bounded.

Thus, when (\ref{sec3:Bmult1}) and (\ref{sec3:Bmult2}) hold, define the following subspaces of $\T^{s,p}_\c$:
\begin{align}
  \J^{s,p}_\c  &= \im \left(\dd_\c: B^{s+1,p}\Omega^0(Y; i\R) \to \T^{s,p}_\c\right)\\
  \J^{s,p}_{\c, \bot} &= \{(-d\xi,\xi\Psi) \in \J^{s,p}_\c : \int_Y\xi = 0\}\\
  \J^{s,p}_{\c, t} &= \{(-d\xi,\xi\Psi) \in \J^{s,p}_\c : \xi|_\Sigma = 0\}\\
  \K^{s,p}_\c  &= \ker \left(\dd_\c^*: \T^{s,p}_\c \to B^{s-1,p}\Omega^0(Y;i\R)\right)\\
  \K^{s,p}_{\c, n} &= \{(b,\psi) \in \K^{s,p}_\c : *b|_\Sigma = 0\}.
\end{align}
Observe that when $\c \in \fC^{s,p}(Y)$, then $\J^{s,p}_{\c}$, $\J^{s,p}_{\bot}$, $\J^{s,p}_{\c,t}$ are the tangent spaces to the gauge orbit of $\c$ in $\fC^{s,p}$ determined by the gauge groups $\G^{s+1,p}(Y)$, $\G^{s+1,p}_{\bot}(Y)$, and $\G^{s+1,p}_{\partial}(Y)$, respectively.  Note that the subscript $t$ appearing in $\J^{s,p}_{\c, t}$ is a label to denote that the (tangential) restriction of $\xi$ to the boundary vanishes; it is not to be confused with a real parameter.  This is consistent with the notation used in (\ref{eq3.1:dt}). Likewise, the subscript $n$ appearing in $\K^{s,p}_{\c, n}$ and (\ref{eq3.1:dn}) denotes that the elements belonging to these spaces have normal components for their $1$-form parts equal to zero on the boundary. We also have the linear Coulomb and Coulomb-Neumann slices:
\begin{align}
  \C^{s,p}_\c & = \{(b,\psi) \in \T_\c^{s,p} : d^*b = 0\}\\
  \C^{s,p}_{\c,n} & = \{(b,\psi) \in \T_\c^{s,p} : d^*b = 0, *b|_\Sigma = 0\}.
\end{align}

The following lemma is essentially the linear version of Lemma \ref{PiCdecomp}.  The statement is only mildly more technical in that one may consider the basepoint $(B,\Psi)$ and the tangent space $\T_{(B,\Psi)}$ in different topologies.  We do this because we will need to consider topologies on $\T_{(B,\Psi)}$ that are weaker than the regularity of $(B,\Psi)$, which occurs, for example, when we apply differential operators to elements of $\T_{(B,\Psi)}^{s,p}$ when $\c \in \fC^{s,p}(Y)$, thereby obtaining spaces such as $\T_\c^{s-1,p}$.  These spaces and their decompositions will become important for us in the next section, when we study the linearized Seiberg-Witten equations and try to recast them in a form in which they become elliptic.

\begin{Lemma}\label{LemmaTdecomp}
Let $s+1 > 1/p$ and let $(B,\Psi) \in \fC^{t,q}(Y)$, where $t > 3/q$, $q \geq 2$ are such that (\ref{sec3:Bmult1}) and (\ref{sec3:Bmult2}) hold.  In particular, if $q = p$, then we need $t \geq s$ and $t > \max(-s,3/p)$.
\begin{enumerate}
  \item We have the following decompositions:
\begin{align}
\T_\c^{s,p} &= \J^{s,p}_{\c,t} \oplus \K^{s,p}_{\c} \label{JtK}\\
\T_\c^{s,p} &= \J^{s,p}_{\c,t} \oplus \C^{s,p}_{\c}. \label{JtC}
\end{align}
 \item If in addition $s > 1/p$, then
 \begin{align}
\T_\c^{s,p} &= \J^{s,p}_{\c,\bot} \oplus \C^{s,p}_{\c,n} \label{Cn}.
\end{align}
If $\Psi \not\equiv 0$, then furthermore
\begin{align}
\T_\c^{s,p} &= \J^{s,p}_\c \oplus \K^{s,p}_{\c,n} \label{Kn}.
\end{align}
\end{enumerate}
\end{Lemma}

\Proof We first prove (\ref{JtK}).  Given $(b,\psi) \in \T^{s,p}_\c$, consider the boundary value problem
\begin{equation}
\left\{\begin{array}{rl}
  \Delta_\c \xi &= f \in B^{s-1,p}(Y; i\R)\\
  \xi|_\Sigma &= 0,
\end{array}\right.  \label{DL2}
\end{equation}
where $f = \bd_\c^*b$ and
\begin{equation}
  \Delta_\c := \bd_\c^*\bd_\c = \Delta + |\Psi|^2. \label{Deltac}
\end{equation}
We have $\bd^*_\c b\in B^{s-1,p}(Y; i\R)$ since we have a bounded multiplication $B^{t,p}(Y)\times B^{s,p}(Y) \to B^{s-1,p}(Y)$ by the hypotheses.  Likewise, since we have a bounded map $B^{t,p}(Y) \times B^{s+1,p}(Y) \to B^{s,p}(Y)$, we see that multiplication by $|\Psi|^2 \in B^{t,p}(Y)$ is a compact perturbation of $\Delta: B^{s+1,p}(Y) \to B^{s-1,p}(Y)$.  Thus, the Dirichlet boundary value problem (\ref{DL2}) is Fredholm for $s+1 > 1/p$ (where the requirement on $s$ is so that Dirichlet boundary conditions make sense, cf. Corollary \ref{CorEBP}).  Moreover, since $|\Psi|^2$ is a positive multiplication operator, a simple computation shows the existence and uniqueness of (\ref{DL2}).  Indeed, if $\Delta \alpha = -|\Psi|^2\alpha$ and $\alpha|_\Sigma = 0$, then repeated elliptic boostrapping for the inhomogeneous Dirichlet Laplacian shows that $\alpha \in B^{t+2,q}(Y) \subseteq B^{2,2}(Y)$ since $t > 0$ and $q \geq 2$.  Then
$$0 = \int_Y(\Delta_\c \alpha, \alpha) = \|\nabla \alpha\|^2_{L^2(Y)} + \|\Psi \alpha\|_{L^2(Y)}^2,$$
which implies $\alpha$ is constant. Hence, $\alpha = 0$ since $\alpha|_\Sigma = 0$.  Thus, (\ref{DL2}) has no kernel and since the adjoint problem of (\ref{DL2}) is itself, we see that (\ref{DL2}) has no cokernel as well.  Thus, the existence and uniqueness of (\ref{DL2}) is established.  Let $\Delta_{\c,t}^{-1}$ denote the solution map of (\ref{DL2}).  We have shown that $\Delta_{\c,t}^{-1}: B^{s-1,p}(Y) \to B^{s+1,p}(Y)$ is bounded.  The projection onto $\J^{s,p}_{\c,t}$ through $\K^{s,p}_{\c}$ is now seen to be given by
\begin{equation}
  \Pi_{\J^{s,p}_{\c,t}} = \bd_\c \Delta_{\c,t}^{-1}\bd_\c^* \label{PiJK}
\end{equation}
and it is bounded on $\T^{s,p}_\c$ since $\bd_\c: B^{s+1,p}\Omega^0(Y;i\R) \to B^{s,p}\Omega^1(Y;i\R)$ is bounded. This gives us the decomposition (\ref{JtK}).  Similarly, we get the decomposition (\ref{JtC}) if we replace $\Delta_\c$ with $\Delta$ in the above.

For (ii), if we consider the inhomogeneous Neumann problem for $\Delta_\c$ instead of the Dirichlet problem, proceeding as above yields (\ref{Kn}), since when $\Psi \not\equiv 0$, a similar computation shows that we get existence and uniqueness.  Here, we need $s > 1/p$ so that $s + 1 > 1+1/p$ and the relevant Neumann boundary condition makes sense.  Similarly, considering the inhomogeneous Neumann problem for $\Delta$ yields (\ref{Cn}).\End

For any $s, t \in \R$, we can define the Banach bundle
\begin{equation}
  \T^{s,p}(Y) \to \fC^{t,p}(Y) \label{Tbundle}
\end{equation}
whose fiber over every $\c \in \fC^{t,p}(Y)$ is the Banach space $\T^{s,p}_\c$.  Of course, all the $\T^{s,p}_\c$ are identical, so the bundle (\ref{Tbundle}) is trivial.  If $s = t$, then (\ref{Tbundle}) is the tangent bundle of $\fC^{t,p}(Y)$.  If $s,t$ satisfy the hypotheses of the previous lemma, decomposing each fiber $\T^{s,p}_\c$ according to the decomposition (\ref{JtK}) defines us Banach subbundles of $\T^{s,p}(Y)$.  This is the content of the below proposition, where we specialize to a range of parameters relevant to the situations we will encounter later, e.g., see Lemma \ref{LemmaHtrans}.

\begin{Proposition}\label{PropBundle}
  Let $s > 3/p$.  If $\max(-s,-1+1/p) < s' \leq s$, then the Banach bundles \label{p:K(Y)}
  \begin{align*}
  \J^{s',p}_t(Y) &\to \fC^{s,p}(Y)\\
  \K^{s',p}(Y) & \to \fC^{s,p}(Y),
  \end{align*}
  whose fibers over $\c \in \fC^{s,p}(Y)$ are $\J^{s',p}_{\c,t}$ and $\K^{s',p}_\c$, respectively, are complementary subbundles of $\T^{s',p}(Y)$.
\end{Proposition}

\Proof The restrictions on $s$ and $s'$ ensure that we can apply Lemma \ref{LemmaTdecomp}.  From this, one has to check that the resulting decomposition
$$\T^{s',p'}_\c = \J^{s',p'}_{\c,t} \oplus \K^{s',p'}_\c$$
varies continuously with $\c \in \fC^{s,p}(Y)$.  For this, it suffices to show that the projection $\Pi_{\J^{s',p}_{\c,t}}$ given by (\ref{PiJK}), with range $\J^{s',p}_{\c,t}$ and kernel $\K^{s',p}_\c$, varies continuously with $\c \in \fC^{s,p}(Y)$.  Once we prove that $\J^{s',p}_t(Y)$ is a subbundle, it automatically follows that $\K^{s',p}(Y)$ is a (complementary) subbundle, since then the complementary projection \label{p:PiK}
\begin{equation}
  \Pi_{\K^{s',p}_{\c}} = 1 - \Pi_{\J^{s',p}_{\c,t}} \label{3.1-PiK}
\end{equation}
onto $\K^{s',p}_\c$ varies continuously with $\c$.

From the multiplication theorem, Theorem \ref{ThmMult}, since
$$\Delta_{\c,t}: \{\xi \in B^{s'+1,p}(Y; i\R) : \xi|_\Sigma = 0\} \to B^{s'-1,p}(Y; i\R)$$
varies continuously with $\c \in B^{s,p}(Y)$ and is an isomorphism for all $\c$, its inverse $\Delta_{\c,t}^{-1}$ also varies continuously. Likewise, $\bd_\c^*: \T^{s',p}_\c \to B^{s'-1,p}\Omega^0(Y; i\R)$ and $\bd_\c: B^{s'+1,p}\Omega^0(Y, i\R) \to \T^{s',p}_\c$ vary continuously with $\c \in \fC^{s,p}(Y)$. This establishes the required continuity of $\Pi_{\J^{s',p}_\c} = \bd_\c\Delta_{\c,t}^{-1}\bd_\c^*$ with respect to $\c$.\End

The Banach bundle $\K^{s',p}(Y)$, with $s' = s-1$ will be used to establish transversality properties of the Seiberg-Witten map $SW_3$, see Theorem \ref{ThmMMan}. 

\subsection{The Linearized Seiberg-Witten Equations}\label{SecLinSW}

In this section, we study the linearization of the Seiberg-Witten map $SW_3$ to prove basic properties concerning the (formal) tangent space to our monopole spaces on $Y$ and their behavior under restriction to the boundary.  If the linearization of the Seiberg-Witten equations were elliptic, this would be quite straightforward from the analysis of elliptic boundary value problems, the relevant results of which are summarized in the appendix.  However, because the Seiberg-Witten equations are gauge-invariant, its linearization is not elliptic and we have to do some finessing to account for the gauge-invariance.  To do this, we make fundamental use of the subspaces and decompositions of the previous section.

Before we get started, let us note that our main theorem of this section, Theorem \ref{ThmLinLag}, proves a bit more than what is needed to prove our main theorems.  Indeed, it is mostly phrased in such a way that the results of this section can be tied into the general framework of the pseudodifferential analysis of elliptic boundary value problems in Appendix \ref{AppSeeley} (see the discussion preceding Theorem \ref{ThmLinLag}). Moreover, some of the consequences of Theorem \ref{ThmLinLag} will only be put to full use in \cite{N2}.  Thus, the reader should regard this section as a general framework for studying the Hessian and augmented Hessian operators, (\ref{2-Hess}) and (\ref{aHessian}), whose kernels are equal to the tangent spaces to $\fM$ and $\M$, respectively, via (\ref{eq3.3:TfM}) and (\ref{eq3.3:TM}). Much of this framework consists in the construction of pseudodifferential type operators associated to the Hessian and augmented  operators, namely the Calderon projection and Poisson operators, see Lemma \ref{LemmatHCP} and Definition \ref{DefThm}.  For the augmented Hessian, an elliptic operator, these operators are defined as in Definition \ref{DefCP}, and for the non-elliptic Hessian, they are defined by analogy in Definition \ref{DefThm}.  In a few words, the significance of these operators is that they relate the kernel of the (augmented) Hessian with the kernel's boundary values in a simple and uniform way across multiple topologies.  This is what allows us to relate the tangent spaces to $\fM$ and $\M$ with the tangent spaces to $\L$, the latter being the boundary values of the kernels of the Hessian operators via (\ref{eq3.3:TL}).  Unfortunately, the infinite dimensional nature of all spaces involved and the presence of multiple topologies makes the work we do quite technical.  As a suggestion to the reader, it would be best to first absorb the main ideas of Appendix \ref{AppSeeley} and to understand the statements of Lemma \ref{LemmatHCP} and Theorem \ref{ThmLinLag} before plunging into the details.\\

Let \label{p:T}
\begin{align}
  \T &= \Omega^1(Y; i\R) \oplus \Gamma(\S)
\end{align}
be a fixed copy of the tangent space $\T_\c = T_\c\fC(Y)$ to any smooth configuration $\c \in \fC(Y)$.\footnote{There is no real distinction between $\T$ and a particular tangent space $\T_{\c}$ to a configuration, since $\fC(Y)$ is an affine space.  However, when we study the spaces $\fM$ and $\M$ as subsets of $\fC(Y)$ in Section 4, we will reintroduce base points when we have a particular tangent space in mind. For now, we drop basepoints to minimize notation.}  Thus, all the subspaces of $\T_\c$, namely $\J_{\c}$, $\K_{\c}$, and their associated subspaces defined in the previous section, may be regarded as subspaces of $\T$ that depend on a configuration $\c \in \fC(Y)$.  We let \label{p:C}
\begin{equation}
  \cC = \{(b,\psi) \in \T: d^*b = 0\}
\end{equation}
denote the Coulomb-slice in $\T$. Likewise, let \label{p:TSigma}
\begin{align}
  \T_\Sigma &= \Omega^1(\Sigma; i\R) \oplus \Gamma(\S_\Sigma)
\end{align}
denote a fixed copy of the tangent space to any smooth configuration of $\fC(\Sigma)$. The restriction map (\ref{rSigma}) on configuration spaces induces a restriction map on the tangent spaces \label{p:rSigma2}
\begin{align}
  r_\Sigma: \T & \to \T_\Sigma \nonumber \\
  (b,\psi) & \mapsto (b|_\Sigma,\psi|_\Sigma).
\end{align}

From (\ref{SW3map}), the linearization of the Seiberg-Witten map $SW_3$ at a configuration $\c \in \fC(Y)$ yields an operator \label{p:Hess}
\begin{align}
  \H_\c: \T & \to \T \nonumber \\
  \H_{\c} &= \begin{pmatrix}*d & 2i\Im \rho^{-1}(\cdot\Psi^*)_0\\ \rho(\cdot)\Psi & D_B \end{pmatrix} \label{2-Hess}
\end{align}
which acts on the tangent space $\T$ to $\c$.  We call the operator $\H_\c$ the \textit{Hessian}.\footnote{On a closed manifold, $\H_\c$ would in fact be the Hessian of the Chern-Simons-Dirac functional, see \cite{KM}.} The Hessian operator is a formally self-adjoint first order operator.  For any monopole $(B,\Psi) \in SW_3^{-1}(0)$, we (formally) have that the tangent spaces to our monopole spaces $\fM$ and $\L$ are given by
\begin{align}
  T_\c\fM &= \ker \H_\c \label{eq3.3:TfM}\\
  T_{r_\Sigma\c}\L &= r_\Sigma(\ker\H_\c). \label{eq3.3:TL}
\end{align}
Indeed, this is just the linearization of (\ref{smoothfM}) and (\ref{sec2:L=rM}). Thus, understanding $\fM$ and $\L$ at the linear level is the same as understanding the kernel of $\H_\c$.

Unfortunately, $\H_\c$ is not elliptic, which follows from a simple examination of its symbol. In fact, this nonellipticity follows a priori from the equivariance of the Seiberg-Witten map under gauge transformations.  In particular, since the zero set of $SW_3$ is gauge-invariant, then the linearization $\H_\c$ at a monopole $\c$ annihilates the entire tangent space to the gauge orbit at $\c$, i.e., the subspace $\J_\c \subset \T$.  Furthermore, even if we were to account for this gauge invariance by say, placing configurations in Coulomb-gauge, i.e., if we were instead to consider the operator $\H_\c \oplus d^*: \T \to \T \oplus \Omega^0(Y;i\R)$, we still would not have an elliptic operator in the usual sense.

However, there is a simple remedy for this predicament. Following \cite{KM}, the operator $\H_\c$ naturally embeds as a summand of an elliptic operator.  Namely, if we enlarge the space $\T$ to the \textit{augmented tangent space} \label{p:tT}
\begin{equation}
  \tT := \T \oplus \Omega^0(Y;i\R),
\end{equation}
then we can consider the \textit{augmented Hessian}\footnote{In \cite{KM}, the operators $\bd_\c$ and $\bd_\c^*$ are used in the definition of $\tH_\c$ instead of $-d$ and $-d^*$, respectively.  Our definition reflects the fact that we will work with Coulomb slices $\C_\c$ instead of the slices $\K_\c$ inside $\T$.  The presence of the minus signs on $-d$ and $-d^*$ in $\tH_\c$ lies in the relationship between $\tH_\c$ and the linearization of the $4$-dimensional Seiberg-Witten equations, see \cite{N2}.  Thus, the augmented Hessian operator is not an ad-hoc extension of the Hessian operator but is tied to the underlying geometry of the problem.}\label{p:aHess}
\begin{align}
  \tH_\c: \tT &\to \tT \nonumber \\
  \tH_\c &= \begin{pmatrix}
    \H_\c & -d\\
    -d^* & 0
  \end{pmatrix}. \label{aHessian}
\end{align}
The augmented Hessian is a formally self-adjoint first order elliptic operator, as one can easily verify.  This operator takes into account Coulomb gauge-fixing via the operator $d^*: \Omega^1(Y;i\R) \to \Omega^0(Y;i\R)$, while ensuring ellipticity by adding in the adjoint operator $d: \Omega^0(Y;i\R) \to \Omega^1(Y;i\R)$.  The advantage of studying the operator $\tH_\c$ is that we may apply the pseudodifferential tools from Appendix \ref{AppSeeley} to understand the kernel of $\tH_\c$ and its boundary values.  Moreover, we have (formally) that
\begin{equation}
  T_{\c}\M = \ker (\tH_\c|_{\T}). \label{eq3.3:TM}
\end{equation}

The space of boundary values for $\tT$ is the space \label{p:tTSigma}
\begin{equation}
  \tT_\Sigma := \T_\Sigma \oplus \Omega^0(\Sigma; i\R) \oplus \Omega^0(\Sigma; i\R). \label{tT_Sigma}
\end{equation}
Indeed, one can see that $\tT|_\Sigma \cong \tT_\Sigma$ via the full restriction map $r: \tT \to \tT_\Sigma$ given by \label{p:r}
\begin{eqnarray}
r : \Omega^1(Y; i\R) \oplus \Gamma(\S) \oplus \Omega^0(Y;i\R) & \to &\Omega^1(\Sigma; i\R) \oplus \Gamma(\S_\Sigma) \oplus \Omega^0(\Sigma; i\R) \oplus \Omega^0(\Sigma; i\R) \nonumber \\
(b,\psi,\alpha) & \mapsto & (b|_\Sigma, \psi|_\Sigma, -b(\nu), \alpha|_\Sigma), \label{rmap}
\end{eqnarray}
where in (\ref{rmap}), the term $b(\nu)$ denotes contraction of the $1$-form $b$ with the outward normal \label{p:nu} $\nu$ to $\Sigma$.  Thus, the two copies of $\Omega^0(\Sigma;i\R)$ in $\tT_\Sigma$ are meant to capture the normal component of $\Omega^1(Y;i\R)$ and the trace of $\Omega^0(Y;i\R)$ along boundary.  The map $r_\Sigma: \T \to \T_\Sigma$ appears as the first factor of the map $r$, and it is the tangential part of the full restriction map.  Since we can regard $\T \subset \tT$, then by restriction, the map $r$ also maps $\T$ to $\tT_\Sigma$.

As usual, we can consider the Besov completions of all the spaces involved.  Thus, we have the spaces
$$\C^{s,p}, \T^{s,p},\; \tT^{s,p},\;  \T_\Sigma^{s,p},\; \tT^{s,p}_\Sigma$$
which we use to denote the $B^{s,p}$ completions of their corresponding smooth counterparts.  The restriction maps $r_\Sigma$ and $r$ extend to Besov completions in the usual way.  We also have the spaces $\J^{s,p}_\c$, $\K^{s,p}_\c$, and their subspaces from the previous section, which we may all regard as subspaces of $\T^{s,p}$.

The plan for the rest of this section is as follows.  First, we investigate the kernel of the elliptic operator $\tH_\c$.  We do this first for smooth $\c$, in which case the tools from Appendix \ref{AppSeeley} apply, and then we consider nonsmooth $\c$, in which case modifications must be made.  Here, one has to keep track of the function space arithmetic rather carefully.  Next, we will relate the kernel of $\tH_\c$ to the kernel of $\H_\c$ and see how these spaces behave under the restriction maps $r$ and $r_\Sigma$, respectively.  For this, we place these results under the conceptual framework of Appendix \ref{AppSeeley} by way of using the Calderon projection and Poisson operator associated to an elliptic operator.  For the Hessian $\H_\c$, the main technical issue here is its non-ellipticity (i.e. gauge-invariance). The results of our analysis are summarized in the main theorem of this section, Theorem \ref{ThmLinLag}.\footnote{The complexity of the function space arithmetic in this section can minimized if one does not care about the symplectic properties of the spaces involved, namely, the Lagrangian properties in Lemma \ref{LemmaHessCD} and Theorem \ref{ThmLinLag}(i). See Remark \ref{Rem1/2}.}

We organize the preceding notation into the below diagram, since it will be used consistently for the rest of this paper:

\begin{align} 
\begin{split}
  \xy
(0,0)*+{\tT}="1";
(50,0)*+{\tT=\T\oplus\Omega^0(Y;i\R)\hspace{.4cm}}="2";
(0,-20)*+{\T}="3";
(50,-20)*+{\T=\Omega^1(Y;i\R)\oplus\Gamma(\S)}="4";
{\ar^(.35){\tH_\c} "1";"2"};
{\ar^(.35){\H_\c} "3";"4"};
{\ar@{^{(}->} (0,-16);"1"};    
{\ar@{^{(}->} (35,-16);(35,-2.5)}; 
(0,-40)*+{\tT}="5";
(58,-40)*+{\tT_\Sigma=\T_\Sigma\oplus\Omega^0(\Sigma;i\R)\oplus\Omega^0(\Sigma;i\R)}="6";
(0,-60)*+{\T}="7";
(58,-60)*+{\T_\Sigma=\Omega^1(\Sigma;i\R)\oplus\Gamma(\S_\Sigma)\hspace{1.3cm}}="8";
{\ar^(.3){r} "5";"6"};
{\ar^(.3){r_\Sigma} "7";"8"};
{\ar@{^{(}->} (0,-56);"5"}; 
{\ar@{^{(}->} (35,-56);(35,-42.5)}; 
\endxy \\
\end{split}\label{CD0} \\ \nonumber
\end{align}

In studying the augmented Hessian operators $\tH_\c$ for smooth $\c \in \fC(Y)$, observe that they all differ by bounded zeroth order operators. Indeed, if we write $(b,\psi) = (B_1,\Psi_1) - (B_0,\Psi_0)$, then
$$\tH_{(B_1,\Psi_1)} - \tH_{(B_0,\Psi_0)} = (b,\psi)\#$$
where $(b,\psi)\#$ is the multiplication operator given by
\begin{align}
  (b,\psi)\#: \T &\to \T \nonumber\\
  (b',\psi') & \mapsto (2i\Im \rho^{-1}(\psi'\psi^*)_0, \rho(b)\psi'). \label{sharp}
\end{align}
In general, we will use $\#$ to denote any kind of pointwise multiplication map.  Let $B_\rf$ be our fixed smooth reference connection.  Define
\begin{align}
  \tH_0 & := \tH_{(B_\rf,0)} \nonumber \\
        & = D_{\dgc} \oplus D_{B_\rf} \label{deftH0}
\end{align}
where $D_{B_\rf}: \Gamma(\S) \to \Gamma(\S)$ is the Dirac operator on spinors determined by $B_\rf$ and $D_{\dgc}$ is the div-grad curl operator
\begin{equation}
  D_{\dgc} := \begin{pmatrix}*d & - d^* \\ -d & 0\end{pmatrix}: \Omega^1(Y; i\R) \oplus \Omega^0(Y; i\R) \circlearrowleft \label{DGCop}
\end{equation}
The operator $D_{\dgc}$ is also a Dirac operator.  Thus, the operator $\tH_0$ is a Dirac operator and every other $\tH_\c$ is a zeroth order perturbation of $\tH_0$.  Our first objective therefore is to understand the operator $\tH_0$.

Let us quickly review some basic properties about general Dirac operators (a more detailed treatment can be found in \cite{BBW}).  Let $\mathsf{D}$ be any Dirac operator acting on sections $\Gamma(E)$ of a Clifford bundle $E$ over $Y$ endowed with a connection compatible with the Clifford multiplication.  Here, by a Dirac operator, we mean any operator equal to ``the" Dirac operator on $E$ (the operator determined by the Clifford multiplication and compatible connection) plus any zeroth order symmetric operator. Let $(\cdot,\cdot)$ denote the (real or Hermitian) inner product on $E$. Working in a collar neighborhood of $[0,\eps] \times \Sigma$ of the boundary, where $t \in [0,\eps]$ is the inward normal coordinate, we can identify $\Gamma(E|_{[0,\eps] \times \Sigma})$ with $\Gamma([0,\eps], \Gamma(E_{\Sigma}))$, the space of $t$-dependent sections with values in $\Gamma(E_\Sigma)$.  Under this idenfication, we can write any Dirac operator $D$ as
\begin{equation}
  \mathsf{D} = \mathsf{J}_t\left(\frac{d}{dt} + \mathsf{B}_t + \mathsf{C}_t\right), \label{Diractan}
\end{equation}
where $\mathsf{J}_t$, $\mathsf{B}_t$, and $\mathsf{C}_t$ are $t$-dependent operators acting on $\Gamma(E_\Sigma)$.  The operator $\mathsf{J}_t$ (which is Clifford multiplication by $d/dt$) is a skew-symmetric bundle automorphism satisfying $\mathsf{J}_t^2 = -\id$, the operator $\mathsf{B}_t$ is a first order elliptic self-adjoint operator, and $\mathsf{C}_t$ is a zeroth order bundle endomorphism.

\begin{Def}\label{DefTBO}
  We call $\mathsf{B}_0: \Gamma(E_\Sigma) \to \Gamma(E_\Sigma)$ the \textit{tangential boundary operator} associated to $D$.
\end{Def}

Observe that the above definition is only well-defined up to a symmetric zeroth order term.  By abuse of terminology, we may also refer to the family of operators $\mathsf{B}_t$ in (\ref{Diractan}) as tangential boundary operators as well.

The significance of the decomposition (\ref{Diractan}) is that the space of boundary values of the kernel of $\mathsf{D}$ is, up to a compact error, determined by the operator $\mathsf{B}_0$.  More precisely, we have the following picture. Since $\mathsf{B_0}$ is a first order self-adjoint elliptic operator, the space $\Gamma(E_\Sigma)$ decomposes as
\begin{equation}
  \Gamma(E_\Sigma) = \cZ^+_{\mathsf{B}_0} \oplus \cZ^-_{\mathsf{B}_0} \oplus \cZ^0_{\mathsf{B}_0},
\end{equation}
the positive, negative, and zero spectral subspaces of $\mathsf{B}_0$, respectively.  Moreover, since the projections onto these subspaces are given by pseudodifferential operators, we get a corresponding decomposition on the Besov space completion:
\begin{equation}
  B^{s,p}(E_\Sigma) = B^{s,p}\cZ^+_{\mathsf{B}_0} \oplus B^{s,p}\cZ^-_{\mathsf{B}_0} \oplus B^{s,p}\cZ^0_{\mathsf{B}_0},
\end{equation}
for all $s \in \R$ and $1<p<\infty$.  If we let $\mathsf{D}: B^{s,p}(E) \to B^{s-1,p}(E)$, then we can consider the boundary values of its kernel $r(\ker \mathsf{D}) \subset B^{s-1/p,p}(E_\Sigma)$.  Then what we have is that the spaces $r(\ker \mathsf{D})$ and $B^{s-1/p,p}\cZ_{\mathsf{B}_0}^+$ are \textit{commensurate}, that is, they differ by a compact perturbation\footnote{More precisely, the range of $r(\ker \mathsf{D})$ and $\cZ^+_{\mathsf{B}_0}$ are each given by the range of pseudodifferential projections, and these projections have the same principal symbol. See e.g. \cite{BBW, Se, Se2}.}  (see Definition \ref{DefComm}).  Furthermore, from Proposition \ref{PropSeLag}, we have that $r(\ker D)$ is a Lagrangian subspace of the boundary data space $B^{s-1/p,p}(E_\Sigma)$, where the symplectic form on the Banach space $B^{s-1/p,p}(E_\Sigma)$ is given by Green's formula\footnote{For a general first order differential operator $A$ acting on sections $\Gamma(E)$ over a manifold $X$, Green's formula for $A$ is the adjunction formula
\begin{equation}
  (u,Av)_{L^2(X)} - (A^*u,v)_{L^2(X)} = \int_{\partial X} (r(u),-Jr(v)), \label{genGF2}
\end{equation}
where $A^*$ is the formal adjoint of $A$.  The map $J: E_{\partial X} \to E_{\partial X}$ is a bundle endomorphism on the boundary and it is determined by $A$. Hence, (\ref{genGF}) is an ``integration by parts" formula for $A$.  If $E$ is a Hermitian vector bundle, we will always take the real part of (\ref{genGF2}) in order to get a real valued pairing on the boundary.} for $\mathsf{D}$:
\begin{equation}
  \int_\Sigma \Re (u,-\mathsf{J}_0v) = \Re(u,\mathsf{D}v)_{L^2(Y)} - \Re(\mathsf{D}u,v)_{L^2(Y)}. \label{genGF}
\end{equation}
Summarizing, we have

\begin{Lemma}\label{LemmaKerP}
  The Cauchy data space $r(\ker \mathsf{D}) \subset B^{s-1/p,p}(E_\Sigma)$ is a Lagrangian subspace commensurate with $B^{s-1/p,p}\cZ^+_{\mathsf{B}_0}$.  Furthermore, for $s > 1/p$, the space $\ker \mathsf{D}$ is complemented in $B^{s,p}(E)$.
\end{Lemma}
The last statement follows from Corollary \ref{CorKerComp}.  Thus, while $r(\ker \mathsf{D})$ is a space determined by the entire operator $\mathsf{D}$ on $Y$, it is ``close" to the subspace $B^{s-1/p,p}\cZ^+_{\mathsf{B}_0}$, which is completely determined on the boundary.\\

Let us now apply the above general framework to our Hessian operators.  Let $\mathsf{B}$ denote the tangential boundary operator for $\tH_0$. By (\ref{deftH0}), $\mathsf{B}$ splits as a direct sum of the tangential boundary operators
\begin{align*}
  \mathsf{B}_\dgc: \Omega^1(\Sigma; i\R) \oplus \Omega^0(\Sigma; i\R) \oplus \Omega^0(\Sigma; i\R) & \circlearrowleft\\
  \mathsf{B}_\S: \Gamma(\S_\Sigma)  & \circlearrowleft,
\end{align*}
for $D_{\dgc}$ and $D_{B_\rf}$, respectively.  For the div-grad-curl operator $D_{\dgc}$, we can compute the tangential boundary operator and its spectrum rather explicitly.  As before, we work inside a collar neighborhood $[0,\eps] \times \Sigma$ of the boundary of $Y$, with the inward normal coordinate given by $t \in [0,\eps]$, and we choose coordinates so that the metric is of the form $dt^2 + g_t^2$, where $g_t$ is a family of Riemannian metrics on $\Sigma$. We can write $b \in \Omega^1(Y)$ as $b =  a + \beta dt$, where $a \in \Gamma([0,\eps), \Omega^1(\Sigma))$ and $\beta \in \Gamma([0,\eps), \Omega^0(\Sigma))$.  Let $\check{*}$ denote the Hodge star on $\Sigma$ with respect to $g_0$, and let $d_\Sigma$ be the exterior derivative on $\Sigma$.

So with the above notation, we have the following lemma concerning $D_\dgc$ (where for notational simplicity, we state the result for real-valued forms):

\setlength{\arraycolsep}{5pt}

\begin{Lemma}\label{LemmaDGC}
 Let $Y$ be a $3$-manifold with boundary $\Sigma$ oriented by the outward normal.  Then with respect to $(a,\beta,\alpha) \in \Gamma([0,\eps), \Omega^1(\Sigma) \oplus \Omega^0(\Sigma) \oplus \Omega^0(\Sigma))$ near the boundary, the div-grad-curl operator can be written as $D_{\dgc} = J_{\dgc}(\frac{d}{dt} + \mathsf{B}_{\dgc,t} + \mathsf{C}_{\dgc,t})$ as in (\ref{Diractan}), where\footnote{Note the signs, since $t$ is the \textit{inward} normal coordinate.}
 \begin{eqnarray}
   J_{\dgc} &=& \begin{pmatrix}
     -\check{*} & 0 & 0 \\ 0 & 0 & -1 \\ 0 & 1 & 0
   \end{pmatrix}, \label{Jdgc}\\
   B_{\dgc} = B_{\dgc,0} &=& \begin{pmatrix}
  0 & d_\Sigma & \check{*}d_\Sigma\\
  d_\Sigma^* & 0 & 0\\
  -\check{*}d_\Sigma & 0 & 0\end{pmatrix}. \label{Bdgc}
 \end{eqnarray}
The positive, negative, and zero eigenspace decompositions for $B_{\dgc}$ are given by
\begin{eqnarray}
\cZ^{\pm}_{\dgc} &=& \cZ_e^{\pm} \oplus \cZ_c^{\pm}\\
&:=& \mathrm{span}\left\{\begin{pmatrix}|\lambda|^{-1}d_\Sigma f_{\lambda^2} \\ \pm f_{\lambda^2} \\ 0\end{pmatrix}\right\} \oplus
\mathrm{span} \left\{ \begin{pmatrix} |\lambda|^{-1}\check{*}d_\Sigma f_{\lambda^2} \\ 0 \\ \pm f_{\lambda^2}\end{pmatrix} \right\}\\
\cZ^0_{\dgc} &=& H^1(\Sigma; \R) \oplus H^0(\Sigma; \R) \oplus H^0(\Sigma; \R),
\end{eqnarray}
where the $f_{\lambda^2}$ span the nonzero eigenfunctions of $\Delta = d_\Sigma^*d_\Sigma$ and $\Delta f_{\lambda^2} = \lambda^2f_{\lambda^2}$.

Let $\Omega^0_\bot(\Sigma) = \{\alpha \in \Omega^0(\Sigma) : \int\alpha = 0\}$ be the span of the nonzero eigenfunctions of $\Delta$.  Then for every $s \in \R$, and $1<p<\infty$, $B^{s,p}\cZ_e^{\pm}$ is the graph of the isomorphism $\pm d_\Sigma\Delta^{-1/2}: B^{s,p}\Omega^0_\bot(\Sigma) \to B^{s,p}\im d_\Sigma$.  Similarly, the spaces $B^{s,p}\cZ_c^{\pm}$ are graphs of the isomorphisms $\pm\check{*}d_\Sigma\Delta^{-1/2}: B^{s,p}\Omega^0_\bot(\Sigma) \to B^{s,p}\im \check{*}d_\Sigma$.
\end{Lemma}

\Proof The proof is by direct computation.\End


Altogether, we have the following spectral decompositions
\begin{align}
    \tT_\Sigma &= \cZ^+ \oplus \cZ^{-} \oplus \cZ^0,\\
    \Omega^1(\Sigma;i\R) \oplus \Omega^0(\Sigma;i\R)\oplus \Om &= \cZ^+_\dgc \oplus \cZ^-_\dgc \oplus \cZ^0_\dgc,\\
    \Gamma(\S_\Sigma) &= \cZ^+_\S\oplus\cZ^-_\S\oplus\cZ^0_\S,
\end{align}
corresponding to the positive, negative, and zero spectral subspaces of $\mathsf{B}$, $\mathsf{B_{\dgc}}$, and $\mathsf{B_\S}$, respectively. Since $\mathsf{B} = \mathsf{B}_\dgc \oplus \mathsf{B}_\S$, we obviously have
\begin{equation}
  \cZ^\bullet = \cZ^\bullet_\dgc \oplus \cZ^\bullet_\S, \qquad \bullet \in \{+,-,0\}.
\end{equation}
In particular, we have
\begin{align}
  \cZ^+ &= \cZ^+_\dgc \oplus \cZ_\S^+\\
  &= \cZ^+_e \oplus \cZ^+_c \oplus \cZ^+_\S, \label{Zspaces}
\end{align}
by Lemma \ref{LemmaDGC}. All the above decompositions hold when we take Besov closures.  In light of Lemma \ref{LemmaKerP}, the explicit decomposition (\ref{Zspaces}) will be important for us in the analysis to come.

Next, we work out the associated symplectic data for $\tH_0$ on the boundary, following the general picture described previously.  Namely, Green's formula (\ref{genGF}) for the Dirac operator $\tH_0$ induces a symplectic form on the boundary data space $\tT_\Sigma$.  Moreover, because $\tH_0$ is a Dirac operator, the endomorphism $-\mathsf{J}_0$ is a compatible complex structure for the symplectic form.  Explicitly, the symplectic form is \label{p:tomega}
\begin{align}
\tomega: \tT_\Sigma \oplus \tT_\Sigma & \to \R \nonumber \\
\tomega((a,\phi,\alpha_1,\alpha_0), (b,\psi,\beta_1,\beta_0)) &= \int_\Sigma a \wedge b + \int_\Sigma \Re(\phi,\rho(\nu)\psi)  - \int_\Sigma (\alpha_1\beta_0 - \alpha_0\beta_1), \label{tomega}
\end{align}
and the compatible complex structure is \label{p:tJSigma}
\begin{align}
 \tJ_\Sigma: \tT_\Sigma & \to \tT_\Sigma \nonumber \\
 (a,\phi,\alpha_1,\alpha_0) & \mapsto (-\check{*}a, -\rho(\nu)\phi, -\alpha_0, \alpha_1).
\end{align}
Observe that since $\tH_0 = \H_{(B_\rf,0)} \oplus -(d + d^*)$, the symplectic form and compatible complex structure above are a direct sum of those corresponding to the operators $\H_{(B_\rf,0)}$ and $(d + d^*)$.  In particular, Green's formula for $\H_{(B_\rf,0)} = *d \oplus D_{B_\rf}$ yields the symplectic form \label{p:omega}
\begin{align}
  \omega: \T_\Sigma \oplus \T_\Sigma & \to \R \nonumber \\
  \omega((a,\phi),(b,\psi)) &= \int_\Sigma a \wedge b + \int_\Sigma \Re(\phi,\rho(\nu)\psi)  \label{omega}
\end{align}
and compatible complex structure \label{p:JSigma}
\begin{align}
  J_\Sigma: \T_\Sigma &\to \T_\Sigma \nonumber\\
  (a,\phi) & \mapsto (-\check{*}a, -\rho(\nu)\phi).
\end{align}
Since the tangent space to $\fC(\Sigma)$ at any configuration is a copy of $\T_\Sigma$, we see that $\omega$ gives us a constant symplectic form on $\fC(\Sigma)$.  This symplectic form extends to $\fC^{0,2}(\Sigma)$, the $L^2$ closure of the configuration space, and since $B^{s,p}(\Sigma) \subseteq B^{0,2}(\Sigma) = L^2(\Sigma)$ for all $s > 0$ and $p \geq 2$, we also get a constant symplectic form on the Besov configuration spaces $\fC^{s,p}(\Sigma)$.  From now on, we will always regard $\fC^{s,p}(\Sigma)$ as being endowed with this symplectic structure.  Likewise, we always regard $\tT^{s,p}_\Sigma$ as being endowed with the symplectic form (\ref{tomega}). Indeed, the symplectic forms $\omega$ and $\tomega$ are the appropriate ones to consider, since they are the symplectic forms induced by the Hessian and augmented Hessian operators, respectively.\\

Having studied the particular augmented Hessian operator $\tH_0 = \tH_{(B_\rf,0)}$, we now study general augmented Hessian operators $\tH_\c$.  Here, $\c \in \fC^{s,p}(Y)$ is an arbitrary possibly nonsmooth configuration.  Suppose we have a bounded multiplication $B^{s,p}(Y) \times B^{t,q}(Y) \to B^{t-1,q}(Y)$, for some $t \in \R$ and $q \geq 2$.  It follows that $\tH_\c: \tT^{t,q} \to \tT^{t-1,q}$ and $\H_\c: \T^{t,q}\to\T^{t-1,q}$ are bounded maps.  To keep the topologies clear, we will often use the notation
\begin{align*}
  \H_\c^{t,q}: \T^{t,q} & \to \T^{t-1,q}\\
  \tH_\c^{t,q}: \tT^{t,q} &\to \tT^{t-1,q},
\end{align*}
so that the superscripts on the operators specify the regularity of the domains.  The next two lemmas tell us that $\ker \tH_\c^{t,q}$ and $r(\ker \tH_\c^{t,q})$ are compact perturbations of $\ker \tH_0^{t,q}$ and $r(\ker\tH_0^{t,q})$, respectively, for $(t,q)$ in a certain range.  We also give a more concrete description of this perturbation using Lemma \ref{LemmaFredproj}.

\begin{Lemma}\label{LemmaKerH}
  Let $s > 3/p$.  Let $\c \in \fC^{s,p}(Y)$ and suppose $t \in \R$ and $q \geq 2$ are such that we have a bounded multiplication map $B^{s,p}(Y) \times B^{t,q}(Y) \to B^{t'-1,q}(Y)$, where $t' > 1/q$ and $t \leq t' \leq t + 1$.
  \begin{enumerate}
    \item We have that $\ker \tH_\c^{t,q}$ is commensurate with $\ker \tH_0^{t,q}$ and the restriction map $r: \ker \tH_\c^{t,q} \to \tT^{t-1/q,q}_\Sigma$ is bounded.  More precisely, we have the decomposition
        \begin{equation}
          \ker \tH^{t,q}_\c = \{x + \wT x : x \in X_0'\} \oplus F, \label{graphofker}
        \end{equation}
        where $X_0' \subseteq \ker \tH_0^{t,q}$ has finite codimension, $\wT: X_0' \to \tT^{t',q}$, and $F \subseteq \tT^{t',q}$ is a finite dimensional subspace.  Moreover, one can choose as a complement for $X_0' \subset \ker \tH_0^{t,q}$ a space that is spanned by smooth elements.
    \item The space $\ker \tH_\c^{t,q}$ varies continuously\footnote{See Definition \ref{DefSS}.} with $\c \in \fC^{s,p}(Y)$.
  \end{enumerate}
\end{Lemma}

\Proof (i) Let $(b,\psi) = (B-B_\rf,\Psi)$.  The multiplication map $(b,\psi)\# = \tH_\c - \tH_0$ given by (\ref{sharp}) yields a bounded map
\begin{equation}
  (b,\psi)\# : \tT^{t',q} \to \tT^{t'-1,q} \label{compactmult}
\end{equation}
by hypothesis.  This map is a compact operator since it is the norm limit of $(b_i,\psi_i)\#$, with $(b_i,\psi_i)$ smooth.  Each of the operators $(b_i,\psi_i)\#$ is compact, since it is a bounded operator on $\tT^{t',q}$ and the inclusion $\tT^{t',q} \hookrightarrow \tT^{t'-1,q}$ is compact by Theorem \ref{ThmEmbed}.  Since the space of compact operators is norm closed, this proves (\ref{compactmult}) is compact.

Since $t' > 1/q$, then $\ker \tH^{t',q}_0$ is complemented in $\tT^{t',q}$ by Corollary \ref{CorKerComp}.  Let $X_1 \subset \tT^{t',q}$ be any such complement. Thus,
\begin{equation}
  \tH_0: X_1 \to \tT^{t'-1,q} \label{tH0X1}
\end{equation}
is injective.  It is also surjective by unique continuation, Theorem \ref{ThmUCPsurj}.  Hence (\ref{tH0X1}) is an isomorphism and the map
\begin{equation}
  \tH_\c: X_1 \to \tT^{t'-1,q}, \label{sec3.3:tHcX1}
\end{equation}
being a compact perturbation of an isomorphism, is Fredholm.  This allows us to write the kernel of $\tH_\c^{t,q}$ perturbatively as follows.

Let $x \in \ker \tH_\c^{t,q}$.  Then $\tH_0x = (b,\psi)\# x \in \T^{t'-1,q}$ and we can define
$$x_1 = -\left(\tH_0^{t',q}|_{X_1}\right)^{-1}(b,\psi)\# x \in X_1 \subset \tT^{t',q}.$$
Then if we define $x_0 = x - x_1 \in \tT^{t,q}$, we have
\begin{align*}
  \tH_0 x_0 &= \tH_0(x-x_1)\\
  &= \left(\tH_\c - (b,\psi)\#\right) x - \tH_0x_1\\
  &= -(b,\psi)\# x + (b,\psi)\# x\\
  & = 0.
\end{align*}
Hence, $x_0 \in \ker \tH^{t,q}_0$.  Thus, we have decomposted $x \in \ker\tH^{t,q}_\c$ as $x = x_0 + x_1$, where $x_0 \in \ker \tH^{t,q}_0$ is in the kernel of a smooth operator and $x_1 \in \tT^{t',q}$ is more regular (for $t' > t$).  We also have
\begin{align}
  0 &= \tH_\c x \nonumber \\
  & = \tH_\c (x_1 + x_0) \nonumber \\
  &= \tH_\c x_1 + (b,\psi)\# x_0. \label{x1x0}
\end{align}
By the above, we know that $\tH_\c: X_1 \to \tT^{t'-1,q}$ is Fredholm.  Thus, from (\ref{x1x0}), we see that there exists a subspace $X_0' \subseteq \ker \tH^{t,q}_0$ of finite codimension such that for all $x_0 \in X_0'$, there exists a solution $x_1 \in X_1$ to (\ref{x1x0}).  This solution is unique up to some finite dimensional subspace $F \subset X_1$; in fact $F$ is just the kernel of (\ref{sec3.3:tHcX1}).  This proves the decomposition
(\ref{graphofker}), where the map $\wT$ is given by
\begin{align}
  \wT: X_0' & \to X_1' \nonumber \\
  x_0 & \mapsto -(\tH_\c|_{X_1'})^{-1}(b,\psi)\# x_0, \label{Tmap}
\end{align}
where $X_1'$ is any complement of $F \subset X_1$.  The map $\wT$ is compact since the map $(b,\psi)\#$ is compact.  The rest of the statement now follows, since the restriction map $r: \ker \tH_0^{t,q} \to \tT^{t-1/q,q}_\Sigma$ is bounded by Theorem \ref{ThmSeeley}(i), and $r: \tT^{t',q} \to \tT^{t'-1/q,q}_\Sigma \subset \tT^{t-1/q,q}_\Sigma$ is bounded since $t' > 1/q$. Moreover, since smooth elements are dense in $\ker \tH^{t,q}_0$ by Corollary \ref{CorKerComp}, any finite dimensional complement for $X_0' \subseteq \ker \tH^{t,q}_0$ can be replaced by a complement that is spanned by smooth elements if necessary.

(ii) Let $\co \in \fC^{s,p}(Y)$.  By Definition \ref{DefSS}, we have to show that $\ker \tH^{t,q}_\c$ is a graph over $\ker \tH^{t,q}_\co$ for $\c$ close to $\co$.  We do the same thing as in (i).  Let $X_2$ be any complement of $\ker \tH^{t',q}_\co$ in $\tT^{t',q}$, which exists since $\ker \tH^{t',q}_\co$ is commensurate with $\ker \tH^{t',q}_0$ by (i), and the latter space is complemented.  Then $\tH_\co: X_2 \to \tT^{t'-1,q}$ is an isomorphism.  For $\c$ sufficiently close to $\co$, the map $\tH_\c: X_2 \to \tT^{t'-1,q}$ is injective, hence surjective (the index is invariant under compact perturbations), and therefore an isomorphism.  Then from the above analysis,
\begin{equation}
  \ker \tH^{t,q}_\c = \{x + \wT_\c x: x \in \ker \tH^{t,q}_\co\}, \label{kergraph2}
\end{equation}
where
\begin{align}
  \wT_\c: \ker \tH^{t,q}_\co & \to X_2\\
  x & \mapsto -(\tH_\c|_{X_2})^{-1}(b,\psi)\# x
\end{align}
and $(b,\psi) = (B-B_0,\Psi-\Psi_0)$. The map $\wT_{\c}$ varies continuously with $\c \in \fC^{s,p}(Y)$ near $\co$. \End

\begin{Rem}\label{Rem1/2}
In applications of the above lemma, instead of $(t,q)$ satisfying the very general hypothesis
\begin{equation}
\left.\begin{split}
\textrm{(i)\, } & t \in \R,\; q \geq 2,  \\
 \textrm{(ii)\, } & \textrm{the multiplication } B^{s,p}(Y) \times B^{t,q}(Y) \to B^{t'-1,q}(Y) \textrm{ is bounded,  where } t' > 1/q \\ & \textrm{and }
 t \leq t' \leq t + 1,
\end{split}\right\}\label{tq1}
\end{equation}
we will primarily only need the cases
\begin{equation}
  (t,q) \in \left\{(s+1,p),\, (s,p),\, (1/2,2)\right\}, \label{tq2}
\end{equation}
with corresponding values
\begin{equation}
  (t',q) \in \left\{(s+1,p),\, (s+1,p),\, (1/2+\eps,2)\right\}, \qquad \eps > 0.
\end{equation}
The last case of (\ref{tq2}) arises because we want to consider the space of boundary values in the $L^2$ topology, i.e., the spaces $\T^{0,2}_\Sigma$ and $\tT^{0,2}_\Sigma$.  In this particular case, the above lemma allows us to conclude that for $\c \in \fC^{s,p}(Y)$, we still get bounded restriction maps $r: \ker \tH_\c^{1/2,2} \to \tT^{0,2}_\Sigma$, just like in the case where $\c$ is smooth via Theorem \ref{ThmSeeley}. The boundedness of this map will be important when we perform symplectic reduction on Banach spaces in the proof of Theorem \ref{ThmLinLag}.  The case $(t,q) = (s+1,p)$ will be important for Proposition \ref{PropInvertH} and its applications in Section \ref{SecMonopoles}. In what follows, we will consider the operators $\tH^{s,p}_\c$ but they equally well apply to $\tH^{t,q}_\c$ in light of the analysis in Lemma \ref{LemmaKerH}, for $t,q$ satisfying (\ref{tq1}).
\end{Rem}

\begin{Lemma}\label{LemmaHessCD}
  Let $s>3/p$. For any $\c \in \fC^{s,p}(Y)$, we have the following:
  \begin{enumerate}
    \item The Cauchy data space  $r(\ker \tH_{\c}^{s,p})$ is a Lagrangian subspace of $\tT_\Sigma^{s-1/p,p}$ commensurate with $B^{s-1/p,p}\cZ^+$ and it varies continuously with $\c$.
    \item We have a direct sum decomposition $\tT_\Sigma^{s-1/p,p} = r(\ker \tH_{\c}^{s,p}) \oplus \tJ_\Sigma r(\ker \tH_{\c}^{s,p})$.
  \end{enumerate}
\end{Lemma}

\Proof (i) For any $\c \in \fC^{s,p}(Y)$, the space $r(\ker \tH_{\c}^{s,p})$ is isotropic since $\tH_{\c}^{s,p}$ is formally self-adjoint.  Since $s > 3/p$, then $\c \in L^\infty(Y)$ and we can apply the unique continuation theorem, Theorem \ref{ThmUCP}, which implies that $r: \ker \tH_\c^{s,p} \to \tT^{s-1/p,p}_\Sigma$ is injective.  In fact, it is an isomorphism onto its image, since this is true for $r: \ker \tH_0^{s,p} \to \tT^{s-1/p,p}_\Sigma$ (by Theorem \ref{ThmSeeley}(i) and unique continuation applied to the smooth operator $\tH_0^{s,p}$) and $\ker \tH_\c^{s,p}$ is a compact perturbation of $\ker \tH_0^{s,p}$ by Lemma \ref{LemmaKerH}. Hence, we get that $\tl^{s-1/p,p}_\c := r(\ker \tH_{\c}^{s,p})$ varies continuously with $\c$, since $\ker \tH_\c^{s,p}$ varies continuously by Lemma \ref{LemmaKerH} and $r: \ker \tH_\c^{s,p} \to \tT^{s-1/p,p}_\Sigma$ is an isomorphism onto its image. For $\c$ smooth, we know that $\tl^{s-1/p,p}_\c \subset \tT^{s-1/p,p}_\Sigma$ is a Lagrangian subspace by Proposition \ref{PropSeLag}.  By continuity then, $\tl^{s-1/p,p}_\c$ is a Lagrangian for all $\c \in \fC^{s,p}(Y)$.  Moreover, all the $\tl^{s-1/p,p}_\c$ are commensurate with one another, in particular, with $r(\ker \tH^{s-1/p,p}_0)$, and this latter space is commensurate with $B^{s-1/p,p}\cZ^+$ by Lemma \ref{LemmaKerP}.

(ii) When $\c$ is smooth, this follows from Proposition \ref{PropSeLag}.  Now we use the continuity of the Lagrangians with respect to $\c \in \fC^{s,p}(Y)$ for the general case.\End

We want to apply the previous results concerning the augmented Hessian $\tH_\c$ to deduce properties about the Hessian $\H_\c$.  To place these results in a context similar to the pseudodifferential picture in Appendix \ref{AppSeeley}, let us recall some more basic properties concerning the smooth operator $\tH_0$.  By Theorem \ref{ThmSeeley}, the operator $\tH_0$, by virtue of it being a smooth elliptic operator, has a Calderon projection $\tP^+_0$ \label{p:tP} and a Poisson operator $\tP_0$. These operators satisfy the following properties.  The map $\tP^+_0$ is a projection of the boundary data $\tT^{s-1/p,p}_\Sigma$ onto $r(\ker \tH_0^{s,p})$, the boundary values of $\ker \tH_0^{s,p}$, and the map $\tP_0$ is a map from the boundary data $\tT^{s-1/p,p}_\Sigma$ into $\ker(\tH_0^{s,p}) \subset \tT^{s,p}$.  Moreover, the maps $r: \ker \tH_0^{s,p} \to r(\ker \tH_0^{s,p})$ and $\tP_0: r(\ker \tH_0^{s,p}) \to \ker \tH_0^{s,p}$ are inverse to one another, and $r\tP_0 = \tP^+_0$.  This implies that the map $\tilde\pi_0 := \tP_0r: \tT^{s,p} \to \ker (\tH_0^{s,p})$ is a projection. We also have that $\im \tP^+_0 = r(\ker\tH_0^{s,p})$ is a Lagrangian subspace of $B^{s-1/p,p}\tT_\Sigma$ by Proposition \ref{PropSeLag}.

For a general nonsmooth $\c \in \fC^{s,p}$, we have $\tH^{s,p}_\c$ is a compact perturbation of the smooth elliptic operator $\tH_0^{s,p}$. The previous lemmas imply that $\ker \tH_\c^{s,p}$ and $r(\ker \tH_\c^{s,p})$ are compact perturbations of $\ker\tH_0^{s,p}$ and $r(\ker\tH_0^{s,p})$, respectively, and moreover, we still have unique continuation, i.e., $r: \ker \tH_\c^{s,p} \to \tT^{s-1/p,p}_\Sigma$ is an isomorphism onto its image.  It follows that there exists a Calderon projection $\tP_\c^+$ and Poisson operator $\tP_\c$ for $\tH_\c^{s,p}$ as well, which satisfy the same corresponding properties (see Lemma \ref{LemmaFredproj}).  We also have a projection $\tilde\pi_\c := \tP_\c r: \tT^{s,p} \to \ker \tH_\c^{s,p}$. We summarize this in the following lemma and diagram:

\begin{Lemma}\label{LemmatHCP}
  Let $s>3/p$ and $\c \in \fC^{s,p}(Y)$.  Then there exists a Calderon projection $\tP^+_\c: \tT^{s-1/p,p}_\Sigma \to r(\ker \tH^{s,p}_\c)$ and a Poisson operator $\tP_\c: \tT^{s-1/p,p}_\Sigma \to \ker\tH^{s,p}_\c$.  The maps $r: \ker \tH_\c \to r(\ker \tH^{s,p}_\c)$ and $\tP_\c: r(\ker\tH^{s,p}_\c) \to \ker \tH^{s,p}_\c$ are inverse to one another, and $r\tP_\c = \tP^+_\c$.
\end{Lemma}

\begin{equation}
  \begin{split}
  \xymatrix@C+.5cm@M+.2cm@!R{
\ker \tH_\c^{s,p} \ar@<-1ex>[d]_r  & \tT^{s,p} \ar[l]_{\hspace{.7cm}\tilde\pi_\c} \ar[d]^r  \ar[r]^{\tH_\c^{s,p}} & \tT^{s-1,p} \\
r(\ker \tH_\c^{s,p})\ar@<-1ex>[u]_{\tP_\c} & \tT_\Sigma^{s-1/p,p} \ar[l]^{\hspace{.6cm}\tP^+_\c}} \label{CD1}
\end{split}
\end{equation}
In studying the Hessian $\H^{s,p}_\c$ we want to establish similar results as in Lemmas \ref{LemmaHessCD} and \ref{LemmatHCP}.  These results are summarized in the main theorem of this section.


\begin{Theorem}\label{ThmLinLag}
  Let $s > \max(3/p,1/2)$ and let $\c \in \fM^{s,p}(Y)$. Suppose $\H_\c: \T^{s,p} \to \K^{s-1,p}_\c$ is surjective.\footnote{This holds under the assumption (\ref{assum}). See Lemma \ref{LemmaHtrans}.} Then we have the following:
  \begin{enumerate}
    \item The space $r_\Sigma(\ker \H_{\c}^{s,p})$ is a Lagrangian subspace of $\T_\Sigma^{s-1/p,p}$ commensurate with $B^{s-1/p,p}(\im d \oplus \cZ^+_\S)$.  Moreover, we have the direct sum decomposition
        \begin{equation}
          \T^{s-1/p,p}_\Sigma = r_\Sigma(\ker \H_{\c}^{s,p}) \oplus J_\Sigma r_\Sigma(\ker \H_{\c}^{s,p}). \label{eq3.3:decomp}
        \end{equation}
    \item Define \label{p:P}
  \begin{align}
    P^+_\c: \T_\Sigma^{s-1/p,p} & \to r_\Sigma(\ker \H_\c^{s,p})
  \end{align}
  to be the projection onto $r_\Sigma(\ker \H_\c^{s,p})$ through $J_\Sigma r_\Sigma(\ker \H_{\c}^{s,p})$ as given by (\ref{eq3.3:decomp}).  Let $\pi^+: \T_\Sigma \to \im d \oplus \cZ^+_\S$ denote the orthogonal projection onto $\im d \oplus \cZ^+_\S$ through the complementary space $\ker d^* \oplus (\cZ^-_\S \oplus \cZ^0_\S)$.  Then $\pi^+$, being a pseudodifferential projection, extends to a bounded map on $\T^{s-1/p,p}_\Sigma$, and it differs from the projection $P^+_\co$ by an operator
  \begin{equation}
    (P^+_\c - \pi^+): \T_\Sigma^{s-1/p,p} \to \T_\Sigma^{s-1/p+1,p}. \label{diffproj}
  \end{equation}
  which smooths by one derivative.
    \item There exists a unique operator
      \begin{align}
    P_\c: \T_\Sigma^{s-1/p,p} \to \ker (\H_\c|_{\C^{s,p}})
  \end{align}
 that satisfies $r_\Sigma P_\c = P^+_\c$.  The maps $r_\Sigma: \ker (\H_\c|_{\C^{s,p}}) \to r_\Sigma(\ker \H_{\c}^{s,p})$ and $P_\c: r_\Sigma(\ker \H_{\c}^{s,p}) \to \ker (\H_\c|_{\C^{s,p}})$ are inverse to one another.
 \end{enumerate}
 Furthermore, let $(B(t),\Psi(t))$ be a continuous (resp. smooth) path in $\fM^{s,p}(Y)$ such that $\H_{(B(t),\Psi(t))}: \T^{s,p} \to \K_{(B(t),\Psi(t))}^{s-1,p}$ is surjective for all $t$.
 \begin{enumerate}\stepcounter{enumi}\stepcounter{enumi}\stepcounter{enumi}
  \item  Then $\ker \H_{(B(t),\Psi(t))}^{s,p}$ and $r_\Sigma(\ker\H_{(B(t),\Psi(t))}^{s,p})$ are continuously (resp. smoothly) varying families of subspaces\footnote{See Definition \ref{DefSS}.}.  Consequently, the corresponding operators $P_{(B(t),\Psi(t)}^+$ and $P_{(B(t),\Psi(t)}$ vary continuously (resp. smoothly) in the operator norm topologies.
  \end{enumerate}
    Keeping $\fM^{s,p}(Y)$ fixed, the statements in (i), (iii), and (iv) remain true if we replace the $B^{s,p}(Y)$ and $B^{s-1/p,p}(\Sigma)$ topologies on all vector spaces with the $B^{t,q}(Y)$ and $B^{t-1/q,q}(\Sigma)$ topologies, respectively, where $t,q$ satisfy (\ref{tq2}) or more generally (\ref{tq1}).  If we do the same for (ii), everything also holds except that the map $(\ref{diffproj})$ smooths by $t'-t$ derivatives.
\end{Theorem}

\noindent The theorem implies we have the following corresponding diagram for the Hessian $\H_\c$:
\begin{equation}
  \begin{split}
  \xymatrix@C+.5cm@M+.2cm@!R{
\ker(\H_\c|_{\C^{s,p}}) \ar@<-1ex>[d]_{r_\Sigma}  & \T^{s,p} \ar[l]_{\hspace{1cm}\pi_\c} \ar[d]^{r_\Sigma}  \ar[r]^{\H_\c^{s,p}} & \T^{s-1,p} \\
r_\Sigma(\ker \H_\c^{s,p}) \ar@<-1ex>[u]_{P_\c} & \T_\Sigma^{s-1/p,p} \ar[l]^{\hspace{.7cm}P^+_\c}} \label{CD2}
\end{split}
\end{equation}
Here, $\pi_\c := P_\c r_\Sigma$ is a projection of $\T^{s,p}$ onto $\ker(\H_\c|_{\C^{s,p}})$.

\begin{Def}\label{DefThm}
 By abuse of language, we call the operators $P^+_\c$ and $P_\c$ defined in Theorem \ref{ThmLinLag} the \textit{Calderon projection} and \textit{Poisson operator} associated to $\H_\c^{s,p}$, respectively (even though $\H_\c^{s,p}$ is not an elliptic operator), due to their formal resemblance to Calderon and Poisson operators for elliptic operators (as seen in the diagrams (\ref{CD1}) and (\ref{CD2})).
\end{Def}

Note that the Calderon projection $P^+_\c$ and Poisson operator $P_\c$ we define above are unique, since we specified their kernels.  In the general situation of an elliptic operator (such as $\tH_\c^{s,p}$ above) one usually only specifies the range of the Calderon projection, in which case, the projection is not unique (see also Remark \ref{RemCP}).  Our particular choice of kernel for $P^+_\c$ is made so that $P^+_\c$ is nearly pseudodifferential, in the sense of the smoothing property (\ref{diffproj}).  This property will be used in \cite{N2}, where analytic properties of the tangent spaces to the Lagrangian $\L^{s-1/p,p}$ and the projections onto them become crucial.

\begin{Rem}\label{RemCont}
  The continuous (resp. smooth) dependence of $P^+_\c$ and $P_\c$ in Theorem \ref{ThmLinLag}(iii) with respect to $\c$, as well as all other continuous dependence statements appearing in the rest of this paper, will only attain their true significance in \cite{N2}.  There, we will consider paths of configurations, and so naturally, we will have to consider time-varying objects.  For brevity, we will only make statements regarding continuous dependence from now on, though they can all be adapted to smooth dependence with no change in argument.
\end{Rem}

Proving Theorem \ref{ThmLinLag} is essentially deducing diagram (\ref{CD2}) from diagram (\ref{CD1}). Let us first make sense of the hypotheses of the theorem.  From Lemma \ref{LemmaTdecomp}, in order for $\K^{s-1,p}_\c$ to be well-defined when $\c \in \fC^{s,p}(Y)$, we need $s > \max(1-s, 3/p)$, which means we need $s > \max(3/p,1/2)$.  This explains the first hypothesis.  Next, observe that for $\c \in \M(Y)$ a smooth monopole, we have
\begin{align}
  \J_{\c} &\subseteq \ker \H_\c, \label{HkillsJ}\\
  \im \H_\c & \subseteq \K_\c. \label{imHinK}
\end{align}
One can verify this directly by a computation or reason as follows.  As previously discussed, the Seiberg-Witten map (\ref{SW3map}) is gauge-equivariant and hence its set of zeros is gauge-invariant.  Thus, the derivative of $SW_3$ along the gauge-orbit of a monopole vanishes.  This is precisely (\ref{HkillsJ}).  For (\ref{imHinK}), observe that the range of $\H_\c$ annihilates $\J_{\c,t}$ by (\ref{HkillsJ}) and since $\H_\c$ is formally self-adjoint.  From the orthogonal decomposition $\T = \J_{\c,t} \oplus \K_\c$, we conclude that $\im\H_\c \subseteq \K_\c$.
We want to establish similar properties on Besov spaces.  Namely, we want
\begin{align}
  \J_{\c}^{s,p} &\subseteq \ker \H_\c^{s,p}, \label{HkillsJ2}\\
  \im \H_\c^{s,p} &\subseteq \K_\c^{s-1,p}. \label{imHinK2}
\end{align}
However, this follows formally from (\ref{HkillsJ}) and (\ref{imHinK}) as long as we can establish on Besov spaces the appropriate mapping properties of the differentiation and multiplication involved in verifying (\ref{HkillsJ2}) and (\ref{imHinK2}) directly. Thus, (\ref{HkillsJ2}) holds because the map $\tH_\c: \T^{s,p} \to \T^{s-1,p}$ is bounded when $s > 3/p$.  Likewise, (\ref{imHinK}) holds because $\bd_\c^*: \T^{s-1,p} \to B^{s-2,p}\Omega^0(Y; i\R)$ is bounded when $s > \max(3/p,1/2)$.  In drawing these conclusions, as done everywhere else in this paper, we make essential use of Corollary \ref{CorDiffMap} and Theorem \ref{ThmMult}.

Thus, from (\ref{imHinK2}), we see that the hypotheses of Theorem \ref{ThmLinLag} make sense. In fact, for $\c \in \fM^{s,p}(Y)$, we have the following result concerning the range of $\H_\c^{s,p}$:

\begin{Lemma}\label{LemmaHessK}
  Let $\c \in \fM^{s,p}(Y)$. Then $\im \H_{\c}^{s,p} \subseteq \K_{\c}^{s-1,p}$ and $\H_{\c}^{s,p}: \T^{s,p} \to \K_{\c}^{s-1,p}$ has closed range and finite dimensional cokernel.
\end{Lemma}

\Proof It remains to prove the final statement. Pick any elliptic boundary condition for the operator $\tH_\c^{s,p}$ such that one of the boundary conditions for $(b,\psi,\alpha) \in \tT^{s,p}$ is $\alpha|_\Sigma = 0$.  Such a boundary condition is possible, since the subspace $B^{s-1/p,p}(\T_\Sigma\oplus \Omega^0(\Sigma;i\R) \oplus 0)$ of $\tT_\Sigma^{s-1/p,p}$ with vanishing $0\oplus0\oplus B^{s-1/p,p}\Omega^0(\Sigma)$ component contains subspaces Fredholm\footnote{See Definition \ref{DefFred}.} with $r(\ker \H_\c^{s,p})$ by Lemma \ref{LemmaKerP}, Lemma \ref{LemmaDGC}, and (\ref{Zspaces}).  For such a boundary condition, observe that $\im (\tH^{s,p}_\c) \cap \K^{s-1,p}_\c \subseteq \im \H^{s,p}_\c$.  This is because if $d\alpha \in \K^{s-1,p}_\c$ with $\alpha|_\Sigma = 0$, then $d\alpha = 0$.  Since we chose elliptic boundary conditions for $\tH^{s,p}_\c$, this means $\im \tH^{s,p}_\c \subseteq \tT^{s-1,p}_\c$ is closed and has finite codimension, which implies $\im (\tH^{s,p}_\c) \cap \K^{s-1,p}_\c$ is also closed and has finite codimension in $\K^{s-1,p}_\c$.  Hence, the same is true for $\im \H^{s,p}_\c$.\End

Next, we relate the kernel of $\H_\c^{s,p}$ to the kernel of $\tH_\c^{s,p}$ along with their respective boundary values.

\begin{Lemma}\label{LemmaParts}
  Let $s > \max(3/p,1/2)$ and $\c \in \fM^{s,p}(Y)$.
  \begin{enumerate}
    \item We have a decomposition
  \begin{align}
    \ker \H_\c^{s,p} &= \ker (\H_\c|_{\C^{s,p}}) \oplus \J^{s,p}_{\c,t} \label{HCJ} \\
    \ker \tH_\c^{s,p} &= \ker (\H_\c|_{\C^{s,p}}) \oplus \Gamma_0^{s,p}, \label{tHCG}
  \end{align}
  where $\Gamma_0^{s,p} \subseteq \tT^{s,p}$ is the graph of a partially defined map $\Theta_0: \ker \Delta \dashrightarrow \T^{s,p}$, where $\Delta$ is the Laplacian on $B^{s,p}\Omega^0(Y; i\R)$, and the domain of $\Theta_0$ has finite codimension.
  \item We have
   \begin{align}
     r_\Sigma(\ker \H^{s,p}_\c) &= r_\Sigma(\ker (\H_\c|_{\C^{s,p}})) \label{kerH=kerHC}\\
     r(\ker \tH_\c^{s,p}) &= r(\ker (\H_\c|_{\C^{s,p}})) \oplus \check\Gamma_0^{s-1/p,p},
   \end{align}
    where $\check\Gamma_0^{s-1/p,p} = r(\Gamma_0^{s,p})$ is the graph of a partially defined map $\check\Theta_0: 0\oplus0\oplus B^{s-1/p,p}\Omega^0(\Sigma;i\R) \dashrightarrow \T^{s-1/p,p}_\Sigma$ and the domain of $\check\Theta_0$ has finite codimension.
  \item We have $r(\ker \H_\c^{s,p})$ is commensurate with $B^{s-1/p,p}(\cZ^+_e \oplus \cZ^+_\S)$ and $r_\Sigma(\ker \H_\c^{s,p})$ is commensurate with $B^{s-1/p,p}(\im d \oplus \cZ^+_\S)$.
  \end{enumerate}
\end{Lemma}

\Proof (i) The first decomposition (\ref{HCJ}) follows from (\ref{JtC}) and $\J^{s,p}_{\c,t} \subset \ker \H^{s,p}_\c$.  For (\ref{tHCG}), observe that $\ker (\H_\c|_{\C^{s,p}_\c}) = \ker \tH_\c|_{\T^{s,p}}$.  Thus, the elements of  $\ker\tH_\c^{s,p}$ that do not lie in $\ker \tH_\c|_{\T^{s,p}}$ have nonzero $B^{s,p}\Omega^0(Y; i\R)$ component.  To find them, we need to solve the equation
\begin{equation}
  \H_\c(b,\psi) - d\alpha = 0, \label{solvealpha}
\end{equation}
with $\alpha$ nonzero.  Since $\im \H_\c^{s,p} \subseteq \K^{s-1,p}_\c$ by the previous lemma, we need $d\alpha \in \K^{s-1,p}_\c$, whence $\alpha \in \ker \Delta$.  Since $\im \H_\c^{s,p}$ has finite codimension in $\K^{s-1,p}_\c$ by Lemma \ref{LemmaHessK}, then (\ref{solvealpha}) has a solution $(b,\psi)$ for all $\alpha$ in some subspace of $\ker \Delta$ of finite codimension.  The $(b,\psi)$ is unique up to an element of $\ker \H_\c^{s,p}$.  Thus, picking a complement\footnote{The reasoning used in the proof of Lemma \ref{LemmaHessK} shows that $\H_\c: \T^{s,p} \to \K^{s-1,p}_\c$ has a right parametrix.  This implies that $\ker \H_\c^{s,p} \subset \T^{s,p}$ is complemented.} of $\ker \H_\c^{s,p}$ in $\T_\c^{s,p}$ specifies for us a map $\Theta_0: \ker \Delta \dashrightarrow \T^{s,p}$ whose graph $\Gamma_0^{s,p}$ is a  complementary subspace of $\ker (\tH_\c|_{\T^{s,p}_\c})$ in $\ker \tH_\c^{s,p}$, and which parametrizes solutions to (\ref{solvealpha}).

(ii) This follows from applying $r_\Sigma$ and $r$ to (i).  The graph property of $\check\Gamma_0^{s-1/p,p}$ comes from noting that any element of $\Gamma_0^{s,p}$ is uniquely determined by the $0 \oplus 0 \oplus B^{s-1/p,p}\Omega^0(\Sigma;i\R)$ component of its image under $r$.  This follows from considering the homogeneous Dirichlet problem for $\Delta$, namely
\begin{eqnarray*}
\Delta\alpha &= 0\\
\alpha|_\Sigma &= \beta.
\end{eqnarray*}
This problem has a unique solution for every $\beta$.

(iii) By Lemma \ref{LemmaKerP}, we have $r(\ker \tH^{s,p}_\c)$ is commensurate with $B^{s-1/p,p}\cZ^+$.  Let
$$\pi_0: \tT^{s-1/p,p}_\Sigma \to 0\oplus 0 \oplus B^{s-1/p,p}\Omega^0(\Sigma; i\R)$$
denote the coordinate projection onto the last $0$-form factor in $\tT^{s-1/p}_\Sigma$.  By Lemma \ref{LemmaDGC} and (ii), we have
\begin{align*}
  \pi_0: B^{s-1/p,p}\cZ^+_c &\to 0\oplus 0\oplus B^{s-1/p,p}\Omega^0(\Sigma; i\R)\\
  \pi_0: B^{s-1/p,p}\check\Gamma_0 &\to 0\oplus 0\oplus  B^{s-1/p,p}\Omega^0(\Sigma; i\R)
\end{align*}
 are Fredholm.  We now apply Lemma \ref{LemmaComm} with $X = \tT^{s-1/p,p}_\Sigma$ and complementary subspaces
\begin{align*}
  X_1 &= \T_\Sigma^{s-1/p,p} \oplus B^{s-1/p,p}\Omega^0(\Sigma; i\R) \oplus 0,\\
  X_0 &= 0 \oplus 0 \oplus B^{s-1/p,p}\Omega^0(\Sigma; i\R).
\end{align*}
Let $U = r(\ker \tH^{s,p}_\c)$ and $V = B^{s-1/p,p}\cZ^+$ in the lemma. Then from that lemma and (ii), we conclude that $r(\ker (\H_\c|_{\C^{s,p}})) = U \cap X_1$ is commensurate with $B^{s-1/p,p}(\cZ^+_e \oplus \cZ^+_\S) =  V \cap X_1.$
This proves the first part of (iii).  For the second part, consider the coordinate projection of $X_1$ onto
$\T_\Sigma^{s-1/p,p}$.  This restricts to an isomorphism of $V \cap X_1$ onto its image $B^{s-1/p,p}(\im d \oplus \cZ^+_\S)$, by Lemma \ref{LemmaDGC}.  It follows that this projection maps $U \cap X_1$ onto a space commensurate with $V \cap X_1$, and this space is precisely $r_\Sigma(\ker \H_\c^{s,p})$.\End


\begin{Corollary}\label{CorSur}
  Let $\c \in \fM^{s,p}(Y)$ and suppose $\H_\c: \T^{s,p} \to \K_\c^{s-1,p}$ is surjective. Then
  \begin{enumerate}
    \item the maps $\Theta_0^{s,p}$ and $\check\Theta_0^{s-1/p,p}$ are defined everywhere;
    \item $r_\Sigma: \ker (\H_\c|_{\C^{s,p}}) \to \T_\Sigma^{s-1/p,p}$ is an isomorphism onto its image.
  \end{enumerate}
\end{Corollary}

\Proof (i) This follows from the constructions of $\Theta_0$ and $\check\Theta_0$ in the previous lemma.

(ii) By unique continuation, the map $r: \ker (\tH^{s,p}_\c) \to \tT^{s-1/p,p}_\Sigma$ is an isomorphism onto its range.  By restriction, it follows that
\begin{equation}
  r: \ker (\tH_\c|_{\C^{s,p}}) \to \tT^{s-1/p,p}_\Sigma \label{rinj2}
\end{equation}
is injective.  To prove (ii), it suffices to show that
\begin{equation}
  r_\Sigma: \ker (\tH_\c|_{\C^{s,p}}) \to \T^{s-1/p,p}_\Sigma \label{rinj}
\end{equation}
is injective, since $\tH_\c|_{\C^{s,p}} = \H_\c|_{\C^{s,p}}$.  So suppose (\ref{rinj}) is not injective.  Since (\ref{rinj2}) is injective, this means there is an element of the form $((0,0),\alpha,0) \in r(\ker (\tH_\c|_{\C^{s,p}}))$ with $\alpha \in B^{s-1/p,p}\Omega^0(\Sigma)$ nonzero.  On the other hand, $r(\ker \tH^{s,p}_\c)$ is a Lagrangian subspace of $\tT_\Sigma^{s-1/p,p}$  by Lemma \ref{LemmaHessCD}.  This contradicts (i), since if $\check\Theta_0^{s-1/p,p}$ is defined everywhere, then $(0,0,\alpha,0)$ cannnot symplectically annihilate $\check\Gamma_0^{s-1/p,p}$.  Indeed, the spaces $0\oplus\Omega^0(\Sigma)\oplus 0$ and $0\oplus0\oplus\Omega^0(\Sigma)$ are symplectic conjugates with respect to the symplectic form (\ref{tomega}).\End

\textbf{Proof of Theorem \ref{ThmLinLag}:}\hspace{0.3cm}
(i) We will apply the method of symplectic reduction, via Theorem \ref{ThmSR} and Corollary \ref{CorSR}.  By Lemma \ref{LemmaKerH}, we may consider the operators $\tH_\c^{1/2,2}$ and $\H^{1/2,2}_\c$, their kernels, and the restrictions of these latter spaces to the boundary.  Indeed, let us verify the hypotheses of Lemma \ref{LemmaKerH}.  Since $p \geq 2$, we have the embedding $B^{s,p}(\Sigma) \subseteq B^{s-\eps,2}(\Sigma)$ for any $\eps > 0$ by Theorem \ref{ThmEmbed}.  Choose $\eps$ small enough so that $s - \eps > 1/2+\eps$.  Then $(t,q) = (1/2,2)$ and $t' = \frac{1}{2}+\eps$ satisfies the hypotheses of Lemma \ref{LemmaKerH} since we have $B^{s-\eps,2}(Y) \times B^{1/2,2}(Y) \to B^{t'-1,2}(Y)$.

Let $U = L^2(\T_\Sigma \oplus \Omega^0(\Sigma;i\R) \oplus 0)$.  It is a coisotropic subspace of the strongly symplectic Hilbert space $L^2\tT_\Sigma = \tT_\Sigma^{0,2}$.  If we apply Theorem \ref{ThmSR} to the Lagrangian $L = r(\ker \tH_\c^{1/2,2})$, the symplectic reduction of $L$ with respect to $U$ is precisely $r_\Sigma(\ker \tH^{1/2,2}_\c)$ by Lemma \ref{LemmaParts}(ii).  It follows that $r_\Sigma(\ker \H^{1/2,2}_\c)$ is a Lagrangian inside $U \cap \tJ_\Sigma U = L^2\T_\Sigma$.  We would like to make the corresponding statement in the Besov topologies. By Lemma \ref{LemmaParts}(iii), we know that $r_\Sigma(\ker \H^{s,p}_\c)$ is commensurate with $B^{s-1/p,p}(\im d \oplus \cZ^+_\S)$.  On the other hand, we have that $B^{s-1/p,p}(\im d \oplus \cZ^+_\S)$ and $J_\Sigma B^{s-1/p,p}(\im d \oplus \cZ^+_\S)$ are Fredholm in $\T^{s-1/p,p}_\Sigma$.  Indeed, the Hodge decomposition implies $\im d$ and $\mathrm{im} *d$ are Fredholm in $B^{s-1/p,p}\Omega^1(\Sigma;i\R)$, and since $\rho(\nu)$ interchanges the positive and negative eigenspaces $\cZ^+_\S$ and $\cZ^-_\S$ of the tangential boundary operator $\mathsf{B}_\S$ associated to the spinor Dirac operator $D_{B_\rf}$, we have that the $B^{s-1/p,p}(\Sigma)$ closures of $\cZ^+_\S$ and $\rho(\nu)\cZ^+_S = \cZ^-_\S$ are Fredholm in $B^{s-1/p,p}\Gamma(\S)$.  That these decompositions are Fredholm in Besov topologies follows from the fact these spaces are given by the range of pseudodifferential projections whose principal symbols are complementary projections, and pseudodifferential operators are bounded on Besov spaces.  We now apply Corollary \ref{CorSR}, with $X = \T^{0,2}_\Sigma$ and $Y = \T^{s-1/p,p}_\Sigma$, to conclude that $r_\Sigma(\ker \H^{s,p}_\c)$ is a Lagrangian subspace of $B^{s-1/p,p}\T_\Sigma$.

(ii) By Lemma \ref{LemmaHessCD} and (i), $r(\ker \tH^{s,p}_0)$ is commensurate with $B^{s-1/p,p}\cZ^+$ and $r_\Sigma(\ker \H^{s,p}_0)$ is commensurate with $B^{s-1/p,p}(\im d \oplus \cZ_S^+)$, respectively.  Since $\tH_0$ is smooth, then we can even say more: there exist pseudodifferential projections onto $r(\ker \tH^{s,p}_0)$ and $B^{s-1/p,p}\cZ^+$ that have the same principal symbol, which means that their difference is a pseudodifferential operator of order $-1$.  It follows that the projection of $r(\ker \tH^{s,p}_0)$ onto any complement\footnote{More precisely, in what follows, when we speak of some unspecified complementary subspace, we mean one defined by a pseudodifferential projection. This is convenient because pseudodifferential operators preserve regularity, i.e., they map $B^{t,q}(\Sigma)$ to itself for all $t,q \in \R$, and so we never lose any smoothness once we have gained it.} of $B^{s-1/p,p}\cZ^+$ is smoothing of order one.  Consequently, letting $U^{s-1/p,p} = B^{s-1/p,p}(\T_\Sigma \oplus \Omega^0(\Sigma;i\R) \oplus 0)$, then the projection of $r(\ker \tH^{s,p}_0) \cap U$ onto any complement of $B^{s-1/p,p}\cZ^+ \cap U$ is smoothing of order one.  (Here, we use the fact that $U^{s-1/p,p} + B^{s-1/p,p}\cZ^+$ has finite codimension in $\tT^{s-1/p,p}_\Sigma$.) Applying symplectic reduction with respect to $U^{s-1/p,p}$, it follows that the projection of $r_\Sigma(\ker \H^{s,p}_0)$ onto any complement of $B^{s-1/p,p}(\im d \oplus \cZ^+_\S)$ is smoothing of order one.

For a nonsmooth configuration $\c$, we also want to show that the projection of $r_\Sigma(\ker \H^{s,p}_\c)$ onto any complement of $B^{s-1/p,p}(\im d \oplus \cZ^+_\S)$ is smoothing of order one.  For then this will imply the corresponding property with respect to the pair of spaces $J_\Sigma(r_\Sigma(\ker \H^{s,p}_\c))$ and $B^{s-1/p,p}J_\Sigma(\im d \oplus \cZ^+_\S)$, the latter being of finite codimension in $B^{s-1/p,p}(\ker d^* \oplus (\cZ^+_\S \oplus \cZ^0_\S))$.  We can then apply Lemma \ref{LemmaFredDecomp}(ii) while noting Remark \ref{RemError}, to conclude that the projection $P^+_\co$ differs from $\pi^+$ by a operator that is smoothing of order one.

Thus, by our first step, it suffices to show that the projection of $r_\Sigma(\ker \H^{s,p}_\c)$ onto any complement of $r_\Sigma(\ker \H^{s,p}_0)$ smooths by one derivative.  This follows however from Lemma \ref{LemmaKerH}.  Indeed, we can take $(t,q) = (s,p)$ and $t' = s+1$ in Lemma \ref{LemmaKerH}, and since there exists a projection of $\ker \tH^{s,p}_\c$ onto a complement of $\ker \tH^{s,p}_0$ that smooths by one derivative, the corresponding statement is true for the spaces $r_\Sigma(\ker \H^{s,p}_\c)$ and $r_\Sigma(\ker\H^{s,p}_0)$.  Here, it is important that all finite dimensional errors involved are spanned by elements that are smoother by one derivative (so that the finite rank projection onto the space spanned by these elements smooths by one derivative), which is guaranteed by Lemma \ref{LemmaKerH}. From these properties, one can now apply Lemma \ref{LemmaFredDecomp}(ii), with
\begin{align*}
    X &= \T^{s-1/p,p}_\Sigma\\
  U_0 &= B^{s-1/p,p}\big(\im d \oplus \cZ^+_\S\big)\\
  U_1 &= B^{s-1/p,p}\big(\ker d^* \oplus (\cZ^-_\S \oplus \cZ^0_\S)\big)\\
  V_1 &= r_\Sigma(\ker \H^{s,p}_\c)\\
  V_2 &= J_\Sigma(r_\Sigma(\ker \H^{s,p}_\c)).
\end{align*}
In our case, we know that $X = U_0 \oplus U_1 = V_0 \oplus V_1$, and that the $U_i$ and $V_i$ are commensurate, $i = 0,1$, where the compact error is smoothing of order one.  Thus, by Remark \ref{RemError}, $P^+_\co = \pi_{V_0,V_1}$ and $\pi^+ = \pi_{U_0,U_1}$ differ by an operator that smooths of order one.

(iii) Let
$$\pi_{SR}: r(\ker \H^{s,p}_\c) \to r_\Sigma(\ker \H_\c^{s,p})$$
be the symplectic reduction as in (i), i.e., the map $\pi_{SR}$ is the map which projects $r(\ker \H^{s,p}_\c) \subset \tT_\Sigma^{s-1/p,p}$ onto $r_\Sigma(\ker\H_\c^{s,p})$, induced by the projection $\tT_\Sigma^{s-1/p,p} \to \T_\Sigma^{s-1/p,p}$ onto the first factor.
This map is an isomorphism by Corollary \ref{CorSur}(ii).  Hence, $\pi_{SR}^{-1}$ exists and is bounded.  Define
\begin{align}
P_{\c} &= \tP_{\c}(\pi_{SR})^{-1}P^+_\c,
\end{align}
where $\tP_\c$ is the Poisson operator of $\tH_\c^{s,p}$. By construction, $P^+_\c: \T^{s-1/p,p}_\Sigma \to r_\Sigma(\ker \H^{s,p})$ and $\tilde P_{\c}(\pi_{SR})^{-1}: r_\Sigma(\ker \H^{s,p}) \to \ker(\H_\c|_{\C^{s,p}})$.  Thus, $P_\c: \T^{s-1/p,p}_\Sigma \to \ker(\H_\c|_{\C^{s,p}})$ and $r_\Sigma P_\c = P^+_\c$.  Moreover, from Corollary \ref{CorSur}(ii), it follows that $P_\c: \T^{s-1/p,p}_\Sigma \to \ker(\H_\c|_{\C^{s,p}})$ and $r_\Sigma: \ker(\H_\c|_{\C^{s,p}}) \to \T^{s-1/p,p}_\Sigma$ are inverse to each other.

(iv) We establish the smooth case, with the continuous case being exactly the same. It is easy to check that all the subspaces and operators involved in the construction of the maps in (ii) and (iii) vary smoothly with $\ct$.  Indeed, since $\K^{s-1,p}(Y)$ is a bundle, by Proposition \ref{PropBundle}, we can locally identify its fibers, i.e., the maps
$$\Pi_{\K^{s-1,p}_\co}: \K^{s-1,p}_{\ct} \to \K^{s-1,p}_{\co}$$
are all isomorphisms for all $\ct$ sufficiently $B^{s,p}(Y)$ close to a fixed $\co$.  Then restricting to $t$ on a small interval for which this is the case, then we have $\ker \H_{\ct} = \ker \big(\Pi_{\K^{s-1,p}_\co}\H_{\ct}\big)$, and $\Pi_{\K^{s-1,p}_\co}\H_{\ct}^{s,p}: \T^{s,p}\to\K^{s-1,p}_\co$ are all surjective for all $t$.  From this, it follows that $\ker \H_{\ct}^{s,p}$ varies smoothly, and since $\J^{s,p}_\ct \in \ker \H_{\ct}^{s,p}$ for all $t$, this implies $\ker (\H_{\ct}|_{\C^{s,p}}$ vary smoothly.  Indeed, one argues as in Lemma \ref{LemmaKerH} for the continuity of $\ker \tH_\c^{s,p}$ with respect to $\c$, only now we have in addition that all objects vary smoothly. Since $r_\Sigma: \ker (\H_{\ct}|_{\C^{s,p}}) \to \T_\Sigma^{s-1/p,p}$ is an isomorphism onto its image for all $t$, it follows that $r_\Sigma(\ker \H^{s,p}_{\ct})$ varies smoothly.  Since this holds for all $t$ on small intervals, it holds for all $t$ along the whole path.

To prove the final statement, we observe that all the above methods apply to $\tH^{t,q}_\c$ and $\H^{t,q}_\c$ without modification in light of Lemma \ref{LemmaKerH}.  See also Remark \ref{Rem1/2}.\End

We conclude this section with some important results that will be used later.

\begin{Lemma}
Let $\c \in \fM^{s,p}(Y)$, assume all the hypotheses of Theorem \ref{ThmLinLag}, and suppose $(t,q)$ satisfies (\ref{tq2}) or more generally (\ref{tq1}).  Then the space
  \begin{equation}
    \mathsf{L}_\c^{t-1/q,q} := J_\Sigma r_\Sigma(\ker \H^{t,q}_\c) \oplus B^{t-1/q,q}\Omega^0(\Sigma) \oplus 0 \label{complag}
  \end{equation}
  is a complementary Lagrangian for $r(\ker \tH^{t,q}_\c)$ in $\tT^{t-1/q,q}_\Sigma$.  The space $\mathsf{L}_\c^{t-1/q,q}$ varies continuously with $\c \in \fM^{s,p}(Y)$ (as long as $\H_\c^{s,p}: \T^{s,p} \to \K^{s-1,p}_\c$ is always surjective).
\end{Lemma}

\Proof By Theorem \ref{ThmLinLag}(i), $J_\Sigma r_\Sigma(\ker \H^{t,q}_\c)$ and $r_\Sigma (\ker \H^{t,q}_\c)$ are complementary Lagrangians in $\T^{t-1/q,q}_\Sigma$.  By Lemma \ref{LemmaParts}(ii) and Corollary \ref{CorSur}(i), it is now easy to see that (\ref{complag}) is a complement of $r(\ker \tH^{t,q}_\c)$ in $\tT^{t-1/q,q}_\Sigma$.  Since $r_\Sigma(\ker \H^{t,q}_\c)$ depends continuously on $\c \in \fM^{s,p}(Y)$ by Theorem \ref{ThmLinLag}(iv), the last statement follows.\End

For $t > 1/q$, define
\begin{equation}
  \tX_\c^{t,q} = \{(b,\psi,\alpha) \in \tT^{t,q} : r(b,\psi,\alpha) \in J_\Sigma r_\Sigma(\ker \H^{t,q}_\c) \oplus B^{t-1/q,q}\Omega^0(\Sigma) \oplus 0\}, \label{Xdomain}
\end{equation}
the subspace of $\tT^{t,q}$ whose boundary values lie in (\ref{complag}).  Likewise, define \label{p:tX}
\begin{equation}
  X_\c^{t,q} = \C^{t,q} \cap \tX_\c^{t,q} \subset \T^{t,q}. \label{Xtq}
\end{equation}
By the above lemma, the domains $\tX_\c^{t,q}$ and $X_\c^{t,q}$ are such that their boundary values under $r$ and $r_\Sigma$ are complementary to the boundary values of $\ker \tH_\c^{t,q}$ and $\ker \H_\c^{t,q}$, respectively.  Thus, we expect these domains to be ones on which the operators $\tH_\c^{t,q}$ and $\H_\c^{t,q}$ are invertible elliptic operators.  This is exactly what the following proposition tells us.

\begin{Proposition}\label{PropInvertH}
Let $\c \in \fM^{s,p}(Y)$ and assume all the hypotheses of Theorem \ref{ThmLinLag}.  Let $t > 1/q$ and $q \geq 2$ satisfy (\ref{tq2}) or more generally (\ref{tq1}).  Then the maps
   \begin{align}
     \tH_\c: \tX^{t,q}_\c & \to \tT^{t-1,q}, \label{tHiso}\\
     \H_\c: X_\c^{t,q} & \to \K^{t-1,q}_\c \label{Hiso}
   \end{align}
   are isomorphisms.  Moreover, we have the commutative diagram
\begin{equation}
  \begin{split}
  \xy
(0,0)*+{\tX^{t,q}_\c}="1";
(30,0)*+{\tT^{t-1,q}_{}}="2";
(0,-17)*+{X^{t,q}_\c}="3";
(30,-17)*+{\K^{t-1,q}_\c}="4";
{\ar^{\tH_\c} "1";"2"};
{\ar^{\H_\c} "3";"4"};
{\ar@{^{(}->} (0,-12);"1"};    
{\ar@{^{(}->} (30,-12)*{};"2"};
\endxy
\end{split} \label{CD3}
\end{equation}
In particular, we can take $(t,q) = (s+1,p)$ in the above.

The previous statements all remain true if $\tX_\c^{t,q}$ and $X^{t,q}_\c$ are replaced with $\tX_{(B',\Psi')}^{t,q}$ and $X_{(B',\Psi')}^{t,q}$, respectively, for $(B',\Psi') \in \fM^{s,p}$ in a sufficiently small $B^{s,p}(Y)$ neighborhood of $\c$.
\end{Proposition}

\Proof The map $\tH_\c: \tT^{t,q} \to \tT^{t-1,q}$ is surjective, by unique continuation, and by restricting to $\tX^{t,q}_\c$, we have eliminated the kernel.  Indeed, $r: \ker \tH^{t,q}_\c \to \tT^{t-1/q,q}_\Sigma$ is an isomorphism onto its image and $r(\ker \tH^{t,q}_\c) \cap \mathsf{L}_\c^{t,q} = 0$, whence $\ker \tH^{t,q}_\c \cap \tX^{t,q} = 0$.  This proves (\ref{tHiso}) is an isomorphism.  For (\ref{Hiso}), the same argument shows that (\ref{Hiso}) is injective. Indeed, $r_\Sigma: \ker (\H_\c|_{\C^{t,q}}) \to \T^{t-1/q,q}_\Sigma$ is injective by Corollary \ref{CorSur}(ii) and Remark \ref{Rem1/2}, and
$$r_\Sigma(\ker \H_\c|_{\C^{t,q}}) \cap r_\Sigma X^{t,q} = r_\Sigma(\ker \H_\c^{t,q}) \cap J_\Sigma r_\Sigma(\ker \H_\c^{t,q})) = 0$$
by Theorem \ref{ThmLinLag}(i).  It remains to show that (\ref{Hiso}) is surjective. We already know that $\H_\c: \T^{t,q} \to \K^{t-1,q}_\c$ is surjective by assumption.  So given any $(a,\phi) \in \T^{s,p}$, we need to find a $(b,\psi) \in X^{t,q}_\c$ such that $\H_\c(b,\psi) = \H_\c(a,\phi)$.  Without loss of generality, we can suppose $(a,\phi) \in \C^{t,q}$ by (\ref{JtC}) and since $J^{t,q}_{\c,t} \subseteq \ker \H^{t,q}_\c$. Since the condition $(b,\psi) \in X^{t,q}_\c$ imposes no restriction on the normal component of $b$ at the boundary, we only need to make sure that $r_\Sigma(b,\psi) \in J_\Sigma r_\Sigma(\ker \H^{t,q}_\c)$.  Since we have a decomposition
$$\T^{t-1/q,q}_\Sigma = r_\Sigma(\ker \H^{t,q}_\c) \oplus J_\Sigma r_\Sigma(\ker \H^{t,q}_\c),$$
we can write $r_\Sigma(a,\phi) = (a_0,\phi_0) + (a_1,\phi_1)$ with respect to the above decomposition.  Now let $(b,\psi) = (a,\phi) - P_\c (a_0,\phi_0)$, where $P_\c$ is the Poisson operator of $\H_\c^{t,q}$ with range equal $\ker (\H_\c|_{\C^{t,q}})$ as given by Theorem \ref{ThmLinLag}.  It follows that $(b,\psi) \in X^{t,q}_\c$, since $r_\Sigma(b,\psi) = (a_1,\phi_1) \in J_\Sigma r_\Sigma(\ker \H^{t,q}_\c)$ and that $(b,\psi) \in \C^{t,q}$ since both $(a,\phi)$ and $P_\c (a_0,\phi_0)$ belong to $\C^{t,q}$.  Thus, $(b,\psi) \in X^{t,q}_\c$ and we have $\H_\c^{t,q}(b,\psi) = \H_\c^{t,q}(a,\phi)$.  So (\ref{Hiso}) is surjective, hence an isomorphism.

The commutativity of the diagram (\ref{CD3}) now readily follows since (\ref{tHiso}) is an isomorphism which extends the isomorphism (\ref{Hiso}). Finally, for the last statement, we know that the space $\tX^{t,q}_\c$ varies continuously with $\c$ since the space $\mathsf{L}^{t-1/q,q}_\c$ varies continuously.  Since $\J^{t,q}_{\c,t} \subseteq \tX^{t,q}_\c$ for all $\c$, it follows that
$$X^{t,q}_\c = \Pi_{\C^{t,q}_\c}\{x \in \T^{t,q} : r(x) \in \mathsf{L}^{t-1/q,q}_\c\},$$
where $\Pi_{\C^{t,q}_\c}$ is the projection of $\T^{t,q}$ onto $\C^{t,q}_\c$ given by (\ref{JtC}).  From this, we see that $X^{t,q}_\c$ varies continuously since $\mathsf{L}^{t,q}_\c$ and $\J^{t,q}_{\c,t}$ vary continuously.  The continuity of $\tilde X^{t,q}_\c$ and $X^{t,q}_\c$ with respect to $\c$ implies the last statement.\End

The above proposition will be important when study the analytic properties of the spaces $\fM^{s,p}(Y)$ and $\M^{s,p}(Y)$ in the next section, where we will need to consider the inverse of the operator (\ref{Hiso}).  The point is that by restricting the domain of the Hessian operator $\H_\c$, it becomes invertible and its inverse smooths by one derivative in a certain range of topologies depending on the regularity of the configuration $\c$.  Thus, the inverse of $\H_\c$ behaves like a pseudodifferential operator of order $-1$ in this range, which is what one would formally expect since $\H_\c$ is a first order operator.  In particular, for $\c$ smooth, we have the following corollary:

\begin{Corollary}\label{CorInvertH}
  If $\c \in \fM$ is smooth, then for all $q \geq 2$ and $t > 1/q$, the maps
   \begin{align}
     \tH_\c: \tX^{t,q}_\c & \to \tT^{t-1,q},\\
     \H_\c: X_\c^{t,q} & \to \K^{t-1,q}_\c
   \end{align}
   are isomorphisms.
\end{Corollary}


\section{The Space of Monopoles}\label{SecMonopoles}

Having studied the linear theory of the Hessian operators $\tH_\c$ and $\H_\c$ in the previous section, we now study the space of Besov monopoles $\fM^{s,p}(Y,\s)$ and $\M^{s,p}(Y,\s)$ on $Y$.  Under suitable hypotheses, we show that these spaces are Banach manifolds and their local coordinate charts obey important analytic properties. Moreover, we show that smooth monopoles are dense in the spaces $\fM^{s,p}(Y,\s)$ and $\M^{s,p}(Y,\s)$, so that these Banach manifolds are Besov completions of the smooth monopole spaces $\fM(Y,\s)$ and $\M(Y,\s)$, respectively. These analytic properties are crucial for the analysis in \cite{N2}.\\

\noindent\textbf{Notation. }Recall that $\T^{s,p}_\c = T_\c\fC^{s,p}(Y)$ is the tangent space to a configuration $\c \in \fC^{s,p}(Y)$.  Since all these tangent spaces are identical, in the previous section we worked within one fixed copy and called it $\T^{s,p}$.  Now that we will work on the configuration space level, it is appropriate to keep track of the basepoint at times and we reintroduce this into our notation, though there really is no gain or loss of information by adding or dropping the basepoint from our notation.\\

Recall that we have fixed a $\spinc$ structure $\s$ from the start, which up to now, has not played any role in the analysis we have done.  We now consider the following assumption:
\begin{equation}
  c_1(\s) \textrm{ is non-torsion or } H^1(Y,\Sigma)=0. \label{assum}
\end{equation}
The following lemma is the fundamental reason we make the above assumption:

\begin{Lemma}\label{LemmaHtrans}
  Suppose (\ref{assum}) holds.  Let $s> \max(3/p,1/2)$.  Then for every $\c \in \fM^{s,p}(Y,\s)$, we have $\H_\c: \T_\c^{s,p} \to \K_\c^{s-1,p}$ is surjective.
\end{Lemma}

\Proof There are two cases $s > 1$ and $s \leq 1$.  We deal with the latter case, with the more regular case $s > 1$ being similar.  So for $s < 1$, there are two main steps.  First, we proceed as in the proof of Theorem \ref{ThmUCPsurj} to show that any element in the cokernel of $\H_\c: \T^{s,p}_\c \to \K^{s-1,p}_\c$ must be more regular, in fact, it must lie in $\K^{s+1,p}_\c$.  This follows because an element in the cokernel of $\H_\c$ satisfies an overdetermined elliptic boundary value problem, and thus we can bootstrap its regularity.  Once we have enough regularity, we can integrate by parts, which shows that any element $(b,\psi) \in \K^{s+1,p}_\c$ in the cokernel of $\H_\c^{s,p}$ must satisfy $\H_\c(b,\psi) = 0$ and $r_\Sigma(b,\psi) = 0$. From here, the second step is to apply the unique continuation theorem, Corollary \ref{CorUCP}, to deduce that the cokernel of $\H_\c^{s,p}$ is zero.

For the first step, by Lemma \ref{LemmaHessK}, we know that $\H_\c: \T_\c^{s,p} \to \K_\c^{s-1,p}$ has closed range and finite dimensional cokernel.  Let $(b,\psi) \in \T^{1-s,p'}_\c$, $p' = p/(p-1)$, be an element in the dual space of $\K_\c^{s-1,p}$ which annihilates $\im \H_\c^{s,p}$.  Indeed, we have that $\T_\c^{1-s,p'}$ is the dual space of $\T_\c^{s-1,p}$ by Theorem \ref{ThmFunMan}.  Next, we have the topological decomposition
\begin{equation}
  \T^{1-s,p'}_\c = \J^{1-s,p'}_{\c,t} \oplus \K^{1-s,p'}_{\c}. \label{Tp'}
\end{equation}
This follows from the decomposition (\ref{JtK}), since one can check that the map (\ref{PiJK}), by duality, is bounded on $\T^{1-s,p'}_\c$.  More precisely, by our choice of $s$, we have the multiplication maps
\begin{align*}
  B^{s,p}(Y) \times B^{s,p}(Y) \to B^{s,p}(Y)\\
  B^{s,p}(Y) \times B^{s-1,p}(Y) \to B^{s-1,p}(Y),
\end{align*}
which by duality means that the multiplications
\begin{align}
  B^{s,p}(Y) \times B^{-s,p'}(Y) &\to B^{-s,p'}(Y), \label{dualmult1}\\
  B^{s,p}(Y) \times B^{1-s,p'}(Y) &\to B^{1-s,p'}(Y). \label{dualmult2}
\end{align}
are also bounded.  Thus, repeating the proof of (\ref{JtK}) shows that there exists a bounded projection of $\T^{1-s,p'}_\c$ onto $\J^{1-s,p'}_\c$ through $\K^{1-s,p'}_{\c}$, for $\c \in \fC^{s,p}(Y)$.  This proves (\ref{Tp'}). Since $\J^{1-s,p'}_{\c,t}$ and $\K^{s-1,p}_{\c,t}$ annihilate each other, we see can choose our annihilating element $(b,\psi) \in \K^{1-s,p'}_\c$ since $\im (\H_\c^{s,p}) \subseteq \K^{s-1,p}_\c$.  Moreover, the fact that $(b,\psi)$ annihilates $\im (\H_\c^{s,p})$ also means that $\H_\c(b,\psi) = 0$ (weakly, i.e., as a distribution).  Altogether then, we see that we have the weak equation
\begin{equation}
  \tH_0(b,\psi) = (B - B_\rf, \Psi)\#(b,\psi). \label{EBPtH0}
\end{equation}
Everything now proceeds as in the bootstrapping argument in Theorem \ref{ThmUCPsurj}, but with modifications since the multiplication term is not smooth.  Because of the multiplication (\ref{dualmult2}), we have $(B - B_\rf, \Psi)\#(b,\psi) \in \tT^{1-s,p'}$.  By Theorem \ref{ThmEBP}(i), $r_\Sigma(b,\psi) \in \T^{1-s-1/p',p'}_\Sigma$ is well-defined.  Applying Green's formula to the symmetric operator $\H_\c$, we obtain for all $(a,\phi) \in \T$ that
\begin{align}
  0 &= (\H_\c(a,\phi),(b,\psi))_{L^2(Y)} - ((a,\phi),\H_\c(b,\psi))_{L^2(Y)} \nonumber\\
  &= -\omega(r_\Sigma(a,\phi), r_\Sigma(b,\psi)). \label{killboundary}
\end{align}
In the first line, we used that $(b,\psi)$ annihilates $\im (\H_\c)$ and $\H_\c(b,\psi) = 0$ (weakly).  In the second line, we use that $r_\Sigma(b,\psi) \in \T^{1-s-1/p',p'}_\Sigma$ is well-defined.  Since (\ref{killboundary}) holds for all $(a,\phi) \in \T$, we have $r_\Sigma(b,\psi) = 0$. This boundary condition together with $(\ref{EBPtH0})$ implies that we have an overdetermined elliptic boundary value problem (cf. Proposition \ref{PropInvertH}, we have $r(b,\psi) \in 0 \oplus B^{1-s-1/p',p'}\Omega^0(\Sigma; i\R) \oplus 0$). By Theorem \ref{ThmEBP}, this means we gain a derivative and so $(b,\psi) \in \T^{2-s,p'}_\c$.  This implies $(B-B_0,\Psi)\#(b,\psi)$ is more regular than an element of $\T^{1-s-1/p',p'}$, and we can elliptic bootstrap again.  We keep on boostrapping until we obtain $(b,\psi) \in \T^{s+1,p}_\c$, which is one derivative more regular than the maximum regularity of (\ref{EBPtH0}) since $(B,\Psi) \in \fM^{s,p}(Y)$.  Thus, $(b,\psi) \in \K^{s+1,p}_\c$ is now a strong solution to $\H_\c(b,\psi) = 0$.

We can now use Corollary \ref{CorUCP}, since $\K^{s+1,p}_\c \subset \K^{1,2}_\c$, as $p\geq 2$.  This theorem implies the following. Either $(b,\psi) = 0$, in which case the cokernel of $\H_\c: \T_\c^{s,p} \to \K_\c^{s-1,p}$ is zero, or else $\c = (B,0)$ and $\psi \equiv 0$, $b \in H^1(Y,\Sigma; i\R)$.  In the former case, our map $\H_\c: \T_\c^{s,p} \to \K_\c^{s-1,p}$ is surjective and we are done.  For the latter case, we apply assumption (\ref{assum}).  In case $c_1(\s)$ is non-torsion, $\det(\s)$ admits no flat connections, hence, we cannot have a reducible configuration $\c = (B,0)$ be a monopole, else $B^t$ would be a flat connection on $\det(\s)$.  In case $H^1(Y,\Sigma) = 0$, then we see $(b,\psi) = 0$ and the Hessian is surjective.  This proves the lemma.\End

\noindent\textbf{Assumption: }For the rest of this paper, we assume (\ref{assum}) holds.\\

So let us fix $Y$ and $\s$ satisfying (\ref{assum}), and write $\fM^{s,p} = \fM^{s,p}(Y,\s)$ and $\M^{s,p} = \fM^{s,p}(Y,\s)$ for short. The conclusion of the lemma guarantees that we have transversality for the monopole equations.  This implies the following theorem:

\begin{Theorem}\label{ThmMMan}
  For $s > \max(3/p,1/2)$, $\fM^{s,p}$ and $\M^{s,p}$ are closed submanifolds of $\fC^{s,p}(Y)$.
\end{Theorem}

\Proof For any smooth $\c \in \fC(Y)$, one can verify directly that $SW_3\c \in \K_\c$.\footnote{This is no coincidence.  On a closed-manifold $Y$, the Seiberg-Witten equations are the variational equations for the Chern-Simons-Dirac functional, see \cite{KM}.  In other words, $SW_3\c$ is the gradient of the Chern-Simons-Dirac functional $CSD$, i.e., the differential of $CSD$ at $\c$ satisfies $D_\c CSD(b,\psi) = (SW_3\c, (b,\psi))$ so that $SW_3\c$ vanishes precisely at the critical points of $CSD$.  When $\partial Y$ is nonempty, we still have $D_\c CSD(b,\psi) = (SW_3\c, (b,\psi)) = 0$ for $(b,\psi)$ vanishing on the boundary, in particular, for $(b,\psi) \in \J_{\c,t}$.  Since $CSD$ is invariant under the gauge group $\G_{\id,\partial}(Y)$, this means $(SW_3\c, (b,\psi)) = 0$ for all $(b,\psi) \in \J_{\c,t}$.  So $SW_3\c \in \K_\c$, the orthogonal complement.} Thus when $\c \in \fC^{s,p}(Y)$, we have $SW_3\c \in \K^{s-1,p}_\c$, since the map $\bd_\c^*: \T^{s-1,p}_\c \to B^{s-2,p}\Omega^0(Y;i\R)$ is still bounded by our choice of $s$.  Proceeding as in \cite[Chapter 12]{KM}, we can therefore think of $SW_3: \fC^{s,p}(Y) \to \K^{s-1,p}(Y)$ as a section of the Banach bundle $\K^{s-1,p}(Y) \to \fC^{s,p}(Y)$ (see Proposition \ref{PropBundle}). The previous lemma shows that $SW_3$ is transverse to the zero section.  More precisely, from Proposition \ref{PropBundle}, we have that $\K^{s-1,p}(Y) \to \fC^{s,p}(Y)$ is Banach bundle complementary to the bundle $\J^{s-1,p}_t(Y) \to \fC^{s,p}(Y)$, which means that for any configuration $\co \in \fC^{s,p}(Y)$, there exists a neighborhood $\frak{U}$ of $\co$ in $\fC^{s,p}(Y)$ such that
\begin{equation}
  \Pi_{\K^{s-1,p}_\co}: \K^{s-1,p}_\c \to \K^{s-1,p}_\co \label{PiK}
\end{equation}
is an isomorphism for all $\c \in \frak{U}$.  Here, $\Pi_{\K^{s-1,p}_\co}: \T^{s-1,p}_\c \to \K^{s-1,p}_{\co}$ is the projection through $\J^{s-1,p}_{\co,t}$ given by (\ref{3.1-PiK}).  Thus, if $SW_3\co = 0$, we consider the map
\begin{equation}
  f = \Pi_{\K^{s-1,p}_\co}SW_3: \frak{U} \to \K^{s-1,p}_\co \label{fSW3}
\end{equation}
Then $f\c = 0$ if and only if $SW_3\c = 0$, and at such a monopole, we have
\begin{equation}
  D_\c f = \Pi_{\K^{s-1,p}_\co}\H_\c^{s,p}: \T^{s,p}_\c \to \K^{s-1,p}_\co.
\end{equation}
By Lemma \ref{LemmaHtrans}, $\H_\c: \T^{s,p}_\c \to \K^{s-1,p}_\c$ is surjective, and so since (\ref{PiK}) is an isomorphism, this means $D_\c f$ is surjective for all $\c \in \frak{U}$.  Thus, we can apply the implicit function theorem to conclude that $f^{-1}(0)$ is a submanifold of $\fC^{s,p}(Y)$.  Since we can apply the preceding local model near every monopole, it follows that $\fM^{s,p} = SW_3^{-1}(0) \subset \fC^{s,p}(Y)$ is globally a smooth Banach submanifold.  Lemma \ref{PiCdecomp} implies that we have the product decomposition
\begin{equation}
  \fM^{s,p} = \G^{s+1,p}_{\id,\partial}(Y) \times \M^{s,p}. \label{fMM}
\end{equation}
Thus $\M^{s,p}$ is also a submanifold of $\fC^{s,p}(Y)$, since $\G^{s+1,p}_{\id,\partial}(Y)$ is a smooth Banach Lie group by Lemma \ref{Ggroup}.  The closedness of $\fM^{s,p}$ and $\M^{s,p}$ readily follows from the fact that these two spaces are defined as the zero set of equations.\En

\begin{Rem}\label{RemLball}
  Note that we can take the open neighborhood $\frak{U} \subset \fC^{s,p}(Y)$ of $\co$ to contain a ball in the $L^2(Y)$ topology (so that $\frak{U}$ is a very large open subset of $\fC^{s,p}(Y)$).  Indeed, this is because $\K^{s-1,p}_{\cx[1]}$ and $\J^{s-1,p}_{\co,t}$ are complementary for any $\cx[1]$ in a sufficiently small $L^2(Y)$ neighborhood of $\co$, and so the map (\ref{PiK}) is an isomorphism for $\c = \cx[1]$.  To show this, it suffices to show that
  \begin{equation}
    \K^{s-1,p}_{\cx[1]} \cap \J^{s-1,p}_{\co,t} = 0. \label{eq3.3:KcapJ=0}
  \end{equation}
  Indeed, this will show that (\ref{PiK}) injective.  However, it must also be an isomorphism, since $\K^{s-1,p}_\c$ varies continuously with $\c \in \fC^{s,p}(Y)$ as a consequence of Proposition \ref{PropBundle}.  Namely, since (\ref{PiK}) is an isomorphism for $\c = \co$, then if it is injective for all $\c = (B(t),\Psi(t))$ along a path in $\fC^{s,p}(Y)$ joining $\co$ to $\cx[1]$, then it must also be an isomorphism for all such $\c$.

  We now show (\ref{eq3.3:KcapJ=0}).  Note that an element of $\K^{s-1,p}_{\cx[1]} \cap \J^{s-1,p}_{\co,t}$ is determined by
  a  $\xi \in B^{s,p}\Omega^0(Y;i\R)$ that solves
  \begin{align}
    \Delta \xi + \Re(\Psi_1,\Psi_0)\xi &= 0 \label{eq3.3:xi1}\\
    \xi|_\Sigma &= 0. \label{eq3.3:xi2}
  \end{align}
  Using elliptic regularity for the Dirichlet Laplacian, we bootstrap the regularity of $\xi$ to obtain $\xi \in B^{2,2}\Omega^2(Y;i\R)$. Writing $\Delta + (\Psi_1,\Psi_0) = \Delta + |\Psi_0|^2 + \Re(\Psi_1-\Psi_0,\Psi_0)$, we see that the operator $\Delta + (\Psi_1,\Psi_0)$ is a perturbation of the operator
  $$\Delta + |\Psi_0|^2: B^{2,2}\Omega^0_t(Y;i\R) \to L^2\Omega^0(Y;i\R),$$
  whose domain $B^{2,2}\Omega^0_t(Y;i\R)$ consists of those $\alpha \in B^{2,2}\Omega^0(Y;i\R)$ such that $\alpha|_\Sigma = 0$.  We showed that this latter operator is invertible in the proof of Lemma \ref{LemmaTdecomp}.  It follows that if the multiplication operator $\Re(\Psi_1-\Psi_0,\Psi_0)$ has small enough norm, as a map from $B^{2,2}(Y)$ to $L^2(Y)$, then the operator $\Delta + \Re(\Psi_1,\Psi_0)$ remains invertible and the only solution to (\ref{eq3.3:xi1})--(\ref{eq3.3:xi2}) is $\xi = 0$.
We have
  \begin{align}
    \|\Re(\Psi_1-\Psi_0,\Psi_0)\alpha\|_{L^2(Y)} & \leq \|\Psi_1-\Psi_0\|_{L^2(Y)}\|\Psi_0\|_{L^\infty(Y)}\|\alpha\|_{L^\infty(Y)} \label{eq:L2small}\\
     & \leq C\|\Psi_1-\Psi_0\|_{L^2(Y)}\|\Psi_0\|_{B^{s,p}(Y)}\|\alpha\|_{B^{2,2}(Y)}. \nonumber
  \end{align}
since both $B^{s,p}(Y)$ and $B^{2,2}(Y)$ embed into $L^\infty(Y)$.  Hence, if $\|\Psi_1-\Psi_0\|_{L^2(Y)}$ is sufficiently small, we see that the only solution to (\ref{eq3.3:xi1})--(\ref{eq3.3:xi2}) is $\xi = 0$, which establishes (\ref{eq3.3:KcapJ=0}).
\end{Rem}

Theorem \ref{ThmMMan} proves the first part of our main theorem.  However, to better understand the analytic properties of these monopole spaces, we want to construct explicit charts for our manifolds $\fM^{s,p}$ and $\M^{s,p}$.  Furthermore, we want to show that smooth monopoles are dense in these spaces.  These properties are not only of interest in their own right but will be essential in \cite{N2}.

In a neighborhood of $\c \in \M^{s,p}$, the Banach manifolds $\fM^{s,p}$ and $\M^{s,p}$ are modeled on their tangent spaces at $\c$, namely $\ker \H^{s,p}_\c$ and $\ker(\H_\c|_{\C^{s,p}}) = \ker (\tH_\c|_{\T^{s,p}})$, respectively.  Moreover, the tangent space to our manifolds at $\c$ are the range of operators which are ``nearly pseudodifferential".  Indeed, in the previous section, we constructed a Poisson operator $P_\c$ whose range is $\ker (\tH_\c|_{\T^{s,p}})$.  Since this operator is constructed from the Calderon projection $P^+_\c$ and the Poisson operator $\tilde P^+_\c$ for the augmented Hessian $\tH_\c$, both of which differ from pseudodifferential operators by a compact operator, it is in this sense that $P_\c$ is close to being pseudodifferential.

Let $\cx[1],\co \in \fC^{s,p}(Y)$ and write $(b,\psi) = (B_1-B_0,\Psi_1-\Psi_0)$.  Then we have the difference equation
\begin{align}
SW_3\cx[1]-SW_3\co &= \H_{\co}(b,\psi) + (\rho^{-1}(\psi\psi^*)_0,\rho(b)\psi),  \label{SWlin} 
\end{align}
which reflects the fact that $SW_3$ is a quadratic map.  The linear part, is of course, given by the Hessian, and its quadratic part is just a pointwise multiplication map.  Thus, we define the bilinear map \label{p:q}
\begin{align}
  \q: \T \times \T & \to \T \nonumber \\
  \q((b_1,\psi_1),(b_2,\psi_2)) & = \left(\rho^{-1}(\psi_1\psi_2^*)_0, \frac{1}{2}\big(\rho(b_1)\psi_2 + \rho(b_2)\psi_1\big)\right) \label{def:q}
\end{align}
which as a quadratic function enters into the Seiberg-Witten map via (\ref{SWlin}).  The map $\q$ extends to function space completions as governed by the multiplication theorems.  Observe that $q$ is a bounded map on $\T^{s,p}$ since $B^{s,p}(Y)$ is an algebra.  This is key, because then the Seiberg-Witten map $SW_3$ is the sum of a first order differential operator and a zeroth order operator, and using Proposition \ref{PropInvertH}, we have elliptic regularity for the linear part of the operator on suitable domains.

From these observations, we can prove the following important lemma which we will need to show that smooth monopoles are dense in $\fM^{s,p}$.

\begin{Lemma}\label{LemmaSmoothMan}
  Let $s > \max(3/p,1/2)$. Let $\co \in \M^{s,p}$. Then $B^{s+1,p}(Y)$ configurations are dense in the affine space $\co + T_\co\M^{s,p}$.
\end{Lemma}

\Proof Pick any smooth $(B_1,\Psi_1) \in \fC_C(Y)$ in Coulomb-gauge with respect to $B_\rf$.  Let $(b,\psi) = (B_1-B_0,\Psi_1-\Psi_0)$. Then from (\ref{SWlin}) together with the Coulomb-gauge condition, we have
\begin{equation}
  \tH_{\co}^{s,p}(b,\psi) = SW_3(B_1,\Psi_1) - \q((b,\psi),(b,\psi)), \label{eq:4.3}
\end{equation}
where on the right-hand side the first term is smooth and the second term is in $\T^{s,p}$. Applying Proposition \ref{PropInvertH} with $(t,q) = (s+1,p)$, we see that $(b,\psi) \in (b',\psi') + \ker (\tH_{\co}|_{\T^{s,p}})$ for some $(b',\psi') \in X^{s+1,p}_\co \subseteq \T^{s+1,p}$.  In other words, if we invert $\tH_{\co}^{s,p}$ in (\ref{eq:4.3}), we find that $(b,\psi)$ is equal to a smoother element $(b',\psi')$, modulo an element of the kernel of $\tH_{\co}^{s,p}|_{\T^{s,p}}$.  It remains to show that $B^{s+1,p}(Y)$ configurations are dense in the latter space.  First, we have $B^{s+1,p}(Y)$ configurations are dense in $\ker\tH_{\co}^{s,p} \subset \tT^{s,p}$ by Corollary \ref{CorKerComp} and Lemma \ref{LemmaKerH}. Similarly, $B^{s+1,p}(Y)$ configurations are dense in $\Gamma_0$, the subspace given by (\ref{tHCG}).  This follows from the construction of $\Gamma_0$.  First, we have $\Gamma_0$ is a graph of the map $\Theta_0$, which is defined over $\ker \Delta \subseteq B^{s,p}\Omega^0(Y; i\R)$, and smooth configurations are dense in $\ker \Delta$ by Corollary \ref{CorKerComp}.  We now apply Proposition \ref{PropInvertH} with $(t,q) = (s+1,p)$, since the map $\Theta_0$ is defined by inverting the Hessian.  Altogether, we see that $B^{s+1,p}(Y)$ configurations are dense in $\Gamma_0$. Because of the decomposition (\ref{tHCG}), it now follows from the density of $B^{s+1,p}(Y)$ configurations in $\ker \tH_\co^{s,p}$ and $\Gamma_0$ that $B^{s+1,p}(Y)$ configurations are dense in $\ker (\tH_\co|_{\T^{s,p}}) = T_\co\M^{s,p}$.

Altogether, we have shown that $\co + T_\co\M^{s,p} = (B_1,\Psi_1) + (b',\psi') + T_\co\M^{s,p}$, where $(B_1,\Psi_1)$ is smooth, $(b',\psi') \in \T^{s+1,p}$, and $B^{s+1,p}(Y)$ configurations are dense in $T_\co\M^{s,p}$.  This proves the lemma.\End

From Theorem \ref{ThmLinLag}, given $\co \in \M^{s,p}$, we have a projection $\pi_\co = P_\co r_\Sigma: \T_\co^{s,p} \to T_\co\fM^{s,p}$ onto the tangent space $T_\co \fM^{s,p}$ for any $\co \in \fM^{s,p}$.  Thus, locally $\fM^{s,p}$ is the graph of a map from $T_\co \fM^{s,p}$ to any complementary subspace in $\T^{s,p}_\co$.  We wish to describe the analytic properties of this local graph model in more detail.  First, we record the following simple lemma which describes for us natural complementary subspaces for $T_\co\fM^{s,p}$.

\begin{Lemma}\label{LemmaXdecomp}
  Let $s > \max(3/p,1/2)$. Given any $\co \in \fM^{s,p}$, we have the direct sum decomposition
  \begin{equation}
    \T^{s,p}_\co = T_\co\fM^{s,p} \oplus X^{s,p}_{\c}  \label{Xdecomp1}
  \end{equation}
  for any $\c \in \fM^{s,p}$ sufficiently $B^{s,p}(Y)$ close to $\co$, where $X^{s,p}_\c$ is defined as in (\ref{Xtq}).
\end{Lemma}

\Proof By Lemma \ref{LemmaHtrans}, $\H_\co: \T^{s,p}_\co \to \K^{s-1,p}_\co$ is surjective.  Thus, (\ref{Xdecomp1}) follows readily from $T_\co\fM^{s,p} = \ker \H_\co^{s,p}$ and $\H_\co: X^{s,p}_\c \to \K^{s-1,p}_\co$ being an isomorphism by Proposition \ref{PropInvertH}. Note also that $X_\c^{s,p}$ is the kernel of the projection $\pi_\c: \T^{s,p} \to T_\c\M^{s,p}$.\End

Using any one of above complementary subspaces for $T_\co\fM^{s,p}$ (we will always use $X^{s,p}_\co$ for simplicity), we can describe the Banach manifold $\fM^{s,p}$ locally as follows.  In the proof of Theorem \ref{ThmMMan}, we introduced the local defining function $f$ in (\ref{fSW3}) on a neighborhood $\frak{U} \subset \fC^{s,p}(Y)$ so that $\fM^{s,p} \cap \frak{U} = f^{-1}(0)$.  In other words, we used the implicit function theorem for $f$ to obtain $\fM^{s,p}$.  On the other hand, we can describe $\fM^{s,p}$ in an equivalent way using the inverse function theorem,  as in the framework of Theorem \ref{ThmIFT}, whereby $\fM^{s,p}$ is given locally by the preimage of an open set under diffeormophism rather than the preimage of a regular value of a surjective map.  This means we need to construct a \textit{local straightening map} as in Definition \ref{DefLDF}.  Following the same ansatz in Theorem \ref{ThmIFT}, we have the following:

\begin{Lemma}\label{LemmaF}
  Let $\co \in \fM^{s,p}$, and let $X = \T_\co^{s,p}$, $X_0 = T_{\co}\fM^{s,p}$, and $X_1 = X_\co^{s,p}$. We have $X = X_0 \oplus X_1$ and define the map
\begin{align}
  F_\co: X_0 \oplus X_1 & \to X \nonumber \\
  x = (x_0,x_1) & \mapsto x_0 + (\H_{\co}|_{X_1})^{-1}\Pi_{\K^{s-1,p}_\co}SW_3\big(\co + x\big).\label{eq:Fmap}
\end{align}
\begin{enumerate}
    \item Then $F_\co(0)=0$, $D_0F_\co = 0$, and $F_\co$ is a local diffeomorphism in a $B^{s,p}(Y)$ neighborhood of $0$.
    \item There exists an open set $V \subset X$ containing $0$ such that for any $x \in V$, we have $\co + x \in \fM^{s,p}$ if and only if $F_\co(x) \in X_0$. We can choose $V$ to contain an $L^2(Y)$ ball, i.e., there exists a $\delta > 0$ such that
    $$V \supseteq \{x \in X: \|x\|_{L^2(Y)} < \delta\}.$$
    Furthermore, we can choose $\delta = \delta\co$ uniformly for all $\co$ in a sufficiently small $L^\infty(Y)$ neighborhood of any configuration in $\fM^{s,p}$.
    \item If  $\co \in \M^{s,p}$ then for any $x \in V \cap \cC^{s,p}_\co$, we have $\co + x \in \M^{s,p}$ if and only if $F_\co(x) \in X_0 \cap \cC^{s,p}_\co$.
  \end{enumerate}
\end{Lemma}

\Proof (i) We have $F_\co(0) = 0$ since $SW_3\co=0$.  Furthermore, the differential of $F_\co$ at $0$ is the identity map by construction; more explicitly,
\begin{align*}
  (D_0F_\co)(x) &= x_0 + (\H_{\co}|_{X_1})^{-1}\Pi_{\K^{s-1,p}_\co}\H_{\co}(x)\\
   &= x_0 + (\H_{\co}|_{X_1})^{-1}\H_{\co}(x) = x_0 + x_1.
\end{align*}
So by the inverse function theorem, $F_\co$ is a local diffeomorphism in a $B^{s,p}(Y)$ neighborhood of $0$.

(ii) Observe that $F_\co(x) \in X_0$ if and only if the second term of (\ref{eq:Fmap}), which lies in $X_1$, vanishes.  Let $(B,\Psi) = \co+x$. Then for $x$ in a small $L^2(Y)$ neighborhood of $0 \in X$, call it $V$, we know by Remark \ref{RemLball} that (\ref{PiK}) is an isomorphism.  Since $SW_3\big(\co+x\big) \in \K^{s-1,p}_\c$ and $(\H_{\co}|_{X_1})^{-1}$ is an isomorphism, it follows that the second term of (\ref{eq:Fmap}) vanishes if and only if $SW_3\big(\co+x\big)$ vanishes, i.e., if and only if $\co + x \in \fM^{s,p}$. Equation (\ref{eq:L2small}) shows that the size of this $L^2(Y)$ ball depends only on $\|\Psi_0\|_{L^\infty(Y)}$, and this implies the continuity statement for $\delta$.

(iii) Since $X_1 \subset \C^{s,p}_{\co}$ by (\ref{Xtq}), we have that $F_\co(x) \in \C^{s,p}_\co$ if and only if $x_0 \in \C^{s,p}_{\co}$.  Then (iii) now follows from the previous steps via intersection with $\C^{s,p}_\co$.\End

Thus, the map $F_\co$ in the above lemma is a local straightening map for $\fM^{s,p}$ (where we translate by the basepoint $\co$ so that we can regard $\fM^{s,p}$ as living inside the Banach space $\T^{s,p}_\co$) such that its restriction to $\cC^{s,p}_\co$ yields a local straightening map for $\M^{s,p}$ if $\co \in \M^{s,p}$.  In Theorem \ref{ThmMcharts}, we will show, in the precise sense of Definition \ref{DefLDF} that $F_\co$ is a local straightening map for $\M^{s,p}$ within a ``large" neighborhood of $\co$, where large means that the open set contains a ball in a topology weaker than the ambient $B^{s,p}(Y)$ topology. First, we need another important lemma, which allows us to redefine $F_\co$ on weaker function spaces:

\begin{Lemma}\label{LemmaQ}
  Let $\co \in \M^{s,p}$ for $s > \max(3/p,1/2)$. \label{p:Fco}
  \begin{enumerate}
    \item If $x \in \T^{s,p}_\co$ , then we can write $F_\co(x)$ as
    \begin{align}
   F_\co(x)  &= x + (\H_{\co}|_{X^{s+1,p}_\co})^{-1}\Pi_{\K^{s,p}_\co}\q(x,x), \label{eq:q} \\
             &=: x + Q_\co(x,x), \label{eq:Q}
             \end{align}
            where $\q$ is the quadratic multiplication map given by (\ref{def:q}).
   \item The map $(\H_{\co}|_{X^{s+1,p}_\co})^{-1}\Pi_{\K^{s,p}_\co}: \T^{s,p} \to \T^{s+1,p}$ extends to a bounded map
   \begin{equation}
     (\H_{\co}|_{X^{s+1,p}_\co})^{-1}\Pi_{\K^{s,p}_\co}: L^q\T \to H^{1,q}\T \label{eq4:HKmap}
   \end{equation}
   for any $1 < q < \infty$.
   \item Let $3 \leq q \leq \infty$.  For $x \in L^q\T_\co$, define $F_\co(x)$ by (\ref{eq:Q}).  Then $F_\co: L^q\T_\co \to L^q\T_\co$ is a local diffeomorphism in a $L^q(Y)$ neighborhood of $0$.
  \end{enumerate}
\end{Lemma}

\Proof (i) With $x = (x_0,x_1)$ as in Lemma \ref{LemmaF}, we have
\begin{align}
  F(x) & = x_0 + (\H_{\co}|_{X_\co^{s,p}})^{-1}\Pi_{\K^{s-1,p}_\co}SW_3\big(\co + x\big)  \nonumber\\
  &= x_0 + (\H_{\co}|_{X_\co^{s,p}})^{-1}\Pi_{\K^{s-1,p}_\co}\big(\H_\co(x) + \q(x,x)\big) \nonumber \\
  &= x_0 + x_1 + (\H_{\co}|_{X_\co^{s,p}})^{-1}\Pi_{\K^{s-1,p}_\co}\q(x,x). \label{eq4:LemmaQ}
\end{align}
Next, since $B^{s,p}(Y)$ is an algebra, then $\q(x,x) \in \T^{s,p}$.  It follows that in (\ref{eq4:LemmaQ}), we may replace
$\Pi_{\K^{s-1,p}_\co}$ with $\Pi_{\K^{s,p}_\co}$.  From Proposition \ref{PropInvertH}, we know that we have isomorphisms
\begin{align*}
  \H_\co: X_\co^{s,p} & \to \K^{s-1,p}_\co\\
  \H_\co: X_\co^{s+1,p} & \to \K^{s,p}_\co,
\end{align*}
from which it follows that if $y \in \K^{s,p}_\co$, then $\H_\co^{-1}(y) \in X_\co^{s+1,p}$.  The decomposition (\ref{eq:q}) now follows.

(ii) First, we note that Lemma \ref{LemmaTdecomp} extends to Sobolev spaces, since its proof, which involves studying elliptic boundary value problems, carries over verbatim to Sobolev spaces (see Appendix \ref{AppEBP}) so long as the requisite function space multiplication works out.  In this case, we want $\Pi_{\K^{s,p}_\co}$ to yield a bounded map on $L^q\T$, in which case, the bounded multiplications that we want are the boundedness of
\begin{align}
  B^{s,p}(Y) \times H^{1,q}(Y) & \to L^q(Y) \\
  B^{s,p}(Y) \times L^q(Y) & \to H^{-1,q}(Y),
\end{align}
cf. (\ref{sec3:Bmult1}) and (\ref{sec3:Bmult2}).  However, these are straightforward, because we have the embedding $B^{s,p}(Y) \hookrightarrow L^\infty(Y)$, and we have the obvious bounded multiplication $L^\infty(Y) \times L^q(Y) \to L^q(Y)$, which therefore trivially imply the above multiplications.

From this, it remains to show that $(\H_{\co}|_{X_\co^{s+1,p}})^{-1}$ extends to a bounded map $L^q\T \to H^{1,q}\T$.  However, the exact same considerations show that this is the case due to the boundedness of the above multiplication maps.

(iii) We have a bounded multiplication map $L^q(Y) \times L^q(Y) \to L^{q/2}(Y)$, and for $q \geq 3$, we have the Sobolev embedding $H^{1,q/2}(Y) \hookrightarrow L^q(Y)$.  Hence (ii) implies that the map $F_\co$ is bounded on $L^q\T$ for $q \geq 3$. Since $F_\co(0) = 0$ and $D_0F_\co = \id$, the inverse function theorem implies $F_\co$ is a local diffeomorphism in a $L^q(Y)$ neighborhood of $0$.\End

Thus, from now on, we may work with the expression (\ref{eq:Q}) for $F_\co$ since it coincides with (\ref{eq:Fmap}) when the latter is well-defined.  

Given the local straightening map $F_\co$ and the various properties it obeys above, we now import the abstract point of view in Appendix \ref{AppIFT} into our particular situation to construct charts for $\fM^{s,p}$.  This gives us the following picture for a neighborhood of the monopole space $\fM^{s,p}$.  At any $\co \in \fM^{s,p}$, letting $X_1 = X^{s,p}_\co$ be a complement of $T_\co\fM^{s,p}$ as in Lemma \ref{LemmaXdecomp}, then near $\co$, the space $\fM^{s,p}(Y)$ is locally the graph of a map, which we denote by $E^1_\co$, from a neighborhood of $0$ in $T_\co\fM^{s,p}$ to $X_1$.  The local chart map we obtain in this way for $\fM^{s,p}(Y)$ is precisely the induced chart map of the local straightening map $F_\co$ above, in the sense of Definition \ref{DefChart}.  In addition, we show that the map $E^1_\co$ is smoothing, due to the fact that the lower order term $Q_\co$ occurring in $F_\co$, as defined in (\ref{eq:Q}), is smoothing, i.e., it maps $\T^{s,p}_\co$ to $\T^{s+1,p}_\co$.  Moreover, for any $q > 3$, we show that $F_\co$ is a local straightening map in some $L^q(Y)$ neighborhood of $\co$.  Consequently, the induced chart maps we obtain yield charts for $L^q(Y)$ neighborhoods of $\fM^{s,p}$, which are large neighborhoods when viewed within the ambient $B^{s,p}(Y)$ topology.  This latter property will be very important in \cite{N2}, and it is the analog of how the local Coulomb slice theorems for nonabelian gauge theory allow for gauge fixing within large neighborhoods (i.e., neighborhoods defined with respect to a weak norm) of a reference connection (see e.g. \cite[Theorem 8.1]{WeUC}).\footnote{Such gauge fixing properties are important for issues related to compactness, since in proving a compactness theorem, one considers a sequence of configurations that are bounded in some norm, hence strongly convergent along a subsequence but with respect to a weaker norm.  If one wants to gauge fix the elements in the convergent subsequence, one therefore needs a gauge fixing theorem on balls defined with respect to the weaker norm.}

We have the following theorem: 

\begin{Theorem}\label{ThmMcharts}Assume $s > \max(3/p,1/2)$.
\begin{enumerate}
    \item Let $\co \in \fM^{s,p}$ and $X_1 = X^{s,p}_\co$ be a complement of $T_\co\fM^{s,p}$ in $\T^{s,p}_\co$.  Then there exists a neighborhood $U$ of $0 \in T_{\co}\fM^{s,p}$ and a map $E^1_{\co}: U \to X_1$ such that the map \label{p:Eco}
        \begin{align}
        E_{\co}: U & \to \fM^{s,p} \nonumber \\
        x & \mapsto \co + x + E^1_{\co}(x)
        \end{align}
        is a diffeomorphism of $U$ onto an open neighborhood of $\co$ in $\fM^{s,p}$. We have $E^1_\co(0) = 0$, $D_0E^1_\co = 0$, and furthermore, the map $E^1_\co$ smooths by one derivative, i.e., $E^1_\co(x) \in \T^{s+1,p}_\co$ for all $x \in U$. \label{p:E1co}
    \item Let $q > 3$. We can choose $U$ such that both $U$ and its image $E_{\co}(U)$ contain $L^q(Y)$ neighborhoods, i.e., there exists a $\delta > 0$, depending on $\co$, such that
        \begin{align*}
        U & \supseteq \{x \in T_\co\fM^{s,p}:\|x\|_{L^q(Y)} < \delta\}\\
        E_\co(U) & \supseteq \{\c \in \M^{s,p}: \|\c - \co\|_{L^q(Y)} < \delta\}.
        \end{align*}
        The constant $\delta$ can be chosen uniformly in $\co$, for all $\co$ in a sufficiently small $L^\infty(Y)$ neighborhood of any configuration in $\fM^{s,p}$.
    \item If $\co \in \M^{s,p}$, then the map $E_{\co}$ restricted to $U \cap \C^{s,p}_\co$ is a diffeomorphism onto a neighborhood of $\co$ in $\M^{s,p}$.
    \item The smooth monopole spaces $\fM$ and $\M$ are dense in $\fM^{s,p}$ and $\M^{s,p}$, respectively.
\end{enumerate}
\end{Theorem}


\begin{figure}
\centering
\includegraphics{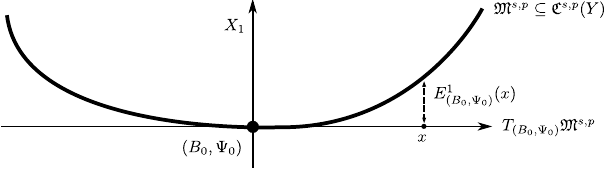}
\caption{A local chart map for $\frak{M}^{s,p}$ at $(B_0,\Phi_0)$.}
\label{Figure1}
\end{figure}

\Proof (i-ii) Given $q > 3$, consider $F_\co$ as defined by (\ref{eq:Q}).  By Lemma \ref{LemmaQ}(iii), $F_\co$ has a local inverse $F_\co^{-1}$ defined in an $L^q$ neighborhood, call it $V_q$, of $0 \in L^q\T_\co$.  First, we show that $F_\co^{-1}$ (equivalently, $F_\co$) is regularity preserving, namely, that $F_\co^{-1}(x) \in \T^{s,p}_\co$ if and only if $x \in \T^{s,p}_\co \cap V_q$.
In one direction, suppose $F_\co^{-1}(x)$ belongs to $\T^{s,p}_\co$. Then since $F_\co$ is continuous on $\T^{s,p}_\co$, then $x = F_\co(F_\co^{-1}(x)) \in \T^{s,p}_\co$.  In the other direction, if $x \in\T^{s,p}_\co$, we apply (\ref{eq:Q}) to obtain
$$F_\co^{-1}(x) = x - Q_\co(F_\co^{-1}(x), F_\co^{-1}(x)).$$
A priori, we only know that $F_\co^{-1}(x) \in L^q\T_\co$.  However, in the above, we have $x \in \T^{s,p}_\co$ and $Q_\co(x) \in H^{1,q/2}\T$ by Lemma \ref{LemmaQ}.  When $q > 3$, then $Q_\co$ always gains for us regularity, and so we can bootstrap the regularity of $F_\co^{-1}(x)$ until it has the same regularity as $x$.  Thus, this shows that $x \in \T^{s,p}$ if and only if $F_\co^{-1}(x) \in \T^{s,p}$.

Shrink $V_q$ if necessary so that $V_q \cap \T^{s,p} \subset V_2$, where $V_2$ is defined to be the open set in Lemma \ref{LemmaF}(ii).  This is possible since $V_2$ contains an $L^2(Y)$ ball and $q > 2$.  Then if we let $V = V_q \cap \T^{s,p}$, then $V$ satisfies the key property of Lemma \ref{LemmaF}(ii), namely if $x \in V$, then $\co + x \in \M^{s,p}$ if and only if $F(x) \in X_0$.  The key step we have done here is that we have shown that $F_\co^{-1}$ is well-defined on the open set $V$, so that $F_\co$ becomes a local straightening map for $\fM^{s,p}$ within the neighborhood $V$ of $\co \in \fM^{s,p}$.  Indeed, with just Lemma \ref{LemmaF}, we would only know that $F_\co$ is a straightening map for a small $B^{s,p}(Y)$ neighborhood of $\co \in \fM^{s,p}$, which is what we get when we apply the inverse function theorem for $F_\co$ as a map on $\T^{s,p}$.  Here, by rewriting $F_\co$ in Lemma \ref{LemmaQ} in a way that makes sense on $L^q$, we get an $L^q$ open set on which we have the inverse $F_\co^{-1}$.  The smoothing property of $Q_\co$ allows us to conclude the regularity preservation property of $F_\co^{-1}$, i.e., it preserves the $B^{s,p}(Y)$ topology, so that altogether, the map $F_\co$ is a straightening map for $\fM^{s,p}$ on the large open set $V \subset \T^{s,p}$.

Once we have the local straightening map $F_\co$, the construction of induced chart maps for $\fM^{s,p}$ now follows from the general picture described in the appendix.  Letting $U = F(V) \cap X_0$, the map $E_\co$ is given by
\begin{equation}
  E_\co(x) = \co + F_\co^{-1}(x), \qquad x \in U. \label{eq4:EfromF}
\end{equation}
The map $E_\co^1(x)$ is just the nonlinear part of $E_\co(x)$, and it is given by
\begin{align}
  E_\co^1(x) &= F_\co^{-1}(x) - x \label{eq:E1-1}\\
  &= -Q_\co(F_\co^{-1}(x),F_\co^{-1}(x)), \qquad x \in U. \label{eq:E1-2}
\end{align}
The smoothing property of $E_\co^1$ now readily follows from the smoothing property of $Q_\co$.  By construction, $U$ contains an $L^q(Y)$ ball since $V$ does.  This implies $E_\co(U)$ contains an $L^q(Y)$ neighborhood of $\co \in \fM^{s,p}(Y)$, since $\fM \cap (\co + V) = E_\co(U)$.

Finally, the local uniform dependence of $\delta$ can be seen as follows.  First, the constant $\delta$ of Lemma \ref{LemmaF} can be chosen uniformly for $\co$ in a small $L^\infty(Y)$ neighborhood of any configuration in $\fM^{s,p}$.  Next, the map $F_\co: L^q\T \to \L^q\T$ varies continuously as $\co$ varies in the $L^\infty(Y)$ topology.  It follows from the construction of $V$ that we can find a fixed $\delta$ such that $V$ contains a $\delta$-ball in the $L^q(Y)$ topology as $\co$ varies inside a small $L^\infty(Y)$ ball. We have now established all statements in (i-ii).

(iii) This follows from the above and Lemma \ref{LemmaF}(iii).

(iv) By Lemma \ref{LemmaSmoothMan} and the smoothing property of $E^1_\co$, we have that $\M^{s+1,p}$ is dense in $\M^{s,p}$.  Iterating this in $s$, we see that $\M$ is dense in $\M^{s,p}$.  Since smooth gauge transformations are dense in the space of gauge transformations, it follows from the decomposition $\fM^{s,p} = \G^{s+1,p}_{\id, \partial}(Y) \times \M^{s,p}$ that $\fM$ is dense in $\fM^{s,p}$ as well.\End

Retracing through the steps in the proof of Theorem \ref{ThmMcharts}, one sees that the chart maps for $\fM^{s,p}$ define bounded maps on weaker function spaces. This allows us to extend these chart maps to $L^q(Y)$ balls inside the closures of the tangent spaces to $\fM^{s,p}$ in weaker topologies.  This yields for us the following important corollary:

\begin{Corollary}\label{CorMCharts}
Let $\co \in \fM^{s,p}$.  Let $1/p \leq t \leq s$ and pick $q \geq 3$ according to the following: for $t = 1/p$, set $q = 3$;  else for $t > 1/p$, choose $q > 3$ such that $B^{t,p} \hookrightarrow L^q(Y)$.  Consider the open subset
$$U^{t,p} = \{x \in B^{t,p}(T_\co\fM^{s,p}) : \|x\|_{L^q(Y)} < \delta\}$$
of $B^{t,p}(T_\co\fM^{s,p})$, the $B^{t,p}$ closure of $T_\co\fM^{s,p}$.
\begin{enumerate}
  \item For $\delta$ sufficiently small, $E_\co$ extends to a bounded map $E_\co: U^{t,p} \to \fC^{t,p}(Y)$.  It is a diffeomorphism onto its image and is therefore a submanifold of $\fC^{t,p}(Y)$ contained in $\fM^{t,p}$.
  \item The constant $\delta$ can be chosen uniformly for $\co$ in a sufficiently small $L^\infty(Y)$ ball around any configuration of $\fM^{s,p}$.
\end{enumerate}
The corresponding results hold also for $\M^{s,p}$.  Finally, all the previous statements hold with the $B^{t,p}(Y)$ topology replaced with the $H^{t,p}(Y)$ topology. 
\end{Corollary}

\Proof We only do the lowest regularity case $t = 1/p$, since the case $t > 1/p$ is simpler and handled in a similar way.  For $t = 1/p$, then in trying to mimic the proof of Theorem \ref{ThmMcharts}, we show that the map $F_\co: L^3\T_\co \to L^3\T_\co$ preserves $B^{1/p,p}(Y)$ regularity on a small $L^3(Y)$ neighborhood of $0$.

In one direction, starting with $x \in \T^{1/p,p}_\co$, we want to show that $F_\co(x) \in \T^{1/p,p}_\co$. This means we must show that $Q_\co$ is bounded on $\T^{t,p}$.  We have the embedding $B^{1/p,p}(Y) \hookrightarrow L^{3p/2}(Y)$. Hence, we have a multiplication map $B^{1/p,p}(Y) \times B^{1/p,p}(Y) \hookrightarrow L^{3p/4}(Y)$.  Next, the projection $\Pi_{\K_\co^{s,p}}$ onto $\K_\co^{s,p}$ extends to a bounded map on $L^{3p/4}\T$ since $\co$ is sufficiently regular (see the proof of Lemma \ref{LemmaQ}).  Finally when we apply the inverse Hessian, we get an element of $H^{1,3p/4}(Y)$ (Proposition \ref{PropInvertH} generalizes to Sobolev spaces, see Remark \ref{RemFun}).  Since we have an embedding $H^{1,3p/4}(Y) \hookrightarrow B^{1-1/p,p}(Y) \subseteq B^{1/p,p}(Y)$, this shows that $Q_\co$ is bounded on $B^{1/p,p}(Y)$.  In the other direction, suppose $x \in L^3\T_\co$ and $F_\co(x) \in \T^{1/p,p}_\co$.  In this situation, we have $Q_\co(x) \in H^{1,3/2}\T_\co$, which embeds into $\T^{1/p,p}_\co$, and so it follows that $x \in \T^{1/p,p}_\co$.  (For $t > 1/p$, we do not have $Q_\co(x) \in \T^{t,p}_\co$, which is why we need $q > 3$ so that we have room to elliptic bootstrap.)

All the steps in Theorem \ref{ThmMcharts} follow through as before to prove the corollary for $t = 1/p$. The arithmetic for the $H^{t,p}$ spaces yields the same result.\En


\subsection{Boundary Values of the Space of Monopoles}

Define the space of tangential boundary values of monopoles
\begin{equation}
  \L^{s-1/p,p}(Y,\s) = r_\Sigma(\fM^{s-1/p,p}(Y,\s)).
\end{equation}
By (\ref{fMM}), we also have
\begin{equation}
  \L^{s-1/p,p}(Y,\s) = r_\Sigma(\M^{s,p}(Y,\s)).
\end{equation}
With $Y$ and $\s$ fixed and satisfying (\ref{assum}), we simply write $\L^{s-1/p,p} = \L^{s-1/p,p}(Y,\s)$.

We know that $\M^{s,p}$ is a manifold for $s > \max(3/p,1/2)$ by Theorem \ref{ThmMMan}.  Under further restrictions on $s$, we will see that $\L^{s-1/p,p}$ is also a manifold and the restriction map $r_\Sigma: \M^{s-1/p,p} \to \L^{s-1/p,p}$ is a covering map with fiber $\G_{h,\partial}(Y)$, which, as defined in (\ref{Ghp}), is the gauge group of harmonic gauge transformations which restrict to the identity on $\Sigma$.  Furthermore, this covering map implies that the chart maps for  $\M^{s,p}$ push forward under $r_\Sigma$ to chart maps for the manifold $\L^{s-1/p,p}$. Consequently, the nice analytic properties of the chart maps for $\M^{s,p}$ in Theorem \ref{ThmMcharts} induce chart maps for $\L^{s-1/p,p}$ that have similar desirable analytic properties.

First, we establish several important lemmas.

\begin{Lemma}\label{Lemma-Imm}
  For $s > \max(3/p,1/2)$, $r_\Sigma: \M^{s,p} \to \fC^{s-1/p,p}(\Sigma)$ is an immersion.
\end{Lemma}

\Proof This is just Corollary \ref{CorSur}(ii).\End

The following important lemma allows us to control the norm of a monopole on $Y$ in terms of the norm of its restriction on $\Sigma$.

\begin{Lemma}\label{Lemma-BI}
Let $s - 1/p > 1/2$ or $s \geq 1$ if $p = 2$. Then there exists a continuous function $\mu_{s,p}: \R^+ \to \R^+$ such that for any $\c \in \M^{s,p}$, we can find a gauge transformation $g \in \G_{h,\partial}(Y)$ such that
\begin{equation}
  \|g^*(B-B_\rf,\Psi)\|_{B^{s,p}(Y)} \leq \mu_{s,p}\left(\|r_\Sigma(B-B_\rf,\Psi)\|_{B^{s-1/p,p}(\Sigma)}\right). \label{BIest}
\end{equation}
\end{Lemma}

\Proof For the moment, assume $(B,\Psi) \in \fC(Y)$ is any smooth configuration. Define the following quantities
\begin{align}
  \E^\an(B,\Psi) &= \frac{1}{4}\int_Y |F_{B^t}|^2 + \int_Y|\nabla_B\Psi|^2 + \frac{1}{4}\int_Y(|\Psi|^2+(s/2))^2 - \int_Y \frac{s^2}{16}\\
  \E^{\mathrm{top}}(B,\Psi) &= -\int_\Sigma(D_B^\partial\Psi,\Psi) + \int (H/2)|\Psi|^2.
\end{align}
Here $s$ is the scalar curvature of $Y$, $H$ is the mean curvature of $\Sigma$, and $D^{\partial}_B$ is the boundary Dirac operator
\begin{align*}
  (D^{\partial}_B\Psi)|_\Sigma = (\rho(\nu)^{-1}D_B\Psi)|_\Sigma - (\nabla_{B,\nu}\Psi)|_\Sigma + (H/2)\Psi|_\Sigma,
\end{align*}
where $\nabla_B$ is the $\spinc$ covariant derivative determined by $B$.  Thus, $D^{\partial}_B$ only involves differentiation along the directions tangential to $\Sigma$.

If we view $(B,\Psi)$ as a time-independent configuration for the four-dimensional Seiberg-Witten equations (see the discussion before Theorem \ref{ThmUCPnonlin}), then the above quantities are the analytic and topological energy of $(B,\Psi)$, respectively, as defined in \cite{KM}. According to \cite[Proposition 4.5.2]{KM}, we have the energy identity
\begin{equation}
\E^{\an}(B,\Psi) = \E^{\tp}(B,\Psi) + \|SW_3(B,\Psi)\|^2_{L^2(Y)}. \label{energy-identity}
\end{equation}
Observe that
\begin{align*}
  \E^{\mathrm{top}}(B,\Psi) & \leq C\left(\|\Psi\|_{B^{1/2,2}(\Sigma)}^2 + \|\Psi\|_{L^3(\Sigma)}^2\|(B-B_\rf)|_{\Sigma}\|_{L^3(\Sigma)} + \|\Psi\|_{L^2(\Sigma)}^2\right)\\
  & \leq  C'\|r_\Sigma(B-B_{\rf},\Psi)\|_{B^{1/2,2}(\Sigma)}^3
\end{align*}
for some constants $C,C'$ independent of $(B,\Psi)$. Here we used the embedding $B^{1/2,2}(\Sigma) \hookrightarrow L^3(\Sigma)$.

In what follows, we will use $x \lesssim y$ to denote $x \leq Cy$ for some constant $C$ that does not depend on the configuration $(B,\Psi)$. Now consider a smooth solution of the three-dimensional Seiberg-Witten equations.  Then we have $SW_3\c=0$, and so it follows that
\begin{equation}
  \E^{\an}(B,\Psi) \lesssim \|r_\Sigma(B-B_{\rf},\Psi)\|_{B^{1/2,2}(\Sigma)}^3. \label{Eanbound}
\end{equation}
From this and the definition of $\E^\an(B,\Psi)$, we get the a priori bound
\begin{equation}
 \|\Psi\|_{L^4(Y)} \lesssim 1+\|r_\Sigma(B-B_{\rf},\Psi)\|_{B^{1/2,2}(\Sigma)}^3. \label{Psi4}
\end{equation}
If $\c \in \fM^{s,p}$ is not smooth, we can approximate $\c$ by smooth configurations by Theorem \ref{ThmMcharts}(iii).  We have $r_\Sigma(B-B_\rf,\Psi) \in \T^{s-1/p,p}_\Sigma \hookrightarrow \T^{s-1/p-\eps,2}_\Sigma$ for any $\eps > 0$ by Theorem \ref{ThmEmbed}.  Since $s > 1/2 + 1/p$, we can choose $\eps$ so that $s - 1/p - \eps > 1/2$.  Thus, we have uniform control over the $B^{1/2,2}(\Sigma)$ norm of the tangential boundary values of an approximating sequence to $\c$.  Thus, taking the limit, we see that (\ref{Psi4}) also holds for $\c \in \fM^{s,p}(Y)$.

Our remaining task is to use the a priori control (\ref{Psi4}) and the elliptic estimates for the Seiberg-Witten equations in Coulomb gauge to bootstrap our way to the estimate (\ref{BIest}).  By Corollary \ref{CorEBP}, we have the following elliptic estimate on $1$-forms $b$:
\begin{equation}
  \|b\|_{B^{t,q}(Y)} \lesssim \|db\|_{B^{t-1,q}(Y)} + \|d^*b\|_{B^{t-1,q}(Y)} + \|b|_\Sigma\|_{B^{t-1/q,q}(\Sigma)} + \|b^h\|_{B^{t-1,q}(Y)} \label{elldd*}
\end{equation}
where $b^h$ is the orthogonal projection of $b$ onto the finite dimensional space
\begin{equation}
  H^1(Y, \Sigma; i\R) \cong \{a \in \Omega^1(Y; i\R) : da = d^*a = 0, a|_\Sigma = 0\}.
\end{equation}
Here $t > 1/q$ and $q \geq 2$.

Now let $(B,\Psi) \in \M^{s,p}$ be any configuration.  Since it is in in the Coulomb slice determined by $B_\rf$, then equation (\ref{elldd*}) implies
\begin{equation}
 \|(B-B_{\rf})\|_{B^{t,q}(Y)} \lesssim \|F_{B^t}-F_{B_\rf^t}\|_{B^{t-1,q}(Y)} + \|r_\Sigma(B - B_\rf)\|_{B^{t-1/q,q}(\Sigma)} + \|(B-B_\rf)^h\|_{B^{t-1,q}(Y)} \label{Bcontrol0}
\end{equation}
where $t,q$ will be chosen later. Since Dirichlet boundary conditions are overdetermined for the smooth Dirac operator $D_{B_\rf}$, we have the elliptic estimate
\begin{equation}
\|\Psi\|_{B^{t,q}(Y)} \lesssim \|D_{B_\rf}\Psi\|_{B^{t-1,q}(Y)} + \|\Psi\|_{B^{t-1/q,q}(\Sigma)}. \label{Psicontrol0}
\end{equation}
There exists an absolute constant $C$ such that for any configuration $(B_0,\Psi_0)$, we can find a gauge transformation $g \in \G_{h,\partial}(Y)$ such that $g^*(B_0,\Psi_0)$ satisfies $\|(g^*(B_0-B_\rf))^h\|_{B^{t-1,q}(Y)} \leq C$, since the quotient of $H^1(Y,\Sigma; i\R)$ by the lattice $\G_{h,\partial}(Y)$ is a torus.  To keep notation simple, redefine $(B,\Psi)$ by such a gauge transformation.  Such a gauge transformation preserves containment in $\M^{s,p}$ since the monopole equations are gauge invariant and the Coulomb-slice is preserved by $\G_{h,\partial}(Y)$.  So using the bound $\|(B-B_\rf)^h\|_{B^{t-1,q}(Y)} \leq C$ and the identity $SW_3\c = 0$, the bounds (\ref{Bcontrol0}) and (\ref{Psicontrol0}) become
\begin{align}
\|(B-B_{\rf})\|_{B^{t,q}(Y)} & \lesssim \|\Psi^2\|_{B^{t-1,q}(Y)} + \|(B - B_\rf)|_\Sigma\|_{B^{t-1/q,q}(\Sigma)} + 1.\label{Bcontrol} \\
\|\Psi\|_{B^{t,q}(Y)} &\lesssim \|\rho(B-B_\rf)\Psi\|_{B^{t-1,q}(Y)} + \|\Psi\|_{B^{t-1/q,q}(\Sigma)}. \label{Psicontrol}
\end{align}
We will use these estimates, bootstrapping in $t$ and $q$ and using the a priori control (\ref{Psi4}), to get the estimate (\ref{BIest}).

Let us first consider the case $p=2$ and $s \geq 1$. Letting $t = 1$ and $q = 2$, (\ref{Bcontrol}) and (\ref{Psi4}) yield \begin{align*}
  \|B-B_{\rf}\|_{B^{1,2}(Y)} & \lesssim 1 + \|\Psi\|_{L^4(Y)} + \|(B - B_\rf)|_\Sigma\|_{B^{1/2,2}(\Sigma)}\\
   & \lesssim 1+\|r_\Sigma(B-B_{\rf},\Psi)\|_{B^{1/2,2}(\Sigma)}^3.
\end{align*}
This yields control over $\|B-B_\rf\|_{L^4(Y)}$ since we have the embedding $B^{1,2}(Y) \hookrightarrow L^6(Y)$.  Using this estimate in (\ref{Psicontrol}) with $t=1$, $q=2$, to control $\rho(B-B_{\rf})$, we have
\begin{align*}
  \|\Psi\|_{B^{1,2}(Y)} & \lesssim \|B-B_\rf\|_{L^4(Y)}\|\Psi\|_{L^4(Y)} + \|\Psi\|_{B^{1/2,2}(\Sigma)}\\
  & \lesssim 1+\|r_\Sigma(B-B_{\rf},\Psi)\|_{B^{1/2,2}(\Sigma)}^6.
\end{align*}
This proves the estimate for $s = 1$.  The estimate (\ref{BIest}) for $s\geq 1$ now follows from boostrapping the elliptic estimates (\ref{Bcontrol}) and (\ref{Psicontrol}) in $t$.  Indeed, once we gain control over $\|(B,\Psi)\|_{B^{t,q}(Y)}$, we can control the quadratic terms $ \|\Psi^2\|_{B^{t'-1,q}(Y)}$ and $\|\rho(B-B_\rf)\Psi\|_{B^{t'-1,q}(Y)}$ for some $t' > t$ as long as $t' \leq s$.  After finitely many steps of bootstrapping, we get (\ref{BIest}), where the function $\mu_{s,p}$ can be computed explicitly if desired.

For $p > 2$ and $s \geq 1$, we use the imbedding $\fC^{s-1/p,p}(\Sigma) \hookrightarrow \fC^{s-1/p-\eps,2}(\Sigma)$, for any $\eps > 0$.  From the previous case, we find that we can control $\|(B,\Psi)\|_{B^{s,2}(Y)}$ in terms of $\|(B,\Psi)\|_{B^{s-1/2,2}(\Sigma)}$.  Since $B^{s,2}(Y) \subseteq B^{1,2}(Y) \hookrightarrow L^6(Y)$, the quadratic terms in (\ref{Bcontrol}) and (\ref{Psicontrol}) lie in $L^3(Y)$.  Since we have the embedding $L^3(Y) \subset B^{0,q}(Y)$, where $q = \max(3,p)$, we can repeat the bootstrapping process (in $t$) as in the previous case to the desired estimate (\ref{BIest}) for any $s \geq 1$ and $p \leq 3$.  Suppose $p > 3$.  Then with $q = 3$ in the previous step, we have established (\ref{Bcontrol}) and (\ref{Psicontrol}) with $t = 1$ and $q = 3$. Since $B^{1,3}(Y) \hookrightarrow L^{q}(Y)$ for any $q < \infty$,  we have control of the quadratic terms of (\ref{Bcontrol}) and (\ref{Psicontrol}) in $L^p$ for any $p < \infty$.  Thus, we have the estimate (\ref{Bcontrol}) and (\ref{Psicontrol}) for $t = 1$ and $q = p$, since $L^p(Y) \subseteq B^{0,p}(Y)$.  We can then bootstrap in $t$ to the estimate (\ref{BIest}) for any $s \geq 1$ and $p < \infty$.  Thus, we have taken care of the case $s \geq 1$ and all $p \geq 2$.

Finally, suppose $s < 1$ and $p > 2$.  We employ the same strategy of bootstrapping in $q$ until we get to $p$.  Since $s - 1/p > 1/2$, we have $B^{s-1/p,p}(\Sigma) \hookrightarrow B^{1/2,2}(\Sigma)$ and so we have control of $\|(B-B_{\rf},\Psi)\|_{B^{1,2}(Y)}$ and $\|(B-B_{\rf},\Psi)\|_{L^6(Y)}$ in terms of $\|(B-B_\rf)|_\Sigma\|_{B^{1/2,2}(\Sigma)}$.  Let $1/2 < t = s < 1$ and $q = \min(3,p)$ in (\ref{Bcontrol}) and (\ref{Psicontrol}).  We have control of the quadratic terms on the right-hand side since $L^3(Y) \subset B^{0,q}(Y) \subset B^{s-1,q}(Y)$, since $s - 1 \leq 0$.  Thus, we have the control (\ref{BIest}) for $p = q$.  If $p \leq 3$, we are done.  Else $p > 3$ and we bootstrap in $q$.  Indeed, starting with $q_1 = 3$, we have a map $B^{s,q_i}(Y) \times B^{s,q_i}(Y) \to B^{2s-3/q_i,q_i}(Y) \hookrightarrow L^{q_{i+1}}(Y) \subseteq B^{s-1,q_{i+1}}(Y)$, where $q_{i+1} = q_i/(2(1-sq_i)) > q_i$.  Using (\ref{Bcontrol}) and (\ref{Psicontrol}), we thus
bootstrap to the estimate (\ref{BIest}) with $p = q_{i+1}$ from the estimate (\ref{BIest}) with $p = q_i$.  The $q_i$ keep increasing until after finitely many steps, we get to the desired $p$, thereby proving (\ref{BIest}).\End

The next lemma tells us that any two monopoles which have the same restriction to $\Sigma$ are gauge equivalent on $Y$.

\begin{Lemma}\label{LemmaCover}
   Let $s > \max(3/p,1/2)$. If $\cx[1],\cx[2] \in \fM^{s,p}$ and $r_\Sigma\cx[1]= r_\Sigma\cx[2]$, then $\cx[1]$ and $\cx[2]$ are gauge equivalent on $\fC^{s,p}(Y)$.
\end{Lemma}

\Proof Because of (\ref{fMM}), without loss of generality, we can suppose $\cx[1],\cx[2] \in \M^{s,p}$.  There are two cases to consider.  In the first case, one and hence both the configurations are reducible.  Indeed, if say $\cx[1]$ is reducible, then $\Psi_2|_\Sigma = \Psi_1|_\Sigma = 0$. Since $D_{B_2}\Psi_2 = 0$, by unique continuation for $D_{B_2}$, we have $\Psi_2 \equiv 0$ so that $\cx[2]$ is also reducible.  In this reducible case, then $B_1$ and $B_2$ are both flat connections and so by (\ref{assum}), we must have $H^1(Y,\Sigma) = 0$.  Since $d(B_1 - B_2) = d^*(B_1 - B_2) = 0$ and by hypothesis $(B_1 - B_2)|_\Sigma = 0$, we must then have $B_1 - B_2 = 0$ since $H^1(Y,\Sigma) = 0$.  So in this case, $B_1$ and $B_2$ are in fact equal. In the second case, neither configuration is reducible.  In this case, consider the $4$-manifold $S^1 \times Y$ and regard $\cx[1]$ and $\cx[2]$ as time-independent solutions to the Seiberg-Witten equations on $S^1 \times Y$.  We now apply Theorem \ref{ThmUCPnonlin}.\End

Piecing the previous lemmas together, we can now prove the rest of our main theorem concerning the monopole spaces:

\begin{Theorem}\label{ThmMtoL}
Let $s > \max(3/p,1/2+1/p)$.  Then $\L^{s-1/p,p}$ is a closed Lagrangian submanifold of $\fC^{s-1/p,p}(\Sigma)$.  Furthermore, the maps
\begin{align}
  r_\Sigma: \fM^{s,p} &\to \L^{s-1/p,p} \label{rfML}\\
  r_\Sigma: \M^{s,p} &\to \L^{s-1/p,p} \label{rML}
\end{align}
are a submersion and a covering map respectively, where the fiber of the latter is the lattice $\G_{h,\partial}(Y) \cong H^1(Y,\Sigma)$.
\end{Theorem}

\Proof By (\ref{fMM}), it suffices to consider the map (\ref{rML}). By Lemma \ref{Lemma-Imm}, the map $r_\Sigma: \M^{s,p} \to \fC^{s-1/p,p}(\Sigma)$ is an immersion, hence a local embedding.  The previous lemma implies that (\ref{rML}) is injective modulo $G := \G_{h,\partial}(Y)$, since the gauge transformations which restrict to the identity and preserve Coulomb gauge are precisely those gauge transformations in $G$. Moreover, $G$ acts freely on $\M^{s,p}(Y)$ by assumption (\ref{assum}), since when there are reducible solutions, we have $G = 1$.

It remains to show that $r_\Sigma: \M^{s,p}/G \to \fC^{s-1/p,p}(\Sigma)$ is an embedding onto its image. Let  $\cx[i] \in \M^{s,p}$, $i \geq 1$, be such that $r_\Sigma\cx[i] \to r\cx[0]$ in $\fC^{s-1/p,p}(\Sigma)$ as $i \to \infty$.  We want to show that given any subsequence of the $\cx[i]$, there exists a further subsequence convergent to an element of the $G$ orbit of $\cx[0]$. This, combined with the fact that (\ref{rML}) is a local embedding will imply that (\ref{rML}) is a global embedding, modulo the covering transformations $G$. Indeed, the local embedding property tells us that there exists a open neighborhood $V_\co \ni \co$ of $\M^{s,p}$ such that $r_\Sigma: V_\co \to \fC^{s-1/p,p}(\Sigma)$ is an embedding onto its image, and moreover, $r_\Sigma(g^*V_\co) = r_\Sigma(V_\co)$ for all $g \in G$.  Proving the above convergence result shows that given a sufficiently small neighborhood $U$ of $r_\Sigma\cx[0]$ in $\fC^{s-1/p,p}(\Sigma)$, then $U \cap \L^{s-1/p,p}$ is contained in the image of any one of the embeddings $r_\Sigma: g^*V_{\co} \to \fC^{s-1/p,p}(\Sigma)$, $g \in G$.  Otherwise, we could find a subsequence $(B_{i'},\Phi_{i'})$ of the $\cx[i]$ such that $r_\Sigma(B_{i'},\Phi_{i'}) \to r_\Sigma\co$ but the $(B_{i'},\Phi_{i'})$ lie outside all the $g^*V_\co$, a contradiction.

Without further ado then, by Lemma \ref{Lemma-BI}, we know we can find gauge transformations $g_i \in G$ such that $g_i^*\cx[i]$ is uniformly bounded in $B^{s,p}(Y)$, since $r_\Sigma\cx[i] \to r_\Sigma\cx[0]$ is uniformly bounded.  For notational simplicity, redefine the $\cx[i]$ by these gauge transformations.  Thus, since the $\cx[i]$ are bounded in $B^{s,p}(Y)$, any subsequence contains a weakly convergent subsequence.  Let $\cx[\infty] \in \M^{s,p}(Y)$ be a weak limit of some subsequence $\cx[i']$.  We have $r_\Sigma\cx[\infty] = r_\Sigma\cx[0]$, and so $\cx[\infty]$ and $\cx[0]$ are gauge equivalent by an element of $G$.  If we can show that that $\cx[i'] \to \cx[\infty]$ strongly in $B^{s,p}(Y)$, then we will be done. Due to the compact embedding $B^{s,p}(Y) \hookrightarrow B^{t,p}(Y)$, for $t < s$, we have $\cx[i] \to \cx[\infty]$ strongly in the topology $B^{s - \eps,p}(Y)$, $\eps > 0$. If we can boostrap this to strong convergence in $B^{s,p}(Y)$, we will be done.  To show this, we use the ellipticity of $\tH_{\cx[\infty]}$.  Let $(b_i,\psi_i) = (B_i-B_\infty,\Psi_i-\Psi_\infty)$.  We have the elliptic estimate
\begin{equation}
  \|(b_i,\psi_i)\|_{B^{s,p}(Y)} \lesssim \|\tH_{\cx[\infty]}(b_i,\psi_i)\|_{B^{s-1,p}(Y)} + \|r_\Sigma(b_i,\psi_i)\|_{B^{s-1/p,p}(\Sigma)}. \label{ellconv}
\end{equation}
This follows because $\tH_{\cx[\infty]}$ is elliptic and the boundary term controls the kernel of $\tH_{\cx[\infty]}|_{\C^{s,p}}$ by Corollary \ref{CorSur}(ii).  The last term of (\ref{ellconv}) tends to zero and for the first term, we have
\begin{equation}
  \tH_{\cx[0]}(b_i,\psi_i) = (b_i,\psi_i)\#(b_i,\psi_i) \label{g^2}
\end{equation}
from (\ref{SWlin}), since $\cx[i], \cx[\infty] \in \M^{s,p}$.  We have a continuous multiplication map $B^{s-\eps,p}(Y) \times B^{s-\eps,p}(Y) \to B^{s-\eps,p}(Y) \subset \B^{s-1,p}(Y)$ for $s - \eps > 3/p$.  Since $\cx[i] \to \co$ strongly in $B^{s-\eps,p}(Y)$, we have that (\ref{g^2}) goes to zero in $B^{s-1,p}(Y)$, which means $(b_{i'},\psi_{i'})$ goes to zero in $B^{s,p}(Y)$ by (\ref{ellconv}).  Thus, $\cx[i] \to \cx[\infty]$ strongly in $B^{s,p}(Y)$.

It now follows that $r_\Sigma: \fM^{s,p}(Y) \to \fC^{s-1/p,p}(\Sigma)$ is a covering map onto a embedded submanifold, where the fiber of the cover is $G$.  Moreover, the proof we just gave also shows that $\L^{s-1/p,p}$ is a closed submanifold, since if $r_\Sigma\cx[i]$ is a convergent sequence, it is convergent to $r_\Sigma\cx[\infty]$ for some $\cx[\infty] \in \M^{s-1/p,p}$.  Finally, Theorem \ref{ThmLinLag}(i) shows that $\L^{s-1/p,p}$ is Lagrangian submanifold of $\fC^{s-1/p,p}(\Sigma)$, since its tangent space at any point is a Lagrangian subspace of $\T^{s-1/p,p}_\Sigma$.\End

Since $r_\Sigma: \M^{s,p} \to \L^{s-1/p,p}$ is a covering, the chart maps on $\M^{s,p}$ push forward and induce chart maps on $\L^{s-1/p,p}$.  Indeed, at the tangent space level, we already know we have isomorphisms
$r_\Sigma: T_{\co}\M^{s,p} \to T_{r_\Sigma\co}\L^{s-1/p,p}$ and $P_\co: T_{r_\Sigma\co}\L^{s-1/p,p} \to T_{\co}\M^{s,p}$ inverse to one another, where recall $P_\co$ is the Poisson operator given by Theorem \ref{ThmLinLag}.  Because $\M^{s,p}$ is locally a graph over $ T_{\co}\M^{s,p}$, then $\L^{s-1/p,p}$ is locally a graph over $T_{r_\Sigma\co}\L^{s-1/p,p}$.  To analyze this properly, we also want to ``push foward" the local straightening map $F_\co$ for $\M^{s,p}$ at a configuration $\co$, defined in Lemma \ref{LemmaF}, to obtain a local straightening map $F_{\Sigma,\co}$ for $\L^{s-1/p,p}$ at $r_\Sigma\co$.

\begin{Lemma}\label{LemmaFSigma}
  Let $s > \max(3/p,1/2+1/p)$ and let $\co \in \M^{s,p}$.  Define the spaces
  $$X_\Sigma = \T_\Sigma^{s-1/p,p}, \qquad X_{\Sigma,0} = T_{r_\Sigma\co}\L^{s-1/p,p}, \qquad X_{\Sigma,1} = J_\Sigma X_{\Sigma,0}.$$
  We have $X_\Sigma = X_{\Sigma,0}\oplus X_{\Sigma,1}$ and we can define the smooth map \label{p:FSigma}
  \begin{align}
  F_{\Sigma,\co}: V_\Sigma & \to X_{\Sigma,0}\oplus X_{\Sigma,1}\nonumber \\
  x = (x_0,x_1) & \mapsto (x_0, x_1 - r_\Sigma E^1_{\co}(P_\co x_0)), \label{eq4:FSigma}
\end{align}
where $V_\Sigma \subset X_\Sigma$ is an open subset containing $0$ and $E^1_\co$ is defined as in Theorem \ref{ThmMcharts}.  For any $\max(1/2,2/p) < s' \leq s-1/p$, we can take $V_\Sigma$ to contain a $B^{s',p}(\Sigma)$ ball, i.e., there exists a $\delta > 0$, depending on $r_\Sigma\co$, $s'$, and $p$, such that
$$V_\Sigma \supseteq \{x \in X_{\Sigma} : \|x\|_{B^{s',p}(\Sigma)} < \delta\}.$$
Moreover, we have the following:
\begin{enumerate}
  \item We have $F_{\Sigma,\co}(0)=0$ and $D_0F_{\Sigma,\co} = \id$.  For $V_\Sigma$ sufficiently small, $F_{\Sigma,\co}$ is a local straightening map for $\L^{s-1/p,p}$ at $r_\Sigma\co$ within the neighorhood $V_\Sigma$.
  \item We can choose $\delta$ uniformly for $r_\Sigma\co$ in a sufficiently small $B^{s',p}(\Sigma)$ neighborhood of any configuration in $\L^{s-1/p,p}$.
\end{enumerate}
\end{Lemma}

\Proof We have $F_{\Sigma,\co}(0)=0$ since $E^1_\co(P_\co(0)) = E^1_\co(0) = 0$, and $D_0F_{\Sigma,\co} = \mathrm{id}$ since $D_0E^1_\co = 0$ by Theorem \ref{ThmMcharts}.  Moreover, we see that $F_{\Sigma,\co}$ can be defined on a $B^{s',p}(\Sigma)$ ball containing $0 \in X_\Sigma$. Indeed, $P_\co$ maps such a ball into a $B^{s'+1/p,p}(Y)$ ball inside $T_{\co}\M^{s,p}$, we have the embedding $B^{s'+1/p,p}(Y) \hookrightarrow L^\infty(Y)$ by our choice of $s'$, and the domain of $E^1_\co$ contains an $L^\infty(Y)$ ball by Theorem \ref{ThmMcharts}.  Take $V_\Sigma$ to be such a $B^{s',p}(\Sigma)$ ball.

It now follows from $\L^{s-1/p,p} \subset \L^{s',p}$ and the fact that $r_\Sigma: \M^{s'+1/p,p} \to \L^{s',p}$ is a covering map onto a globally embedded submanifold (by Theorem \ref{ThmMtoL}) that $F_{\Sigma,\co}$ is a local straightening map for $\L^{s-1/p,p}$ within a $B^{s',p}(\Sigma)$ neighborhood of $0 \in X_\Sigma$.  (Shrinking $V_\Sigma$ if necessary, let this neighborhood be $V_\Sigma$.)  In more detail, if $x \in V_\Sigma$ and $F_{\Sigma,\co}(x) \in X_{\Sigma,0}$, then
\begin{equation}
  r_\Sigma\co + x = r_\Sigma\co + (x_0, r_\Sigma E^1_{\co}(P_\co x_0)), \label{eq4:LemmaFSigma}
\end{equation}
which means that
\begin{equation}
  r_\Sigma\co + x = r_\Sigma\Big(\co + P_\co x_0 + E^1_\co(P_\co x_0)\Big), \label{eq4:MtoLcharts}
\end{equation}
where $\co + P_\co x_0 + E^1_\co(P_\co x_0) \in \M^{s,p}$ by Theorem \ref{ThmMcharts}.  Thus, $r_\Sigma\co + x \in r_\Sigma(\M^{s,p}) = \L^{s-1/p,p}$.  Conversely, if $x \in V_\Sigma$ is such that $r_{\Sigma\co} + x \in \L^{s-1/p,p}$, then (having chosen $V_\Sigma$ small enough) we must have
 \begin{equation}
   r_\Sigma\co + x  = r_\Sigma\Big(\co + x_0' + E^1_\co(x_0')\Big).
 \end{equation}
for some $x_0' \in T_{\co}\M^{s'+1/p,p}$ since $\L^{s-1/p,p} \subset \L^{s',p}$ and $r_\Sigma: \M^{s'+1/p,p} \to \L^{s',p}$ is a local diffeomorphism from a neighborhood of $\co \in \M^{s'+1/p,p}$ onto a neighborhood of $r_\Sigma\co \in \L^{s',p}$.  Since $P_\co: T_{r_\Sigma\co}\L^{s',p} \to T_{\co}\M^{s'+1/p,p}$ is an isomorphism, then $x_0' = P_\co x_0$ for some $x_0 \in T_{r_\Sigma\co}\L^{s',p}$ and so (\ref{eq4:MtoLcharts}), hence (\ref{eq4:LemmaFSigma}) must hold.  By definition of $F_{\Sigma,\co}$, which extends to a well-defined map on the $B^{s'+1/p,p}(\Sigma)$ topology, (\ref{eq4:LemmaFSigma}) implies $x_0 = F_{\Sigma,\co}(x)$.  But $F_{\Sigma,\co}$ acting on a neighborhood of $0$ in $\T^{s',p}_\Sigma$ preserves the $B^{s-1/p,p}(\Sigma)$ topology, so $x_0 \in T_{r_\Sigma\co}\L^{s-1/p,p}$ since $x \in \T^{s-1/p,p}_\Sigma$.  Thus, $F_{\Sigma,\co}(x) \in X_{\Sigma,0}$.   Moreover, both $F_{\Sigma,\co}$ and $F^{-1}_{\Sigma,\co}$ are invertible when restricted to $V_\Sigma$, since the inverse $F^{-1}_{\Sigma,\co}$ is simply given by
\begin{equation}
  F_{\Sigma,\co}^{-1}(x) = (x_0, x_1 + r_\Sigma E^1_{\co}(P_\co x_0)). \label{FSigma-1}
\end{equation}
Altogether, this shows that $F_{\Sigma,\co}$ is a local straightening map for $\L^{s-1/p,p}$ within $V_\Sigma$.

(ii) This follows from the uniformity statement of Theorem \ref{ThmMcharts}(ii) and the continuous dependence of $P_\co: \T^{s',p}_\Sigma \to \T^{s'+1/p,p} \hookrightarrow L^\infty\T$ with respect to $\co$ (see Theorem \ref{ThmLinLag}(iv)).  Here, we use the fact that if $r_\Sigma\co \in \L^{s-1/p,p}$ varies continuously in a small $B^{s',p}(\Sigma)$ neighborhood, then one can choose $\co \in \M^{s,p}$ continuously in a small $B^{s'+1/p,p}(Y)$ neighborhood, since $r_\Sigma: \M^{s'+1/p,p} \to \L^{s',p}$ and $r_\Sigma: \M^{s,p} \to \L^{s-1/p,p}$ are covers.\End

With the above lemma, we have local straightening maps for our Banach manifold $\L^{s-1/p,p}$.  Then from Theorem \ref{ThmMcharts} and the general framework of Appendix \ref{AppIFT}, we have the following theorem for the local chart maps for $\L^{s-1/p,p}$.

\begin{Theorem}\label{ThmLcharts} Let $s > \max(3/p,1/2+1/p)$.
\begin{enumerate}
  \item Let $\co \in \M^{s,p}$. Then there exists a neighborhood $U \subset T_{r_\Sigma\co}\L^{s-1/p,p}$ of $0$ and a map $E^1_{r_\Sigma\co}: U \to X_{\Sigma,1}$, where $X_{\Sigma,1}$ is as in Lemma \ref{LemmaFSigma}, such that the map \label{p:Erco}
\begin{align}
  E_{r_\Sigma\co}: U & \to \L^{s-1/p,p} \nonumber \\
  x &\mapsto r_\Sigma\co + x + E^1_{r_\Sigma\co}(x) \label{eq4:ESigma}
\end{align}
is a diffeomorphism of $U$ onto a neighborhood of $r_\Sigma\co$ in $\L^{s-1/p,p}$.  Furthermore, the map $E^1_{r_\Sigma\co}$ smooths by one derivative, i.e. $E^1_{r_\Sigma\co}(x) \in \T_\Sigma^{s+1-1/p,p}$ for all $x \in U$. \label{p:E1rco}
\item For any $\max(1/2,2/p) < s' \leq s-1/p$, we can choose $U$ such that both $U$ and $E_{r_\Sigma\co}(U)$ contain $B^{s',p}(U)$ neighborhoods, i.e., there exists a $\delta > 0$, depending on $r_\Sigma\co$, $s'$, and $p$, such that
\begin{align*}
  U & \supseteq \{x \in T_{r_\Sigma\co}\L^{s-1/p,p} : \|x\|_{B^{s',p}(\Sigma)} < \delta\}\\
  E_{r_\Sigma\co}(U) & \supseteq \{\c \in \L^{s-1/p,p} : \|\c - r_\Sigma\co\|_{B^{s',p}(\Sigma)} < \delta\},
\end{align*}
The constant $\delta$ can be chosen uniformly in $r_\Sigma\co$, for all $r_\Sigma\co$ in a sufficiently small $B^{s-1/p,p}(Y)$ neighborhood of any configuration in $\L^{s-1/p,p}$.
 \item Smooth configurations are dense in $\L^{s-1/p,p}$.
\end{enumerate}
\end{Theorem}

\begin{figure}
\centering
\includegraphics{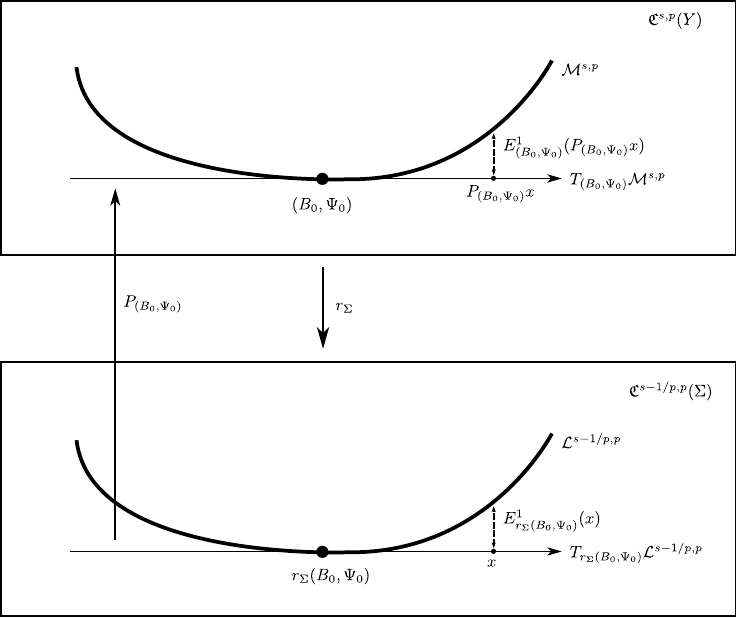}
\caption{A chart map for $\M^{s,p}$ at $(B_0,\Phi_0)$ induces a chart map for $\L^{s-1/p,p}$ at 
$r_\Sigma(B_0,\Phi_0)$.}
\label{Figure 2a}
\end{figure}

\Proof (i) As in (\ref{eq4:EfromF}), the chart map $E_{r_\Sigma\co}$ is determined by restricting $F_{\Sigma,\co}^{-1}$, the inverse of the local straightening map $F_{\Sigma,\co}$, to a neighborhood of $0$ in the tangent space $T_{r_\Sigma\co}\L^{s-1/p,p}$.  Thus, we have
\begin{equation}
    E_{r_\Sigma\co}(x) = r_\Sigma\co + F_{\Sigma,\co}^{-1}(x), \qquad x \in U:= F_{\Sigma,\co}(V_\Sigma) \cap T_{r_\Sigma\co}\L^{s-1/p,p}, \label{eq4:ESigmaFSigma}
\end{equation}
where $V_\Sigma$ is defined as in Lemma \ref{LemmaFSigma}. The expression for $F_{\Sigma,\co}^{-1}$ is given by (\ref{FSigma-1}). Thus, (\ref{eq4:ESigmaFSigma}) and the definition of $E^1_\co$ in (\ref{eq4:ESigma}) yields
\begin{equation}
E^1_{r_\Sigma\co}(x) = r_\Sigma E^1_{\co}(P_\co x) \label{E1L}
\end{equation}
The mapping properties of $E^1_{r_\Sigma\co}$ now follow from Theorem \ref{ThmMcharts}.

(ii) This is a direct consequence of Lemma \ref{LemmaFSigma}(ii).

(iii) This just follows from $r_\Sigma: \M^{s,p}\to\L^{s-1/p,p}$ being a cover and the density of smooth configurations in $\M^{s,p}$ by Theorem \ref{ThmMcharts}.\En

\begin{Corollary}\label{CorLpLag}Suppose $\co \in \M^{s,p}$.
\begin{enumerate}
  \item If $U$ is a sufficiently small $L^p(\Sigma)$ neighborhood of $0$ in $L^pT_{r_\Sigma\co}\L^{s-1/p,p}$, then $E_{r_\Sigma\co}$ extends to a bounded map
 \begin{align}
   E_{r_\Sigma\co}: U \to L^p\fC(\Sigma). \label{eq4:E1Lp}
 \end{align}
 The map (\ref{eq4:E1Lp}) is a diffeomorphism onto its image and hence $E_{r_\Sigma\co}(U)$ is an $L^p$ submanifold of $L^p\fC(\Sigma)$ contained in $L^p\L$.
 \item The $L^p(\Sigma)$ topology above can be replaced with $B^{t,p}(\Sigma)$ for any $0 \leq t \leq s-1/p$ and $H^{t,p}(\Sigma)$ for any $0 \leq t \leq s - 1/p$.
\end{enumerate}

\end{Corollary}

\Proof We use Corollary \ref{CorMCharts} to show that $E^1_{r_\Sigma\co}$ is bounded on the $L^p(\Sigma)$ topology.  We only prove the lowest regularity case $s = 0$, since the other cases are similar (and more easily handled).  We have the inclusion $L^p \subset B^{0,p}$ since $p \geq 2$.  By Theorem \ref{ThmLinLag}, the Poisson operator $P_\co$ maps $B^{0,p}(\Sigma)$ to $B^{1/p,p}(Y)$ for $s \geq 0$ since $\co$ is sufficiently regular.  In the proof of Corollary \ref{CorMCharts}, we showed that $E^1_\co$ maps $B^{1/p,p}(Y)$ to $H^{1,3p/4}(Y)$.  Hence when we apply $r_\Sigma$, we find altogether from (\ref{E1L}) that $E^1_{r_\Sigma\co}(x)$ belongs to $B^{1-4/3p,3p/4}(\Sigma) \hookrightarrow L^p(\Sigma)$.  Thus, $E^{1}_\co$ is bounded on $L^p(\Sigma)$ and $B^{0,p}(\Sigma)$.  In the above calculation, we implicitly used the fact that $H^{1/p,p}(Y) \hookrightarrow L^3(Y)$, so that $P_\co$ maps a small $L^p(\Sigma)$ ball into the domain of $E^1_\co$, which contains an $L^3(Y)$ ball by Corollary \ref{CorMCharts}.\En

\begin{Rem}\label{RemFun}
  We have mentioned before that since our analysis works on a variety of function spaces, it is merely a matter of convenience that we worked primarily with Besov spaces on $Y$.  In the above corollary and elsewhere, we see how the usual $L^p$ spaces can be employed as well.  Corollary \ref{CorLpLag} will be significant in \cite{N2}, since we will need to consider, locally, $L^p$ closures of $\L$.  We conclude by noting that every instance in which the $B^{s,p}(Y)$ topology is used in this paper, the topology $H^{s,p}(Y)$ may be used instead.  These spaces, known as the Bessel potential spaces, are defined in Appendix \ref{AppFun}.  For $s$ a nonnegative integer and $1<p<\infty$, we have $H^{s,p}(Y) = W^{s,p}(Y)$, the Sobolev space of functions having $s$ derivatives belonging to $L^p(Y)$.  When $p = 2$, $H^{s,2}(Y) = B^{s,2}(Y)$ for all $s$.  Furthermore, the spaces $H^{s,p}$ and $B^{s,p}$ are ``close" to each other in the sense that $H^{s_1,p}(Y) \subseteq B^{s_2,p}(Y) \subseteq H^{s_3,p}(Y)$ for all $s_1 > s_2 > s_3$.
  Moreover, one sees that all the foundational analysis in Appendix \ref{AppEBP} applies equally to Bessel potential and Besov spaces.

  We should note that two particular places where it is important that Sobolev spaces may be used in addition to Besov spaces are Lemma \ref{LemmaTdecomp} and Proposition \ref{PropInvertH}.  Indeed, their proofs rely only on function space arithmetic and elliptic estimates arising from elliptic boundary value problems.  For both of these, Sobolev spaces can be used all the same, and so we can replace every occurrence of the $B^{\bullet,\bullet}(Y)$ topology with the $H^{\bullet,\bullet}(Y)$ topology in Lemma \ref{LemmaTdecomp} and Proposition \ref{PropInvertH}.  One can now check that the statements of all our lemmas and theorems concerning Besov spaces on $Y$ also hold for their Sobolev counterparts.

  For the purposes of \cite{N2}, it is also important that we can replace Besov spaces on $\Sigma$ with Sobolev spaces on $\Sigma$ as well, but with some care, since the space of boundary values of a Sobolev space is still a Besov space.  We already saw how to do this in Corollary \ref{CorLpLag}.  We should note that for the Calderon projection $P_\co^+$ in Theorem \ref{ThmLinLag}, where $\co \in B^{s,p}(Y)$, one also has that
  \begin{equation}
    P^+_\co: H^{t-1/p,p}\T_\Sigma \to H^{t-1/p,p}\T_\Sigma, \qquad t < s + 1.
  \end{equation}
  is bounded.  This follows from the fact that $\pi^+$ is bounded on $H^{t-1/p,p}\T_\Sigma$, as it is a pseudodifferential operator, and
$$(P^+_\c - \pi^+): H^{t-1/p,p}\T_\Sigma \subseteq \T_\Sigma^{t-1/p,p} \to \T_\Sigma^{\min(s-1/p+1,\; t-1/p+1),p} \subset H^{t-1/p,p}\T_\Sigma$$
  by Theorem \ref{ThmLinLag}(iv) and Theorem \ref{ThmEmbed}.
\end{Rem}

From the above remark, the Sobolev version of our main theorem, with the $H^{s,p}(Y)$ topology replaced with the $B^{s,p}(Y)$ topology, holds.  In fact, one can see from this that the Besov monopole space $\M^{s,p}$ is actually equal to the Sobolev monopole space $H^{s,p}\M$.\footnote{From the density of smooth configurations, it suffices to show that the tangent space to $H^{s,p}\M$ and $\M^{s,p}$ at a smooth monopole $\c$ are both equal.  However, this follows from the fact that the kernel of an elliptic operator (in our case, the operator $\tH_\c$) in the $H^{s,p}$ and $B^{s,p}$ topologies are equivalent.  This follows from the results of Appendix \ref{AppSeeley}, which shows that these spaces are isomorphic (modulo a finite dimensional subspace) to their space of boundary values, which is a fixed subspace of $B^{s-1/p,p}$ on the boundary.} Finally, let us remark that the proof of our main corollary easily follows from the work we have done.\\

\noindent\textit{Proof of Main Corollary: } For every coclosed $1$-form $\eta$, the zero set of $SW_3(B,\Phi)=(\eta,0)$ is  gauge-invariant.  Thus, all the methods of Section \ref{SecLinSW} apply to the linearization of the monopole spaces associated to the perturbed equations.  Next, we still have the transversality result Lemma \ref{LemmaHtrans} so long as we modify the assumption (\ref{assum}) to $c_1(\s) \neq 2[*\eta]$ or $H^1(Y,\Sigma)=0$.  The energy estimates in Lemma \ref{Lemma-BI} still hold in the perturbed case since we still have (\ref{energy-identity}) and the uniform bound $\|SW_3(B,\Psi)\|_{L^2(Y)} = \|\eta\|_{L^2(Y)}$.  Finally, the unique continuation results from Appendix \ref{AppUCP} still apply, since we always apply these results to the difference of solutions to the perturbed Seiberg-Witten equations, and the equation satisfied by the difference is independent of $\eta$.  Thus, all our methods and hence results carry through in the perturbed case.


\appendix

\section{Some Functional Analysis}\label{AppFA}

\subsection{Subspaces and Projections}

Here, we collect some properties about projections and subspaces of Banach spaces. Given a Banach space $X$, a projection $\pi: X\to X$ is a bounded operator such that $\pi^2 = \id$.  A closed subspace $U \subset X$ is complemented if there exists another closed subspace $V \subset X$ such that $X = U \oplus  V$ as topological vector spaces.  We call $X = U \oplus V$ a topological decomposition of $X$.  A closed subspace $U$ is complemented if and only if there exists a projection $\pi$ whose range is $U$.  In this case, we have the decomposition
$$X = \im \pi \oplus \ker \pi.$$
Recall that any finite dimensional subspace of a Banach space is complemented.  Likewise, any subspace of finite codimension is also complemented.  Thus, if $Y \subset X$ has finite (co)dimension, we may always regard $X/Y$ as a subspace of $X$ (though unless $X$ is a Hilbert space, there is in general no canonical embedding $X/Y \hookrightarrow X$).

The following simple lemma tells us that if a projection $\pi$ restricted to a subspace $U' \subset X$ yields a Fredholm map $\pi: U' \to \im \pi$, then $U'$ is essentially a graph over $\im \pi$. More precisely, we have the following:

\begin{Lemma}\label{LemmaFredproj}
 Let $X = U \oplus V$ and let $\pi$ be the projection onto the first factor.  Let $U'$ be a subspace of $X$ and suppose $\pi: U' \to U$ is Fredholm.  Then
 \begin{equation}
   U' = \{x + Tx: x \in \pi(U')\} \oplus F, \label{app:graphT}
 \end{equation}
 where $F = \ker (\pi |_{U'})$ is finite dimensional and $T: \pi(U') \to V$ is a bounded operator. Consequently, $U'$ is also a complemented subspace of $X$, in particular, it is the range of a projection.
\end{Lemma}

\Proof Since $F$ is finite dimensional, it is the range of a bounded projection $\pi_F: X \to F$.  Since $F \subset U'$, then $(1-\pi_F): U' \to U'$ maps $U'$ into itself and its range is a complement of $F$ in $U'$.  It follows that $\pi: (1-\pi_F)(U') \to \pi(U')$ is an isomorphism.  Let $\bar\pi$ denote this isomorphism.  Thus,
\begin{align*}
  U' &= \{(\bar\pi)^{-1}x : x \in \pi(U')\} \oplus F\\
     &= \{x + (-1 + (\bar\pi)^{-1})x : x \in \pi(U')\} \oplus F.
\end{align*}
Let $T = (-1 + (\bar\pi)^{-1}): \pi(U') \to X$.  Since $\im T \subset \ker \pi$, we see that $\im T \subset V$.  This gives us the desired decomposition of $U'$.  One can now explicitly write a projection onto $U'$.  Since $\pi(U') \subset U$ has finite codimension, it has a complement $C \subset U$ along with a projection $\pi_C: U \to C$ such that the complementary projection $(1 -\pi_C)$ has range equal to $\im \bar\pi$.  A projection from $X$ onto $U'$ is now easily seen to be given by the map
\begin{equation}
  (1+T)(1 - \pi_C)\pi(1-\pi_F) \oplus \pi_F.\label{eq:Fredproj}
\end{equation}
This proves the lemma.\End

Given a complemented subspace $U \subset X$, to compare other subspaces $U' \subset X$ with $U$, then we should not only study the projection of $U'$ onto $U$ but also onto a complement of $U$.

\begin{Def}\label{DefComm}
  Let $X$ be a Banach space.  Two projections $\pi$ and $\pi'$ on $X$ are \textit{commensurate} if $\pi - \pi'$ is compact.  Given a complemented subspace $U \subset X$, then a subspace $U'$ is \textit{commensurate} with $U$ if its projection onto $U$ is Fredholm and its projection onto some complement of $U$ is compact. It follows from Lemma \ref{LemmaFredproj} that this definition is independent of the choice of complement of $U$.  We will also say that $U'$ is a \textit{compact perturbation} of $U$.
\end{Def}

\begin{Corollary}
  Let $U$ and $U'$ be as in Lemma \ref{LemmaFredproj}.  Then the subspace $U'$ is commensurate with $U$ if and only if the map $T$ in (\ref{app:graphT}) is compact.  In this case, the space $U$ is also commensurate with $U'$.
\end{Corollary}

Hence, being commensurate is a symmetric relation, and we may simply speak of two subspaces $U$ and $U'$ as being commensurate.  The notion of commensurability obviously captures the notion of two subspaces being ``close" to one another in a functional analytic sense. On the opposite spectrum, one may consider pairs of subspaces that form a direct sum decomposition modulo finite dimensional subspaces.  More precisely, we have the following definition:

\begin{Def}\label{DefFred}
  A pair of complemented subspaces $(U,V)$ of a Banach space $X$ is \textit{Fredholm} if $U \cap V$ is finite dimensional and the algebraic sum $U + V$ is closed and has finite codimension.  In this case, we say that $(U,V)$ form a Fredholm pair, or more simply, that $U$ and $V$ are Fredholm (in $X$).
\end{Def}

Together, the notion of a pair of subspaces being either commensurate or Fredholm will be very important in what we do. Next, we record the following technical lemmas concerning topological decompositions:

\begin{Lemma}\label{LemmaDenseFred}
  Let $X$ and $Y$ be Banach spaces, with $Y \subset X$ dense.  Suppose $X = X_1 \oplus X_0$ and $Y \cap X_i \subseteq X_i$ is dense for $i = 0,1$.  Then if $Y \cap X_1$ and $Y \cap X_0$ are Fredholm in $Y$, then in fact $Y = (Y \cap X_1) \oplus (Y \cap X_0)$.
\end{Lemma}

\Proof The hypotheses imply $Y = (Y \cap X_1) \oplus (Y \cap X_0) \oplus F$ where $F$ is some finite dimensional subspace of $Y$.  If we take the closure of this decomposition in $X$, we have $X \supseteq X_1 \oplus X_0 \oplus F$, which means $F = 0$.\En

\begin{Lemma}\label{LemmaComm}
  Let $X = X_1 \oplus X_0$ be a topological decomposition of $X$ and let $\pi_i: X \to X_i$ be the coordinate projections.  Let $U, V \subset X$ be subspaces and let $V = V_1 \oplus V_0$, where $V_1 = V \cap X_1$ and $\pi_0: V_0 \to X_0$ is Fredholm.  If $U$ is commensurate with $V$, then we can write $U = U_1 \oplus U_0$, where $U_1 = U\cap X_1$, $\pi_0: U_0 \to X_0$ is Fredholm, and $U_i$ is commensurate with $V_i$, $i=0,1$.
\end{Lemma}

\Proof
By the preceding analysis, since $U$ is commensurate with $V$, there exist finite dimensional subspaces $F_1 \subset X$ and $F_2 \subset V$ and a compact operator $T: V/F_2 \to X$ such that $U = \{x + Tx: x\in V/F_2\} \oplus F_1$.  For notational simplicity, let us suppose $F_1=F_2 = 0$, since the conclusion is unaffected by finite dimensional errors.  So then
\begin{align*}
    U &= \{x + Tx: x\in V\}\\
    &= \{x + Tx: x \in V_1\} + \{x + Tx: x \in V_0\}\\
    &=: U_1' + U_0'.
\end{align*}
Since $T$ is compact, then $U_0'$ is commensurate with $V_0$ and since $\pi_0: V_0 \to X_0$ is Fredholm, so is $\pi_0: U_0' \to X_0$.  Thus the map
$$\pi_0' = \pi_0: U_0'/\ker\pi_0 \to \pi_0(U_0')$$
is an isomorphism.  Let $V_1' \subset V_1$ be the subspace of finite codimension defined by
$$V_1' := \{x \in V_1 : \pi_0(Tx) \in \pi_0(U_0')\},$$
In other words, $V_1'$ is the subspace of $V_1$ such that the space $\{x + Tx: x \in V_1'\} \subseteq U_1'$ differs from an element of $X_1$ by an element of $U_0'$.  Indeed, we have
$$\{x + Tx - {\pi_0'}^{-1}\pi_0T(x): x \in V_1'\} \subseteq X_1$$
since it is annihilated by $\pi_0$. We thus have
$$U_1 = U \cap X_1 = \{x + Tx - {\pi_0'}^{-1}\pi_0T(x): x \in V_1'\} + \ker(\pi_0|_{U_0'}).$$
From this expression for $U_1$, it follows that $U_1$ is commensurate with $V_1$. Letting $U_0 = U_0' + \{x + Tx: x \in V_1/V_1'\}$, then $U = U_1 \oplus U_0$ and all the properties are satisfied.\End

\begin{Lemma}\label{LemmaFredDecomp}Let $U_0$, $U_1$, $V_0$, and $V_1$ be subspaces of a Banach space $X$ such that we have the topological decompositions
  \begin{align}
    X &= U_0 \oplus U_1 = V_0 \oplus V_1 \label{UVdecomp1}\\
      &= V_0 \oplus U_1 = U_0 \oplus V_1. \label{UVdecomp2}
  \end{align}
Let $\pi_{U_0,U_1}$ denote the projection onto $U_0$ through $U_1$ and similarly for other pairs of complementary spaces in the above.  Since $\pi_{U_0,U_1}: V_0 \to U_0$ and $\pi_{U_1,U_0}: V_1 \to U_1$ are isomorphisms, then $V_i$ is the graph of a map $T_{V_i}: U_i \to U_{i+1}$, $i = 0, 1\; \mathrm{mod}\; 2$.
\begin{enumerate}
  \item We have the following formulas:
    \begin{align}
    \pi_{V_0,U_1} &= (1+T_{V_0})\pi_{U_0,U_1} \label{proj1}\\
    \pi_{U_1,V_0} &= 1 - \pi_{V_0,U_1} \label{proj2}\\
    \pi_{V_0,V_1} &= (1+T_{V_0})\pi_{U_0,U_1}(1+T_{V_1}\pi_{U_1,V_0})^{-1} \label{proj3}
      \end{align}
      and likewise with the $0$ and $1$ indices switched.
  \item If $V_i$ is commensurate with $U_i$, for $i = 0,1$, then $\pi_{V_0,V_1}$ is commensurate with $\pi_{U_1,U_0}$.  This remains true even if we drop the assumption (\ref{UVdecomp2}).
\end{enumerate}
\end{Lemma}

\Proof (i) Formula (\ref{proj1}) is just the definition of $T_{V_0}$; indeed, this is the special case of the formula (\ref{app:graphT}) when the map $\pi: U' \to U$ is an isomorphism.  Formula (\ref{proj2}) is tautological since $U_1$ and $V_0$ are complementary.  It remains to establish (\ref{proj3}).  Let $\Lambda: X \to X$ be the isomorphism of $X$ which maps $V_0$ identically to $V_0$ and $U_1$ to $V_1$ using the graph map $T_{V_1}$.  In other words, $\Lambda$ is given by
\begin{align*}
  \Lambda &= \pi_{V_0,U_1} + (1+T_{V_1})\pi_{U_1,V_0}\\
  &= 1 + T_{V_1}\pi_{U_1,V_0}.
\end{align*}
The map $\pi_{V_0,V_1}$ is now easily seen to be given by $\pi_{V_0,U_1}\Lambda^{-1}$, which yields (\ref{proj3}).  By symmetry, these formulas hold with $0$ and $1$ indices reversed.

(ii) In this case, the maps $T_{V_i}$ are compact, $i = 0,1$.  It follows from (\ref{proj3}) that $\pi_{V_0,V_1} - \pi_{U_0,U_1}$ is compact.  If (\ref{UVdecomp2}) does not hold, we proceed as follows.  Let $F$ denote the finite dimensional space spanned by the kernel and cokernel of the Fredholm maps $\pi_{U_0,U_1}: V_0 \to U_0$ and $\pi_{U_1,U_0}: V_1 \to U_1$.  Let $\bar X = X/F$ be regarded as a subspace of $X$ and let $\bar \pi: X \to \bar X$ be the projection through $F$.  It follows that we can choose finite codimensional subspaces $U_i' \subseteq U_i$ and $V_i' \subseteq V_i$ such that, letting $\bar U_i = \pi_i(U_i')$ and $\bar V_i = \pi_i(V_i')$, we have
$$\bar X = \bar U_0 \oplus \bar U_1 = \bar V_0 \oplus \bar V_1.$$
By construction of $\bar X$, we also have
$$\bar X = \bar V_0 \oplus \bar U_1 = \bar U_0 \oplus \bar V_1,$$
since now $\bar V_i$ is a graph over $\bar U_i$.  On $\bar X$, we can therefore conclude that the projections $\pi_{\bar V_1,\bar V_0}$ and $\pi_{\bar U_1,\bar U_0}$ are commensurate.  These operators also act on $X$ since we can define them to be zero on $F$, in which case, $\pi_{V_0,V_1}$ and $\pi_{U_0,U_1}$ are finite rank perturbations of $\pi_{\bar V_0, \bar V_1}$ and $\pi_{\bar U_0,\bar U_1}$, respectively.  It now follows that $\pi_{V_1, V_0}$ and $\pi_{U_1,U_0}$ are also commensurate.\End

\begin{Rem}\label{RemError}
  In all applications, our Banach space $X$ under consideration will be a function space of configurations on a manifold, and the compact operators that arise will be maps that smooth by a certain number of derivatives $\sigma \geq 0$ (e.g. the operator maps a Besov space $B^{s,p}$ to a more regular Besov space $B^{s+\sigma,p}$).  In this way, if additionally we have that all finite dimensional subspaces which arise in the above analysis are spanned by elements that are smoother than those of $X$ by $\sigma$ derivatives, one can ensure that all compact perturbations occurring in the projections constructed in the above lemmas continue to be operators that are smoothing of order $\sigma$.  In other words, the amount of smoothing is preserved in all our constructions.
\end{Rem}

The notion of commensurability of two spaces is one qualitative way of measuring two spaces as being close.  Alternatively, we may regard two subspaces $V_1$ and $V_2$ of $X$ as being close if $V_2$ is the graph over $V_1$ of a map with small norm, i.e. $V_2 = \{x + Tx: x \in V_1\}$ where $V_1^\bot$ is any fixed complement of $V_1$ and $T: V_1 \to V_1^\bot$ is an operator with small norm.  If the norm of $T$ is small enough, we can replace $V_1^\bot$ with $X$.  This motivates the following definition:

\begin{Def}\label{DefSS}
\begin{enumerate}
\item A \textit{continuous family of subspaces} $\{V(\sigma)\}_{\sigma \in \frak{X}}$ of $X$, where $\frak{X}$ is a topological space, is a collection of complemented subspaces $V(\sigma)$ of $X$ such that the following local triviality condition holds: for any $\sigma_0 \in \frak{X}$, there exists an open neighborhood $U \ni \sigma_0$ in $\frak{X}$ such that for all $\sigma \in U$, there exists a map $T_{\sigma_0}(\sigma): V \to X$ such that the induced map
  \begin{align}
    V(\sigma_0) &\to V(\sigma) \nonumber \\
    x & \mapsto x + T(\sigma)x \label{graphT}
  \end{align}
is an isomorphism.  The map $T_{\sigma_0}(\sigma)$ varies continuously in the operator norm topology with respect to $\sigma \in U$.
\item  A \textit{smooth family of subspaces} $V(t)$ of $X$, $t \in \R$, is a continuous family of subspaces for which $\frak{X} = \R$ and the maps $T(t)$ in (\ref{graphT}) vary smoothly in operator norm topology.
    \end{enumerate}
\end{Def}

This definition is such that one can construct operators associated to a continuously varying family subspaces in a continuous way, e.g., projections onto such subspaces.  Likewise for the smooth situation. To illustrate this, we state the following trivial lemma for small time intervals:

\begin{Lemma}\label{ThmStraight}
  Let $V(t)$ be a continuous (smooth) family of subspaces of $X$, for $t \in \R$.  Then for any $t_0 \in \R$, we can find an open interval $I$ containing $t_0$, and a continuous (smooth) family of isomorphisms $\Phi(t): X \to X$, $t \in I$, such that $\Phi(t)(V(t_0)) = V(t)$ for all $t \in \R$, with $\Phi(0)=\id$.
\end{Lemma}

\Proof Without loss of generality, let $t_0 = 0$ and suppose we are in the smooth case, with the continuous case being the same.  Let $V(0)^\bot$ be any complement of $V(0)$ in $X$. Then for small enough $t$, $V(t)$ is also a complement of $V(0)^\bot$, and we can define
\begin{align*}
  \Phi(t): V(0)\oplus V(0)^\bot & \to V(t) \oplus V(0)^\bot\\
  (x, y) & \mapsto (x + T(t)x, y),
\end{align*}
where $x \mapsto x + T(t)x$ is the isomorphism from $V(0)$ to $V(t)$ given by the definition of $V$ being a smooth family of subspaces of $X$.  The maps $\Phi(t)$ are smooth since the $V(t)$ are.\End

In other words, the family of spaces $V(t)$ has local trivializations given by the $\Phi(t)$ which identify the $V(t)$ with $V(t_0)$, for $t \in I$.  Given a family of spaces complementary to the $V(t)$ and which vary smoothly, one can construct the $\Phi(t)$ for all $t$, but the above local result will suffice for our purposes.

\subsection{Symplectic Banach Spaces}\label{SecSymp}

Let $X$ be a real Banach space endowed with a skew-symmetric bilinear form $\omega$.  Then $X$ is a \textit{(weakly) symplectic Banach space} if $\omega$ is nondegenerate, i.e., the map $\omega: X \to X^*$ which assigns to $x \in X$ the linear functional $\omega(x,\cdot)$ is injective. If $X$ is a Hilbert space and there exists an automorphism $J: X \to X$ such that $J^2=-\id$ and $\omega(\cdot,J\cdot)$ is the inner product on $X$, we say that $X$ is a \textit{strongly symplectic Hilbert space} and that $J$ is the compatible complex structure. (As a word of caution, many other authors define a symplectic Banach space to be one for which $\omega: X \to X^*$ is an isomorphism, but that will never be the case for us unless $X$ is a strongly symplectic Hilbert space.)

Given any subspace $V$ of a symplectic Banach space $X$, let $\Ann(V) \subset X$ denote its annihilator with respect to the symplectic form.  A (co)isotropic subspace $V$ is one for which $V \subseteq (\supseteq)\; \Ann(V)$.  A \textit{Lagrangian subspace} $L$ is an isotropic subspace which has an isotropic complement.  This implies $L$ is also coisotropic by the nondegeneracy of $\omega$.  In case $X$ is a strongly symplectic Hilbert space, then in fact, an isotropic subspace is Lagrangian if and only if it is coisotropic, see \cite{Wein}.  In this latter case, any Lagrangian subspace $L$ has an orthogonal Lagrangian complement $JL$.

The following procedure, known as symplectic reduction, is well-known in the context of Hilbert spaces (see e.g. \cite[Proposition 6.12]{KL}):

\begin{Theorem}\label{ThmSR}(Symplectic Reduction) Let $(X,\omega)$ be a strongly symplectic Hilbert space with compatible complex structure $J$.  Let $U \subseteq X$ be a closed coisotropic subspace.  Let $L \subset X$ be a Lagrangian subspace such that $L + \Ann(U)$ is closed.  Then $U \cap JU$ is a strongly symplectic Hilbert space and the orthogonal projection $\pi_{U \cap JU}$ onto $U \cap JU$, yields a map
\begin{equation}
  \pi_{U \cap JU}: L \cap U \to U \cap JU \label{eq:app-SR}
\end{equation}
whose image $\pi_{U \cap JU}(L \cap U)$ is a Lagrangian subspace of $U \cap JU$.
\end{Theorem}

We call the map (\ref{eq:app-SR}) the symplectic reduction of $L$ with respect to $U$.  For symplectic reduction on Banach spaces, we can generalize the above result as follows:

\begin{Corollary}\label{CorSR}
    Given the notation from the previous theorem, let $Y$ be a Banach space with $Y \subseteq X$ dense. Given any subspace $V \subset X$, define $V_Y := Y \cap V$. Suppose $\pi_{U \cap JU}$ and $J$ map $Y$ into itself and that $L_Y$ and $U_Y$ are dense in $L$ and $U$, respectively. Suppose
    $\pi_{U \cap JU}(L_Y)$ and $J\pi_{U \cap JU}(L_Y)$ are Fredholm in $U_Y \cap JU_Y$.  Then $\pi_{U \cap JU}(L_Y)$ and $J\pi_{U \cap JU}(L_Y)$ are complementary Lagrangian subspaces of the symplectic Banach space $U_Y \cap JU_Y$.
\end{Corollary}

\Proof This follows from the previous theorem and Lemma \ref{LemmaDenseFred}.\End

\section{Banach Manifolds and The Inverse Function Theorem}\label{AppIFT}

Taking the usual definition of a finite dimensional manifold, one may replace all occurrences of Euclidean space with some other fixed Banach space, thereby obtaining the notion of a (smooth) Banach manifold.  In other words, a Banach manifold, modeled on a Banach space $X$, is a Hausdorff topological space that is locally homeomorphic to $X$ and whose transition maps are all diffeomorphisms\footnote{A map of Banach spaces is smooth if it is infinitely Fr\'echet differentiable.  A diffeomorphism is a smooth map that has a smooth inverse.}.

In a similar way, one also obtains the notion of a (smooth) Banach submanifold of a Banach space.  More precisely, we have the following definition:


\begin{Definition}\label{DefSub}
  Let $X$ be a Banach space.  A \textit{Banach submanifold} $M$ of $X$ is a subspace of $X$ (as a topological space) that satisfies the following.  There exists a closed Banach subspace $Z \subset X$ such that at every point $u \in M$, there exists an open set $V$ in $X$ containing $u$ and a diffeomorphism $\Phi$ from $V$ onto an open neighborhood of $0$ in $X$ such that $\Phi(V \cap M) = \Phi(V) \cap Z$.  We say that $M$ is modeled on the Banach space $Z$.
\end{Definition}

We almost always drop explicit reference to the model Banach space $Z$, since it will be clear what this space is in practice.  Of course, one can consider abstract Banach manifolds that do not come with a global embedding into a Banach space, but such a situation will not occur for us.  The above definition coincides with the usual definition of a submanifold when $X$ is a Euclidean space.

In the general situation above, we have no information about the local chart maps $\Phi$.  However, if $M$ is defined in some natural way, say as the zero set of some function, one can construct a more concrete local model for $M$.  The tools we use for this are the inverse and implicit function theorems in the general setting of Banach spaces. Below, we record these theorems, mostly to fix notation in applications (proofs can be found in e.g. \cite{La}). Let $X = X_0 \oplus X_1$ be a direct sum of Banach spaces and $f: X \to Y$ a smooth map of Banach spaces.  For any $x \in X$, let $D_xf: X \to Y$ denote the Fr\'echet derivative of $f$ at $x$.

\begin{Theorem}\label{ThmIFT} Suppose $D_0f: X \to Y$ is surjective, with $D_0 f: X_1 \to Y$ an isomorphism and $X_0 = \ker D_0 f$.
\begin{enumerate}
  \item (Implicit Function Theorem) Choose $V$ to be an open neighborhood of $0$ in $X$ such that $D_xf: X_1 \to Y$ remains an isomorphism for all $x \in V$. Then $M := f^{-1}(0) \cap V$ is a Banach submanifold of $X$ modeled on $X_0$.
  \item (Inverse Function Theorem) Define the smooth map $F: X_0 \oplus X_1 \to X$ by
  $$F(x_0,x_1) = (x_0, (D_0f|_{X_1})^{-1}f(x_0,x_1)).$$
  Then $F(0) = 0$, $D_0F = \id$, and shrinking $V$ if necessary, we can arrange that both $F$ and $F^{-1}$ are diffeomorphisms onto their images when restricted to $V$.  In this case, we have $M \subseteq F^{-1}(F(V) \cap X_0)$.
\end{enumerate}
\end{Theorem}

\begin{Def}\label{DefLDF}Let $M \subset X$ be a Banach submanifold and let $u \in M$ be any element, which without loss of generality, we let be $0$.  Given a function $f: X \to Y$ as in (i) above, we say that $f$ is a \textit{local defining function} for $M$ near $u$ if there exists a neighborhood $V$ of $u \in X$ such that $M_u := M \cap V$ is a Banach submanifold of $X$ and satisfies $M_u = f^{-1}(0) \cap V$. In this case, the function $F$ associated to $f$ in Theorem \ref{ThmIFT}(ii) is said to be a \textit{local straightening map} for $M$ at $u$.  If we wish to emphasize our choice of $V$, we will say that $F$ is a local straightening map \textit{within the neighborhood $V$}.
\end{Def}

The names we give for $f$ and $F$ are natural given their role in describing $M$.  Namely, the manifold $M_u$, which is an open neighborhood of $u$ in $M$, is the subset of $V$ that lies in the preimage under $f$ of the regular value $0 \in Y$.  On the other hand, the map $F$ is a local diffeomorphism of $X$ which straightens out $M_u$ to an open neighborhood $U := F(M_u)$ inside the tangent space $X_0 = T_uM$.  Consequently, $F^{-1}: U \to M$ is a diffeomorphism of $U$ onto its image $M_u$, an open neighborhood of $0 \in M$.

\begin{Def}\label{DefChart}
  With the above notation, we call $F^{-1}: U \to M$ the \textit{induced chart map} of the local straightening map $F$.
\end{Def}

Thus, while a Banach submanifold has no distinguished choice of local charts near any given point, a local straightening map gives us a canonical choice for one.  We will be consistently using this choice when constructing local chart maps for the Banach submanifolds we study.

\section{Function Spaces}\label{AppFun}

In this section, we define the various function spaces needed for our analysis.  We establish enough of their properties so that we may apply them in the context of elliptic boundary value problems and nonlinear partial differential equations.  In some sense, for the purposes of this paper, the precise definitions of the function spaces below are not as important as their formal properties under such operations as restricting to a hyperplane and multiplication (see Theorems \ref{ThmTrace} and Theorem \ref{ThmMult}).


\subsection{The Classical Function Spaces}

We define the classical Sobolev, Bessel potential, and Besov spaces.  These spaces along with their basic properties are well documented, e.g., see \cite{ET}, \cite{Tr}, and \cite{Tr1}.  The proofs of all the statements here can be found in those references.

\subsubsection*{\ref{AppFun}.1.1 \hspace{.3cm}Function Spaces on $\R^n$}

We begin by defining our spaces on $\R^n$ with coordinates $x_j$, $1 \leq j \leq n$.  Let $\S(\R^n)$ be the space of rapidly decaying Schwartz functions and let $\S'(\R^n)$ be its dual space, the space of tempered distributions.  Given $f \in \S(\R^n)$, we have the Fourier transform
$$\F f(\xi) = \int e^{i\xi\cdot x}f(x)dx.$$
The Fourier transform extends to $\S'(\R^n)$ by duality.  Given a multi-index $\alpha = (\alpha_1,\ldots, \alpha_n) \in \Z_+^n$ of nonnegative integers, we let $D^\alpha f = \partial_{x_1}^{\alpha_1}\cdots\partial_{x_n}^{\alpha_n}f$ be the corresponding partial derivatives of $f$ in the sense of distributions.

Next, we consider a dyadic partition of unity as follows.  Let $\psi(\xi)$ be a smooth bump function, $0 \leq \psi(\xi) \leq 1$, with $\psi(\xi)$ equal to $1$ on $|\xi| \leq 1$ and $\psi$ identically zero on $|\xi| \geq 2$.  Let
\begin{align*}
  \varphi_0(\xi) &= \psi(\xi)\\
  \varphi_j(\xi) &= \psi(2^{-j}\xi) - \psi(2^{-j+1}\xi), \quad j\geq 1.
\end{align*}
Then we have $\mathrm{supp}\, \varphi_j \subseteq [2^{j-1},2^{j+1}]$ for $j \geq 1$ and $\sum_{j=0}^\infty \varphi_j(\xi) \equiv 1$.

Given a tempered distribution $f$, we let
$$f_j = \F^{-1}\varphi_j\F f$$
be its $j$th dyadic component.  The decomposition of $f$ into its dyadic components $\{f_j\}_{j=0}^\infty$ is known as the Littlewood-Paley decomposition.

On $\R^n$, let $L^p(\R^n)$ and $C^\alpha(\R^n)$ denote the usual Lebesgue and Holder spaces of order $p$ and $\alpha$, respectively, where $1\leq p \leq \infty$ and $\alpha \geq 0$.  In addition to these, we have the following classical function spaces:

\begin{Definition}(i) For $s \in \Z_+$ a nonnegative integer and $1 \leq p \leq \infty$, define the \textit{Sobolev spaces}
 \begin{align}
   W^{s,p}(\R^n) & = \{f \in \S'(\R^n) : \quad \|f\|_{W^{s,p}} = (\sum_{|\alpha| \leq s}\|D^\alpha f\|_{L^p}^p)^{1/p} < \infty\}, \qquad p < \infty\\
   W^{s,\infty}(\R^n) & = \{f \in \S'(\R^n) : \|f\|_{W^{s,\infty}} = \sup_{|\alpha| \leq s}\|D^\alpha f\|_{L^\infty} < \infty\}.
  \end{align}
(ii) For $s \in \R$ and $1 < p < \infty$, define the \textit{Bessel potential spaces}
   \begin{align}
   H^{s,p}(\R^n) := \Big\{f \in \S'(\R^n) : \|f\|_{H^{s,p}} = \Big\|\Big(\sum_{j=0}^\infty |2^{sj}f_j|^2\Big)^{1/2}\Big\|_{L^p}<\infty\Big\}. \label{defHsp}
 \end{align}
(iii) For $s \in \R$, $1 < p < \infty$, define the \textit{Besov spaces}\footnote{The classical Besov spaces are usually denoted with two parameters $B^s_{p,q}$.  We take $p = q$.  There are also many other equivalent norms that can be used to define the Besov spaces.  Our choice of norm reflects their similarity with $H^{s,p}$.} \label{p:funspace}
\begin{align}
  B^{s,p}(\R^n) &= \Big\{f \in \S'(\R^n): \|f\|_{B^{s,p}} = \Big\|\Big(\sum_{j=0}^\infty |2^{sj}f_j|^p\Big)^{1/p}\Big\|_{L^p}<\infty\Big\}. \label{defBsp}
\end{align}
(iv) Define $A^{s,p}$ to be shorthand for either $H^{s,p}$ or $B^{s,p}$.  The spaces $A^{s,p}$ are also a special case of what are known as \textit{Triebel-Lizorkin spaces}.
\end{Definition}

Of the above Banach spaces, the Sobolev spaces $W^{s,p}$ are the ones most naturally occurring for many of the basic problems in analysis.  The Bessel potential spaces $H^{s,p}$ arise from (complex) interpolation between the Sobolev spaces, where we may think of $f \in H^{s,p}$ as having $s$ derivatives in $L^p$.  This is most clearly illustrated when $p = 2$, where then $H^{s,p}$ is usually just denoted as $H^s$.   For general $p$, we have the following result:

\begin{Theorem}\cite[Theorem 2.3.3]{Tr} \label{ThmHW}
  For $1 < p < \infty$, $H^{s,p}(\R^n) = W^{s,p}(\R^n)$ for $s$ a nonnegative integer.
\end{Theorem}
Indeed, when $s = 0$, then Theorem \ref{ThmHW} tells us that
$$\left\|\Big(\sum_{j=0}^\infty |f_j|^2\Big)^{1/2}\right\|_{L^p} \sim \|f\|_{L^p}.$$
This is the classical Littlewood-Paley Theorem.

The Besov spaces naturally arise because they are the boundary values of Sobolev spaces.  More precisely, let $\R^{n-1} \subset \R^n$ be the hyperplane $x_n = 0.$  Given a fixed $m \in \Z_+$ and $f$ a function on $\R^n$, let
\begin{equation}
  r_m f = \left(f|_{\R^{n-1}}, \partial_{x_n}f|_{\R^{n-1}}, \ldots, \partial_{x_n}^{m}f|_{\R^{n-1}}\right) \label{rest-m}
\end{equation}
be the trace of $f$ of order $m$ along the hyperplane $\R^{n-1}$.  We have the following theorem:

\begin{Theorem}\label{ThmTrace}
\begin{enumerate}
  \item For $s > m + 1/p$ and $m \in \Z_+$, the trace map $r_m$ extends to a bounded operator
\begin{equation}
  r_m: H^{s,p}(\R^n) \to \oplus_{j=0}^{m-1} B^{s-1/p-j}(\R^{n-1})
\end{equation}
\item For any $s \in \R$ and $m \in \Z_+$, there exists an extension map $e_m: \oplus_{j=0}^{m-1} B^{s-1/p-j,p}(\R^{n-1}) \to H^{s,p}(\R^n)$ such that for $s > m + 1/p$, we have $r_m e_m = \id$.
\item $H^{s,p}(\R^n)$ may be replaced with $B^{s,p}(\R^n)$ in the above.
\end{enumerate}
\end{Theorem}

When $p = 2$, we have
\begin{equation}
  B^{s,2}(\R^n) = H^s(\R^n)
\end{equation}
for all $s$, and so the above theorem is a generalization of the fact that the trace of an element of $H^s(\R^n)$ lies in $H^{s-1/2}(\R^n)$ for $s > 1/2$.  Furthermore, because $\ell^p \subseteq \ell^q$ whenever $p \geq q$, we have the trivial inclusions
\begin{alignat*}{2}
  B^{s,p}(\R^n) & \subseteq H^{s,p}(\R^n) & \qquad & p \leq 2\\
  H^{s,p}(\R^n) & \subseteq B^{s,p}(\R^n) & \qquad &  p \geq 2.
\end{alignat*}

For $s > 0$, we can also write the Besov space norm in terms of finite differences in space rather than in terms of the Littlewood-Paley decomposition in frequency space.  For any $h \in \R^n$, define the operator
$$\delta_h f = f(x + h) - f(x).$$
Using this operator, we have the following proposition:


\begin{Proposition}\label{PropBesovDiff}
  For $s > 0$ and $1 < p < \infty$, let $m$ be any integer such that $m > s$.  Then an equivalent norm for $B^{s,p}(\R^n)$ is given by
  \begin{equation}
    \|f\|_{B^{s,p}(\R^n)} = \|f\|_{L^p(\R^n)} + \left(\int_{\R^n}^\infty \Big\||h|^{-s}\delta_h^mf\Big\|_{L^p(\R^n)}^p \frac{1}{|h|^n}dh\right)^{1/p}. \label{app.eq:besov}
  \end{equation}
\end{Proposition}

\begin{Rem}\label{RemHB}
  The spaces $H^{s,p}(\R^n)$ and $B^{s,p}(\R^n)$ satisfy
  $$H^{s_1,p}(\R^n) \subseteq B^{s_2,p}(\R^n) \subseteq H^{s_3,p}(\R^n)$$
  for all $s_1 > s_2 > s_3$, for $1<p<\infty$.  This is a simple consequence of the definitions (\ref{defHsp}) and (\ref{defBsp}).  Thus, we see that the most important features of the $B^{s,p}$ and $H^{s,p}$ spaces are determined by the exponents $s,p$, with the distinction between the Besov and Bessel potential topologies for fixed $s$ and $p$ being a more refined property.  In this sense, for most purposes, the spaces $B^{s,p}$ and $H^{s,p}$ are ``nearly identical", and many results concerning one of these spaces implies the same result for the other.  This is why we adopt the common notation of using $A^{s,p}$ to denote either $H^{s,p}$ or $B^{s,p}$.  Whenever, $A^{s,p}$ appears in multiple instances in a statement or formula, we always mean that all instances of $A^{s,p}$ are either $H^{s,p}$ or $B^{s,p}$.
\end{Rem}


We have the following fundamental properties:

\begin{Proposition}
  Let $s \in \R$ and $1 < p < \infty$.  Then the space of compactly supported functions $C^\infty_0(\R^n)$ is dense in $A^{s,p}(\R^n)$. Moreover, $A^{-s,p'}(\R^n)$ is the dual space of $A^{s,p}(\R^n)$, where $1/p+1/p' = 1$.
\end{Proposition}

\begin{Proposition}\label{LiftProp}(Lift Property)
  Let $s \in \R$ and $1<p<\infty$.  Then
  $$A^{s,p}(\R^n) = \{f \in A^{s-1,p}(\R^n) : \frac{\partial f}{\partial x^i} \in A^{s-1,p}(\R^n), 1 \leq i\leq n\}.$$
\end{Proposition}

\subsubsection*{\ref{AppFun}.1.1 \hspace{.3cm}Function Spaces on an Open Subset of $\R^n$}

Let $\Omega$ be an open subset of $\R^n$.  Unless otherwise stated, we assume for simplicity that $\Omega$ is bounded and has smooth boundary, though many of the results that follow carry over for more general open sets. Given any tempered distribution $f \in \S'(\R^n)$, we can consider its restriction $r_\Omega(f)$ to $(C^\infty_0(\Omega))'$.  Then we have the corresponding function spaces on $\Omega$:

\begin{Definition}
  For $s \in \Z_+$ and $1 \leq p \leq \infty$, the space $W^{s,p}(\Omega)$ is the space of restrictions to $\Omega$ of elements of $W^{s,p}(\R^n)$, where the norm on $W^{s,p}(\Omega)$ is given by
  $$\|f\|_{W^{s,p}(\Omega)} = \inf_{g : r_\Omega(g) = f}\|g\|_{W^{s,p}(\R^n)}.$$
  For $s \in \R$ and $1<p<\infty$, the spaces $H^{s,p}(\Omega)$ and $B^{s,p}(\Omega)$ are defined similarly.
\end{Definition}
If we consider the function space
$$\tilde A^{s,p}(\Omega) := \{f \in A^{s,p}(\R^n) : \supp f \subset \overline{\Omega}\},$$
then an equivalent definition of $A^{s,p}(\Omega)$ is
\begin{equation}
  A^{s,p}(\Omega) = A^{s,p}(\R^n)/\tilde A^{s,p}(\R^n\setminus\bar\Omega). \label{R^n_+}
\end{equation}
Furthermore, we have the following:
\begin{Proposition}
  Let $-\infty < s < \infty$ and $1 < p < \infty$. Then $C^\infty_0(\Omega)$ is dense in $\tilde A^{s,p}(\Omega)$.  Moreover, $A^{-s,p'}(\Omega)$ is the dual space of $\tilde A^{s,p}(\Omega)$, where $1/p+1/p'=1$.
\end{Proposition}

Define the upper half-space
$$\R^n_+ = \{(x_1,\ldots,x_n) \in \R^n : x_n > 0\}.$$
We have the following extension property:

\begin{Theorem}\label{ThmExt}
  Let $1 < p < \infty$.  For any $k \in \N$, there exists an extension operator
  $$E_k: A^{s,p}(\R^n_+) \to A^{s,p}(\R^n)$$
  for $|s| < k$.
\end{Theorem}

\subsubsection*{\ref{AppFun}.1.1 \hspace{.3cm}Function Spaces on Manifolds}

Ultimately, the function spaces which are important for us are those which are defined on manifolds (with and without boundary). Let $X$ be a compact $n$-manifold or an open subset of it.  We can assign to $X$ the data of an atlas $\{(U_i,\varphi_i,\Phi_i)\}$, where: (1) the $U_i$ are a finite open cover of $X$; (2) the $\varphi_i$ are a partition of unity with $\supp \varphi_i \subset U_i$; (3) each $\Phi_i$ is a map from $U_i$  to $\R^n$, where $\Phi_i$ is a diffeomorphism onto an open subset of $\R^n$ if $\bar U_i \subset \overset{\circ}{X}$ or otherwise, $\Phi_i$ is a diffeomorphism onto an open subset of $\overline{\R^n_+}$ with $\Phi_i(U_i \cap \partial X) \subset \partial\overline{\R^n_+}$.  With this data, we can define the function spaces $A^{s,p}(X)$ in terms of the function spaces on $\R^n$ and $\R^n_+$.

\begin{Definition}
  Let $X$ be a compact manifold or an open subset of it.  Let $\{(U_i,\varphi_i,\Phi_i)\}$ be an atlas as above. Then for $-\infty < s <\infty$ and $1<p<\infty$, we define $A^{s,p}(X)$ to be those distributions $f$ on $X$ such that
  $$\|f\|_{A^{s,p}(X)} = \left(\sum_{U_i \subset \overset{\circ}{X}}\|\Phi_i^*(\varphi_if)\|_{A^{s,p}(\R^n)}^p + \sum_{U_i \cap \partial X \neq \emptyset}\|\Phi_i^*(\varphi_if)\|_{A^{s,p}(\R^n_+)}^p\right)^{1/p} < \infty.$$
We define $W^{s,p}(X)$ for $s \in \Z_+$ and $1 \leq p \leq \infty$ similarly.
\end{Definition}

If we have two different atlases, the following proposition implies that we obtain equivalent norms:

\begin{Proposition}
  Let $f \in A^{s,p}(\R^n)$, $s \in \R$ and $1 < p <\infty$. (i) If $\varphi \in C^\infty_0(\R^n)$ then $\varphi f \in A^{s,p}(\R^n)$. (ii) If $\Phi$ is a diffeomorphism of $\R^n$ which is equal to the identity outside a compact set, then $\Phi^*(f) \in A^{s,p}(\R^n)$.
\end{Proposition}

In particular, if $X$ is a bounded open subset of $\R^n$ with smooth boundary, the above furnishes a definition of $A^{s,p}(X)$.  On the other hand, we also defined $A^{s,p}(X)$ to be the restrictions to $X$ of $A^{s,p}(\R^n)$.  These two definitions of $A^{s,p}(X)$ yield equivalent norms by the extension theorem, Theorem \ref{ThmExt}.  Consequently, if $X$ is a compact manifold and $\tilde X$ is a closed manifold containing $X$, we have the following:

\begin{Proposition}
  For $-\infty<s<\infty$ and $1<p<\infty$,  $A^{s,p}(X)$ is the space of restrictions to $X$ of $A^{s,p}(\tilde X)$.
\end{Proposition}

\begin{Corollary}\label{CorDiffMap}
  Let $X$ be a compact manifold (with or without boundary) or Euclidean space.  If $D$ is a differential operator of order $m$, then $D: A^{s,p}(X) \to A^{s-m,p}(X)$ for all $s \in \R$ and $1 < p < \infty$.
\end{Corollary}

Because function spaces defined on manifolds are locally the function spaces defined on Euclidean space, many of the properties of the latter carry over to the manifold case.  For instance, if $\tilde X$ is any closed manifold containing the manifold $X$, we can define
\begin{equation}
  \tilde A^{s,p}(X) = \{f \in A^{s,p}(\tilde X) : \supp f \subseteq X\}.
\end{equation}
We have the following theorem:

\begin{Theorem}\label{ThmFunMan}
  Let $X$ be a compact manifold.  We have that $C^\infty(X)$ is dense in $A^{s,p}(X)$ and multiplication by a smooth function defines a bounded operator.  Moreover, for any $s \in \R$, $A^{-s,p'}(X)$ is the dual space of $\tilde A^{s,p}(X)$, where $p'=p/(p-1)$. If $X$ is closed or $s < 1/p$, then $\tilde A^{s,p}(X) = A^{s,p}(X)$.
\end{Theorem}

The trace theorem, Theorem \ref{ThmTrace}, readily generalizes to manifolds with boundary:

\begin{Theorem}\label{ThmTraceMan} Let $X$ be a compact manifold with boundary $\partial X$.
\begin{enumerate}
  \item For $s > m + 1/p$ and $m \in \Z_+$, the trace map (\ref{rest-m}) extends to a bounded operator
\begin{align}
  r_m: H^{s,p}(X) & \to \oplus_{j=0}^{m-1} B^{s-1/p-j}(\partial X).
\end{align}
 \item For any $s \in \R$ and $m \in \Z_+$, there exists an extension map $e_m: \oplus_{j=0}^{m-1} B^{s-1/p-j,p}(\partial X) \to H^{s,p}(X)$ such that for $s > m + 1/p$, we have $r_m e_m = \id$.
\item $H^{s,p}(X)$ may be replaced with $B^{s,p}(X)$ in the above.
\end{enumerate}
\end{Theorem}

\subsection{Further Properties}

In the following, $X$ is a compact manifold (with or without boundary).

\begin{Theorem}\label{ThmEmbed} (Embedding Theorem)
  Let $-\infty < t < s < \infty$ and $1 < p,q < \infty$ with
    \begin{equation}
    s - n/p \geq t - n/q \label{embedst}
  \end{equation}
\begin{enumerate}
  \item We have embeddings
  \begin{align}
    B^{s,p}(X) & \hookrightarrow B^{t,q}(X) \cap H^{t,q}(X) \label{Bembed}\\
    H^{s,p}(X) & \hookrightarrow H^{t,q}(X) \cap B^{t,q}(X) \label{Hembed}.
  \end{align}
  If the inequality (\ref{embedst}) is strict, these embeddings are compact.
\item We have the monotonicity property
    \begin{equation}
      H^{s,p}(X) \subseteq H^{s,q}(X), \qquad p > q.
    \end{equation}
   \item If $t > 0$ is not an integer, then
    \begin{align*}
      H^{n/p+t,p}(X) & \hookrightarrow C^t(X)\\
      B^{n/p+t,p}(X) & \hookrightarrow C^t(X).
    \end{align*}
\end{enumerate}
\end{Theorem}

Next, we have a multiplication theorem.  Namely, given two functions $f$ and $g$, we wish to know in which space their product $fg$ lies (where it is assumed that $f$ and $g$ are sufficiently regular so that their product makes sense as a distribution).

\begin{Theorem} (Multiplication Theorem) \label{ThmMult}
    \begin{enumerate}
     \item For all $s > 0$, we have $A^{s,p}(X) \cap L^\infty(X)$ is an algebra. Moreover, we have the estimate
     $$\|fg\|_{A^{s,p}} \leq C(\|f\|_{A^{s,p}}\|g\|_{L^\infty} + \|f\|_{L^{\infty}}\|g\|_{A^{s,p}}).$$
     In particular, if $s > n/p$, then $A^{s,p}(X)$ is an algebra.
     \item Let $s_1 \leq s_2$ and suppose $s_1 + s_2 > n \max(0,\frac{2}{p}-1)$. Then we have a continuous multiplication map
     \begin{align*}
       A^{s_1,p}(X) \times A^{s_2,p}(X) & \to A^{s_3,p}(X),
     \end{align*}
    where
    $$s_3 = \begin{cases}
      s_1 & \mathrm{if }\;\, s_2 > n/p\\
      s_1+s_2-n/p & \mathrm{if }\;\, s_2 < n/p.
    \end{cases}$$
   \end{enumerate}
\end{Theorem}
Both statements are standard facts, whose proofs involve the paraproduct calculus.  For (i), see e.g. \cite{Tay}.  For (ii), see \cite{RS}.\End

\section{Elliptic Boundary Value Problems}\label{AppEBP}

Let $X$ be a compact $n$-manifold with boundary ${\partial X}$ and $E$ and $F$ two vector bundles over $X$.  Let $A: \Gamma(E) \to \Gamma(F)$ be an $m$th order elliptic differential operator mapping smooth sections of $E$ to smooth sections of $F$. For all $s \in \R$ and $1 < p < \infty$, the operator $A$ extends to a map on the Triebel-Lizorkin spaces $A^{s,p}(E) \to A^{s-m,p}(E)$ by Corollary \ref{CorDiffMap}, where for $s < m$, $A$ acts in the sense of distributions.  To keep notation concrete, on $X$ we work with the Bessel potential spaces $H^{s,p}(E)$ for the moment\footnote{In addition, we want to avoid overusing the letter $A$ in our notation.}, though most of what we do applies to $B^{s,p}(E)$ as well (subsequent theorems will be stated for both the $H^{s,p}$ and $B^{s,p}$ spaces).\footnote{In most applications, one wants elliptic estimates on Sobolev spaces on the interior of the manifold, and hence one works with the $H^{s,p}$ spaces on $X$ (one still must work with Besov spaces on the boundary, since these are the space of boundary values of the $H^{s,p}$ spaces).  In this paper however, we will also work with Besov spaces on $X$ since we will need such Besov space estimates in \cite{N2}. Since the space of boundary values of a Besov space is still a Besov space, working with Besov spaces on $X$ is permissible.}

Fix a collar neighborhood $[0,\eps) \times \partial X$ of $X$, where $t\in [0,\eps)$ is the inward normal coordinate and $x$ denotes the coordinates on $\partial X$. In this neighborhood, write the principal part of $A$ as $\sum_{j=0}^m A_j\partial_t^{m-j}$ where $A_j = A_j(x,t)$ are differential operators of degree $j$ in the tangential variables.  Let $(x,\xi) \in T^*{\partial X} \setminus \{0\}$.  Consider the vector space of solutions $f: \R^+ \to \mathbb{C}$ to the ordinary differential equation
\begin{equation}
  \left(\sum_{j=0}^m A_j(x,0,\xi)\partial_t^{m-j}\right)f(t) = 0, \qquad t \in \R, \label{LSode}
\end{equation}
obtained by ``freezing" $A$ at $(x,0,\xi)$.  Here, $A_j(x,0,\xi)$ is the symbol of $A_j$ at $t = 0$.  Let $M^{\pm}_{x}(\xi)$ denote the vector space of solutions to (\ref{LSode}) which decay exponentially as $t \to \pm\infty$.  The assumption that $A$ is elliptic implies that the solution space of (\ref{LSode}) decomposes as a direct sum $M^+_{x}(\xi) \oplus M^-_{x}(\xi)$, for all $(x,\xi')\in T^*\partial X\setminus 0$.  Thus, we have an isomorphism $M^+_{x}(\xi) \oplus M^-_{x}(\xi) \cong E_x^m$ given by taking the full Cauchy data of a solution, $f(t) \mapsto (f(0), \ldots, \partial_t^{m-1}f(0))$.  Via this isomorphism, we can identify $M^{\pm}_x(\xi) \subset E_x^m$.

\begin{Definition}\label{Defpi}
  For $(x,\xi) \in T^*\partial X\setminus 0$, define $\pi^+_A(x,\xi): E_x^m \to E_x^m$ to be the projection onto $M^+_x(\xi)$ through $M^-_x(\xi)$.
\end{Definition}

From this projection, we can define what it means for a boundary condition $B$ to be elliptic.  Such a notion is very classical; for further reading, see \cite{Ho}, \cite{WRL}.

Suppose the operator $B$ is given by a tuple $B = (B_1,\ldots, B_\ell)$ of operators on the boundary, where
\begin{align}
 \hspace{1in}B_k: \Gamma(E_{\partial X})^m &\to \Gamma(V_k) \nonumber \\
 B_kU &:= \sum_{j=0}^{m-1}B_{kj}U_j, \qquad U = (U_j)_{j=0}^{m-1} \in \Gamma(E_{\partial X})^m,\; 1 \leq k \leq \ell, \label{Bops}
\end{align}
and where the $B_{kj}$ are (classic) pseudodifferential operators mapping $\Gamma(E_{\partial X})$ to $\Gamma(V_k)$ for some vector bundle $V_k$ over ${\partial X}$.  The total boundary operator gives us a map
\begin{equation}
  B: \Gamma(E_{\partial X})^m \to \oplus_{k=1}^\ell \Gamma(V_k). \label{BC}
\end{equation}

Given a pseudodifferential operator $T$, let $\sigma_p(T)$ denote the principal symbol of $T$.  If $T$ is a matrix of  pseudodifferential operators (where different entries of the matrix correspond to different vector bundles), then the order of $T$ and hence its principal symbol need to be carefully defined.  In the case of $B$ above, since each $B_k$ represents a boundary condition of order $\beta_k$, we define the principal symbol of $B$ to be the following symbol-valued matrix:
$$\sigma_p(B) = (\sigma_{\beta_k - j}B_{kj})_{1 \leq k \leq \ell, \atop 0 \leq j \leq m-1},$$
where $\sigma_i(B_{kj})$ is the usual principal symbol of $B_{kj}$ if $B_{kj} \in Op^i(\partial X)$ and zero if $B_{kj} \in Op^{i'}(\partial X)$, $i' < i$.

\begin{Def}
  Suppose the boundary operator $B$ in (\ref{BC}) is such that $\sigma_p(B)(x,\xi): (E_x)^m \to \oplus_{k=1}^\ell V_k$ restricted to $\im \pi_A^+(x,\xi)$ is an isomorphism onto $\im \sigma_p(B)(x,\xi)$ for all $(x,\xi) \in T^*\partial X\setminus 0$.  Then $B$ is an \textsl{elliptic boundary condition} for the operator $A$.  In this case, we say that the pair $(A,B)$ is elliptic.
\end{Def}

\noindent\textit{Examples. }Standard examples of elliptic boundary conditions include the Dirichlet and Neumann boundary conditions for the Laplacian.  Likewise, for Dirac operators, the APS boundary condition (which is the positive spectral projection for the tangential boundary operator, see Definition \ref{DefTBO}) and other similar pseudodifferential boundary conditions are elliptic.  A natural example of a boundary condition of mixed order is the following situation, which occurs, e.g., in \cite{SaWe}.  Let $X = Y$ be a $3$-manifold and let $E = T^*Y \oplus \R$. Let $A = D_{\dgc}$ be the div-grad-curl operator (\ref{DGCop}) acting on $\Gamma(E) =  \Omega^1(Y) \oplus \Omega^0(Y)$.  Adopting the notation and identifications made in Lemma \ref{LemmaDGC}, we have $\Gamma(E_\Sigma) = \Omega^1(\Sigma) \oplus \Omega^0(\Sigma) \oplus \Omega^0(\Sigma)$, and we can define the boundary operator
\begin{equation}
  B = d_\Sigma \oplus \id \oplus 0: \Gamma(E_\Sigma) \to \Omega^2(\Sigma) \oplus \Omega^0(\Sigma),
\end{equation}
One can check that $B$ is a elliptic boundary condition for $D_\dgc$.  It is in fact a local boundary condition, i.e., $B$ is a differential operator.  One can replace the operator $d_\Sigma: \Omega^1(\Sigma) \to \Omega^2(\Sigma)$ with a projection $\pi_L: \Omega^1(\Sigma) \to \Omega^1(\Sigma)$ whose kernel is a Lagrangian subspace $L$ of $\Omega^1(\Sigma)$, where the symplectic form on $\Omega^1(\Sigma)$ is the form $\omega(a,b) = \int_\Sigma a \wedge b$.  Since $\ker \pi_L$ and $\ker d$ differ by a finite dimensional subspace, the operator $B_L = \pi_L \oplus \id \oplus 0$ also induces an elliptic boundary condition for $D_{\dgc}$.  The linear operators considered in \cite{SaWe} are based on this example, where various Lagrangians $L$ are considered.\\

On a manifold with boundary, the map $A: H^{s+m,p}(E) \to H^{s,p}(F)$ has infinite dimensional kernel.  An elliptic boundary condition $B$ allows us to control the kernel, which means we can obtain an elliptic estimate for the full mapping pair $(A,B)$.  Moreover, if we consider the restricted operator $A_B$, whose domain consists of those elements annihilated by the boundary condition, the elliptic estimate for the full mapping pair $(A,B)$ then gives us one for $A_B$.

More precisely, given $0 \leq \kappa \leq m - 1$, let $s + m > \kappa + 1/p$ and $1<p<\infty$, and consider the map
\begin{equation}
  r_\kappa: H^{s+m,p}(E) \to \oplus_{j=0}^\kappa B^{s+m-1/p-j,p}(E_{\partial X}), \label{rest-kappa}
\end{equation}
the trace map onto the Cauchy data up to order $\kappa$.  If $B$ is an elliptic boundary condition that depends only on the Cauchy data up to order $\kappa$, i.e., we have $B_{kj} = 0$ in (\ref{Bops}) for $j > \kappa$, then
\begin{equation}
B: \oplus_{j=0}^{\kappa} B^{s+m-1/p-j,p}(E_{\partial X}) \mapsto \oplus_{k=1}^\ell B^{s+m-1/p-\beta_k,p}(V_k)  \label{BC2}
\end{equation}
is bounded. For brevity, for any $s \in \R$, let
\begin{align}
  \cV^{s-1/p,p}_{\bar\beta} &= \oplus_{k=1}^\ell B^{s-1/p-\beta_k,p}(V_k), \qquad \bar\beta = (\beta_1,\ldots,\beta_k).
\end{align}

We have the following fundamental theorem for elliptic boundary value problems. In our applications, we will only consider $m = 1,2$ and $\kappa = 0,1$, see Corollary \ref{CorEBP}.

\begin{Theorem}\cite[Theorem 4]{Se2} \label{ThmEBP} Let $X$ be a compact manifold with boundary $\partial X$ and let $A: \Gamma(E) \to \Gamma(F)$ be an $m$th order elliptic differential operator.  Suppose $B$ is an elliptic boundary condition satisfying (\ref{BC2}) which depends only on the Cauchy data up to order $\kappa$. Let $1<p<\infty$ and $s + m > \kappa + 1/p$.
\begin{enumerate}
\item Let $u \in H^{t,p}(E)$, $t \in \R$, and suppose $Au \in H^{s,p}(F)$.  Then $r_\kappa u \in \oplus_{j=0}^\kappa B^{t-1/p-j,p}(E)$ is a well-defined element.  Suppose we have the additional estimate $Br_\kappa u \in \V^{s+m-1/p,p}_{\bar\beta}$ on the boundary. Then $u \in H^{s+m,p}(E)$ and
\begin{equation}
\|u\|_{H^{s+m,p}(E)} \leq C(\|Av\|_{H^{s,p}(F)} + \|Br_\kappa u\|_{\V^{s+m-1/p,p}_{\bar\beta}} + \|u\|_{H^{t,p}(E)}). \label{Ereg}
\end{equation}
\item The map $A_B: \{u \in H^{s+m,p}(E): Br_\kappa u = 0\} \to H^{s,p}(E)$ is Fredholm.  Its kernel and cokernel are spanned by finitely many smooth sections.
\item If $\sigma_p(B): E^m \to \oplus_{k=1}^\ell V_k$ is surjective, then the full mapping pair
\begin{align*}
  (A,B): H^{s+m,p}(E) & \to H^{s,p}(F) \oplus \V^{s+m-1/p,p}_{\bar\beta}\\
  u & \mapsto (Au, Br_\kappa u)
\end{align*}
is a Fredholm operator.
\item The above statements remain true if $H^{\bullet,p}(E)$ is replaced with $B^{\bullet, p}(E)$.
\end{enumerate}
\end{Theorem}

\begin{Rem}\label{RemKerEBP}
  By a standard argument, the lower order term $\|u\|_{H^{t,p}(E)}$ in (\ref{Ereg}) can be replaced with $\|\pi u\|$, where $\pi$ is any projection onto the finite dimensional space of solutions to $Au = Br_\kappa u = 0$, and $\|\cdot\|$ is any norm on that space.  In other words, we need only control the kernel of the operator $(A,B)$ to get the estimate (\ref{Ereg}).  In particular, if $(A,B)$ has no kernel, the term $\|u\|_{H^{t,p}(E)}$ can be omitted.
\end{Rem}

\begin{Rem}\label{RemEBP}
  Seeley only proves Theorem \ref{ThmEBP} for $p=2$, but because all the maps involved are pseudodifferential or involve taking traces or extensions, the generalization to $1<p<\infty$ is automatic.  Furthermore, the generalization to Besov spaces in (iv) follows from the general results in the previous section on function spaces, since the trace of a Besov space is still a Besov space, and the boundary extension map of a Besov space lies not only in a Bessel potential space but also in a Besov space (see Theorem \ref{ThmTrace}).  This is another instance of the principle that $B^{s,p}$ and $H^{s,p}$ are ``nearly identical".
\end{Rem}

We will not need Theorem \ref{ThmEBP} in its full generality.  For convenience, we summarize the particular applications we have in mind below:

\begin{Corollary}\label{CorEBP}
We have the following elliptic boundary value problems:
\begin{enumerate}
    \item Let $A = \Delta$ be the Laplacian acting on scalar functions.  Then the Dirichlet and Neumann boundary conditions are elliptic boundary conditions.  For the Dirichlet problem, we have the elliptic estimate
  \begin{align}
    \|u\|_{B^{s+2,p}(X)} & \leq C(\|\Delta u\|_{B^{s,p}(X)} + \|r_0u\|_{B^{s+2-1/p,p}(\partial X)})
  \end{align}
  for $s + 2 > 1/p$. For the Neumann problem, we have the elliptic estimate
  \begin{align}
    \|u\|_{B^{s+2,p}(X)} & \leq C\left(\|\Delta u\|_{B^{s,p}(X)} + \|r_0(\partial_\nu u)\|_{B^{s+1-1/p,p}(X)} + \left|\int_X u\right |\right)
  \end{align}
  for $s + 2 > 1+1/p$, where $\partial_\nu u$ denotes the derivative of $u$ with respect to the outward unit normal to $\partial X$.
\item Let $A = d+d^*$ be the Hodge operator acting on $\oplus_{i=0}^n\Omega^i(X)$, the exterior algebra of differential forms on $X$.  Then the tangential component\footnote{See footnote \ref{footnotedf}.} $a|_{\partial X}$ and normal component $*a|_{\partial X}$ are elliptic boundary conditions.  In particular, if $a \in \Omega^1(X)$, then we have the elliptic estimate
    \begin{equation}
      \|a\|_{B^{s+1,p}\Omega^1(X)} \leq C(\|da\|_{B^{s,p}\Omega^2(X)} + \|d^*a\|_{B^{s,p}\Omega^0(X)} + \|a^h\|_{B^{s,p}\Omega^1(X)})
    \end{equation}
    for $s+1>1/p$, where $a^h$ denotes the orthogonal projection of $a$ onto the space
    \begin{equation}
      H^1(X,\partial X; \R) \cong \{a \in \Omega^1(X): da = d^*a = 0, a|_{\partial X} = 0\}. \label{H1X}
    \end{equation}
\item Let $A: \Gamma(E) \to \Gamma(F)$ be a Dirac operator.  If $B$ is any pseudodifferential projection onto $r(\ker A)$, then $B$ is an elliptic boundary condition for $A$.  We have the elliptic estimate
    \begin{equation}
      \|u\|_{B^{s+1,p}(E)} \leq C(\|Au\|_{B^{s,p}(F)} + \|Bru\|_{B^{s+1-/1p,p}(E_{\partial X})}). \label{DiracEBP}
    \end{equation}
    for $s +1 > 1/p$.
\end{enumerate}
\end{Corollary}

\Proof For (i), a standard computation shows that the kernel of the Dirichlet Laplacian is zero and the kernel of the Neumann Laplacian is spanned by constant functions.  We now apply Remark \ref{RemKerEBP}.  For (ii), the kernel of $d+d^*$ on $\Omega^1(Y)$ with vanishing tangential component is the space (\ref{H1X}).  We now apply the previous theorem and Remark \ref{RemKerEBP}.  Observe that for the Dirichlet Laplacian we took $\kappa = 0$ in Theorem \ref{ThmEBP}. For (iii), there is no term to account for the kernel due to unique continuation, Theorem \ref{ThmUCP}, which implies that $r$ maps $\ker A$ isomorphically onto its image (hence $Br$ maps $\ker A$ isomorphically onto its image).\End

\subsection{The Calderon Projection}\label{AppSeeley}

Our main theorem for elliptic boundary value problems, Theorem \ref{ThmEBP}, is quite strong (in fact stronger than most variants of the theorem in the literature) due to its very large range of admissible parameters\footnote{Namely, in Theorem \ref{ThmEBP}, one can work with very low, even negative regularity (i.e. $s < 0$), and furthermore, one can work with just a subset of the full Cauchy data (i.e. $\kappa < m - 1$) in specifying boundary conditions.}.  This theorem, due to Seeley, uses the construction of  special pseudodifferential operators associated to elliptic boundary value problems which we now describe. As with Remark \ref{RemEBP}, all statements in this section remain true if we replace $H^{\bullet,p}(E)$ with $B^{\bullet,p}(E)$.

Let $A: \Gamma(E) \to \Gamma(F)$ be an $m$th order elliptic operator, which for simplicity, we take to be first order, though everything we discuss here generalizes straightforwardly for $m > 1$.  Informally, the general picture is the following. We have two subspaces of interest, $\ker A$ and its restriction to the boundary $r(\ker A)$,
where $r: \Gamma(E)\to\Gamma(E_\Sigma)$ is the restriction map. What we have is that there exist a pseudodifferential operator $P^+: \Gamma(E_\Sigma) \to \Gamma(E_\Sigma)$ acting on boundary sections and a map $P: \Gamma(E_\Sigma) \to \Gamma(E)$ mapping boundary sections into the interior such that $P^+$ is a projection onto $r(\ker A)$ and the range of $P$ is contained in $\ker A$.  Furthermore, we have $rP = P^+$.

More precisely, and assigning the appropriate topologies to the spaces involved, let $s \in \R$ and $1<p<\infty$, and let
\begin{equation}
  Z^{s,p}(A) \subset H^{s,p}(E)
\end{equation}
be the kernel of the operator $A: H^{s,p}(E) \to H^{s-1,p}(E)$.  Let $Z_0(A)$ be the subset of $Z^{s,p}(A)$ consisting of those elements $z$ with vanishing boundary values, i.e., $r(z) = 0$.  Theorem \ref{ThmEBP} implies $Z_0(A) \subset C^\infty(E)$ and is finite dimensional. The map $r$ extends to a bounded map $H^{s,p}(E) \to B^{s-1/p,p}(E_\Sigma)$ only when $s > 1/p$. However, if we restrict $r$ to the kernel of $A$, it turns out that no such restriction on $s$ is necessary.  This is the content of the following very important theorem:

\begin{Theorem}\cite{Se, Se2} \label{ThmSeeley}
  Let $s \in \R$ and $1<p<\infty$.
  \begin{enumerate}
    \item We have a bounded map $r: Z^{s,p}(A) \to B^{s-1/p,p}(E_{\partial X})$, and furthermore, its range is closed.  In particular, if $Z_0(A) = 0$, then $r$ is an isomorphism onto its image.
    \item There is a pseudodifferential projection $P^+$ which projects $B^{s-1/p,p}(E_{\partial X})$ onto $r(Z^{s,p}(A))$. Furthermore, the principal symbol $\sigma_0(P^+)$ of $P^+$ is equal to the symbol $\pi_A^+$ (see Definition \ref{Defpi}).
    \item There is a map $P: B^{s-1/p,p}(E_{\partial X}) \to Z^{s,p}(A)$ whose range has $Z_0(A)$ as a complement.  Furthermore, $PP^+ = P$ and $rP = P^+$.
  \end{enumerate}
\end{Theorem}

Thus, in particular, the above theorem tells us that elements in the kernel of $A$ of any regularity have well-defined restrictions to the boundary.  In fact, the first part of Theorem \ref{ThmEBP}(i) relies crucially on this fact.

\begin{Def}\label{DefCP}
  The operators $P^+$ and $P$ in Theorem \ref{ThmSeeley} are called a \textit{Calderon projection} and \textit{Poisson operator} of $A$, respectively.
\end{Def}

\begin{Rem}\label{RemCP}
  (i) From the definitions, it follows that $P^+$ is an elliptic boundary condition for $A$. (ii) A projection is defined not only by its range but also by its kernel.  Thus, we have \textit{a} Calderon projection and Poisson operator for $A$, since their kernels are not uniquely defined.  When we speak of these operators then, we usually have a particular choice of these operators in mind.  Seeley, for instance, has a particular construction of $P^+$ and $P$.  However, it is usually only the range of $P$ and $P^+$ that are of main interest to us, and these are uniquely specified by the above definitions.  Hence, a Calderon projection is often times referred to as \textit{the} Calderon projection in the literature.
\end{Rem}

Altogether, $P^+$ is a projection onto the Cauchy data of the kernel of $A$, and $P$ is a map from the Cauchy data of the kernel into the kernel.  The latter map is an isomorphism when $Z_0(A) = 0$.  Furthermore, we have

\begin{Corollary}\label{CorKerComp}
  For all $s \in \R$, smooth configurations are dense in $Z^{s,p}(E)$.  Furthermore, suppose $s > 1/p$.  Then $Z^{s,p}(E) \subset H^{s,p}(E)$ is complemented.  Moreover, if $Z_0(A) = 0$, then $Pr: H^{s,p}(E) \to Z^{s,p}(A)$ is a bounded projection onto $Z^{s,p}(A)$.
\end{Corollary}

\Proof We have that $Z^{s,p}(A)$ is the direct sum of $Z_0(A)$ and the image of $P: B^{s-1/p,p}(E_{\partial X}) \to H^{s,p}(E)$.  The first statement now follows since the space $Z_0(A)$ is spanned by smooth sections and smooth sections are dense in $B^{s-1/p,p}(E_{\partial X})$.  Now consider $s > 1/p$.  Then the map $Pr: H^{s,p}(E) \to Z^{s,p}(A)$ is a projection onto the image of $P$ (since $PrP = P$ by Theorem \ref{ThmSeeley}), which is of finite codimension in $Z^{s,p}(A)$. From this, one can construct a projection of $H^{s,p}(E)$ onto $Z^{s,p}(A)$, which means $Z^{s,p}(A)$ is a complemented subspace.  If $Z_0(A) = 0$, then the range of $P$ is all of $Z^{s,p}(A)$, whence $Pr$ is a projection onto $Z^{s,p}(A)$.\End

We present an important application of these operators. Let $A$ be a first order formally self-adjoint elliptic operator.  Then the operator $J:= A_0$ in (\ref{LSode}) is a skew-symmetric automorphism on the boundary, and Green's formula (\ref{genGF}) for $A$ defines for us a symplectic form
$$\omega(u,v) = \int_\Sigma (u,-Jv)$$
on boundary sections $u,v \in \Gamma(E_{\partial X})$.  This symplectic form extends to a well-defined symplectic form on $B^{s,p}(E_{\partial X})$ for $(s,p) = (0,2)$ and for $s > 0$, $p \geq 2$, and the map $-J$ is a compatible complex structure with respect to this symplectic form.  Indeed, for this range of $s$ and $p$, we have $B^{s,p}(E_{\partial X}) \hookrightarrow L^2(E_{\partial X})$, with the latter a strongly symplectic Hilbert space.

We say that a closed subspace of $B^{s,p}(E_{\partial X})$ is Lagrangian if it is isotropic with respect to $\omega$ and it has an isotropic complement.  Observe that if $L \subset L^2(E_{\partial X})$ is Lagrangian, then $JL$ is a Lagrangian complement of $L$.

\begin{Proposition}\label{PropSeLag}\cite{BBLZ}
  Let $A$ be a Dirac operator.  Then $\im P^+$ and $J\im P^+$ are complementary Lagrangian subspaces of $B^{s,p}(E_{\partial X})$, where $(s,p) = (0,2)$ or $s > 0$, $p \geq 2$.
\end{Proposition}

\Proof In \cite{BBLZ}, it is shown that $\im P^+$ and $J\im P^+$ define complementary Lagrangian subspaces of $L^2(E_{\partial X})$.  Here, it is essential that one uses the trick of constructing an ``invertible double" for the operator $A$. Since $P^+$ is a pseudodifferential projection, it is bounded on $B^{s,p}(E_{\partial X})$.  Without loss of generality, we can suppose $P^+$ is an orthogonal projection (making a projection into an orthogonal projection preserves the property of being pseudodifferential).  Define $P^- = JP^+J^{-1}$.  Then $\im P^- = J\im P^+$ and its principal symbol agrees with the principal symbol of $1 - P^+$.  It follows that $\im P^+ \oplus J\im P^+$ is a closed subspace of $B^{s,p}(E_{\partial X})$ of finite codimension.  We now apply Lemma \ref{LemmaDenseFred}, which tells us that $\im P^+ \oplus J\im P^+$ is in fact all of $B^{s,p}(E_{\partial X})$.\End

\section{Unique Continuation}\label{AppUCP}

Let $A: \Gamma(E) \to \Gamma(F)$ be a smooth Dirac operator acting between sections of the Hermitian vector bundles $E$ and $F$ over a compact manifold $X$ (with or without boundary). The operator $A$ is said to obey the unique continuation property if every $u$ that solves $Au = 0$ and which vanishes on an open subset of $X$ vanishes identically.  It is well-known that Dirac operators obey the unique continuation property. If $X$ is a manifold with boundary, we can replace the condition that $u$ vanish on an open set with the condition $ru = 0$, where $r = r_0$ is the restriction map to the boundary. This is because one can extend the operator $A$ to a Dirac operator $\tilde A$ on an open manifold $\tilde X$ that contains $X$ in its interior, and one can extend $u$ to $\tilde X$ by zero outside of $X$.  Since $\tilde A$ is a first order operator, then $\tilde A \tilde u = 0$ on $\tilde X$.  Since $\tilde u$ vanishes on an open set, then $\tilde u \equiv 0$ on $\tilde X$ and so $u \equiv 0$ on $X$.

The following is a well-known general result (see \cite[Chapter 8]{BBW}):

\begin{Theorem}\label{ThmUCP}
  Let $X$ be a compact manifold with boundary, let $D$ be a smooth Dirac operator on $\Gamma(E)$, and let $V$ be an $L^\infty$ multiplication operator.  Then $D+V$ has the unique continuation property.  More precisely, if $u \in B^{1,2}(E)$ satisfies $(D+V)u = 0$ and $ru = 0$, then $u \equiv 0$.\footnote{One can start with $u$ of lower regularity than $B^{1,2}(E)$, say $L^2(E)$, since by elliptic bootstrapping, such a $u$ will necessarily be of regularity $B^{1,2}(E)$, see the proof of \ref{ThmUCPsurj}.}
\end{Theorem}

One application of this theorem is to show that such an operator $D+V$, acting between suitable function spaces, is surjective on a manifold with boundary.  This is in contrast to when $X$ is closed, in which case $D+V$ is only Fredholm, in which case it may have a finite dimensional cokernel. We have the following theorem:

\begin{Theorem}\label{ThmUCPsurj}
  Let $X$ be a compact manifold with boundary.  Let $2 \leq p <\infty$, $s > 1/p$ and let $D + V: B^{s,p}(E) \to B^{s-1,p}(F)$ where $V$ is a sufficiently smooth\footnote{To keep the function space arithmetic simple, we suppose $V$ is smooth in the proof in the theorem, though the necessary modifications can be made for $V$ nonsmooth but bounded as a map between suitable function spaces, depending on $s$,$p$. What mainly needs to carry through is the bootstrapping argument in (\ref{bootstrap}).  In all applications, we will always have $V \in B^{t,p}(Y)$ where $t$ is sufficiently large so that the statement remains true with $V$ of this regularity class.  If $s \geq 1$, one can check that $V \in L^\infty(E)$ suffices.  If $s < 1$, one wants $V$ to have some regularity so that it can act via multiplication on functions of low regularity.} multiplication operator.  Then $D+V$ is surjective.
\end{Theorem}

\Proof Since $D+V$ is a smooth elliptic operator, it has a right (pseudodifferential) parametrix.  This shows that $D+V$ has closed range and finite dimensional cokernel.  It remains to show that the cokernel is zero. There are two cases to distinguish, the cases $s > 1$ and $s \leq 1$.  Let us deal with the latter case, with the case $s > 1$ similar.  Suppose $u \in (B^{s-1,p}(F))^* = B^{1-s,p'}(F)$, $p' = p/(p-1)$, belongs to the dual space of $B^{s-1,p}(F)$ and annihilates $\im(D+V) \subseteq B^{s-1,p}(F)$. We want to show that $u = 0$, which combined with the fact that $\im(D+V)$ is closed means that $\im(D+V)$ is all of $B^{s-1,p}(F)$.  The condition that $u$ annihilate $\im(D+V)$ means that we have the (weak) equation $(D + V)^*u = 0$, and thus, $Du = -V^*u$ (here we think of dual operator $D^*$ acting on the linear functional $u$ as being the same as $D$, since a Dirac operator is formally self-adjoint).  We have $V^*u\in B^{1-s,p'}$, since multiplication by a smooth function is bounded on all Besov spaces. By Theorem \ref{ThmEBP}(i), we have a well-defined trace $r(u) \in B^{1-s-1/p',p'}(F_{\partial X})$.  Thus, for all $v \in B^{s,p}(E)$, we have Green's formula (\ref{genGF2}), which tells us that
\begin{align}
  0 &= (v,(D+V)u) - ((D+V^*)v, u) \nonumber \\
  & = \int_{\partial X}(r(v),-Jr(u)). \label{GFkill}
\end{align}
The first line follows since $u$ annihilates $\im(D+V)$ and $(D+V^*)u = 0$.  The second line is well-defined since $Jr(v) \in B^{s-1/p,p}(F)$ and $B^{s-1/p,p}(F)$ is the dual space of $B^{-s+1/p,p'}(F) = B^{1-s-1/p',p'}(F)$.  Since (\ref{GFkill}) holds for all $v \in B^{s,p}(E)$, it follows that $r(u) = 0$.  The system $(D+V^*)u = 0$ and $r(u) = 0$ is overdetermined which means that we have an elliptic estimate for $u$ via Theorem \ref{ThmEBP}.  That is, since $Du = -V^*u$, we have an estimate of the form
\begin{equation}
  \|u\|_{B^{t+1,q}} \leq C(\|V^*u\|_{B^{t,q}} + \|u\|_{B^{t,q}}). \label{bootstrap}
\end{equation}
for all $t,q$ such that the right-hand side is finite, $t + 1 > 1/q$. Since $u, V^*u \in B^{1-s,p'}(E)$ we have $u \in B^{2-s,p'}(E)$, where $2-s > 1/p'$ since $s \leq 1$.  Feeding this back into (\ref{bootstrap}) and using that $V$ is smooth, we see that we can boostrap $u$ to any desired regularity.  Thus $u$ is smooth.  (In general, for $V$ not smooth, we want $V$ sufficiently regular so that the above steps allow us to bootstrap to $u \in B^{1,2}(E)$).  Furthermore, $r(u) = 0$.  We now apply Theorem \ref{ThmUCP} to conclude $u = 0$.  Thus, $D+V$ is surjective.\End

The next two unique continuation theorems are ones that are specific to the Seiberg-Witten equations. In \cite{KM}, unique continuation theorems are proved for the (linearized and nonlinear) $4$-dimensional Seiberg-Witten equations.  Because solutions to the three-dimensional Seiberg-Witten equations on $Y$ can be regarded as time-independent solutions to the four-dimensional Seiberg-Witten equations on $S^1\times Y$ (see \cite{KM}), these methods carry over to yield unique continuation for the three-dimensional equations. More precisely, if $(B(t),\Psi(t))$ is a path of configurations in $\fC(Y)$, then $(B(t),\Psi(t))$ solve the Seiberg-Witten equations on $S^1 \times Y$ if and only if
\begin{align*}
 \frac{d}{dt}(B(t),\Psi(t)) = - SW_3(B(t),\Psi(t)).
\end{align*}
Thus, if $(B(t),\Psi(t)) \equiv (B,\Psi)$ is independent of time and satisfies $SW_3\c = 0$, then $(B(t),\Psi(t))$ satisfies the Seiberg-Witten equations on $S^1 \times Y$ and we may apply the methods of \cite{KM}.  We sketch proofs for the below results; a more detailed proof will be left to \cite{N0}.

\begin{Theorem}\label{ThmUCPnonlin}
  Assume (\ref{assum}) and $s > \max(3/p,1/2)$.  If $\cx[1],\cx[2] \in \fM^{s,p}$ are irreducible and satisfy $r_\Sigma\cx[1] = r_\Sigma\cx[2]$, then $\cx[1]$ and $\cx[2]$ are gauge equivalent on $Y$.
\end{Theorem}

\textbf{Proof (Sketch)}\hspace{0.3cm} Regard $\cx[1]$ and $\cx[2]$ as solutions to the $4$-dimensional Seiberg-Witten equations on $S^1 \times Y$ as above.  If we can apply \cite[Proposition 7.2.1-2]{KM}, then we will be done, since \cite[Proposition 7.2.1]{KM} implies $\cx[1]$ and $\cx[2]$ are gauge-equivalent on a tubular neighborhood $[0,1] \times (S^1 \times \Sigma)$ of the boundary of $S^1\times Y$, and then \cite[Proposition 7.2.2]{KM} tells us that once the irreducible $\cx[1]$ and $\cx[2]$ are gauge equivalent on an open set, then they are gauge equivalent on all of $S^1 \times Y$.  Restricting back to $Y$ yields the desired result.  In both these propositions, the unique continuation method of \cite[Lemma 7.1.3]{KM} is used to show that $(b(t),\psi(t)) := \cx[1] - \cx[2]$ is zero on $S^1 \times Y$, where in our case, $(b(t),\psi(t))$ is independent of $t \in S^1$.  One can check that the analysis used in \cite[Lemma 7.1.3]{KM} works for configurations that belong to $L^\infty(S^1\times Y) \cap B^{1,2}(S^1 \times Y)$ (i.e., we do not need $C^2(S^1 \times Y)$ smoothness).  Thus, if we can show that $(b,\psi) \in B^{1,2}(S^1\times Y)$, we will be done.  Since $(b,\psi)$ is independent of time, it suffices to show that $(b,\psi) \in B^{1,2}(Y)$.  For this, we use that $(b,\psi)$ satisfies an equation of the form $\tH_\c(b,\psi) = (b,\psi)\#(b,\psi)$, since $\cx[1]$ and $\cx[2]$ are monopoles in Coulomb-gauge.  Moreover, $r_\Sigma(b,\psi) = 0$.  Thus, as in the proof of Theorem \ref{ThmUCPsurj}, we can bootstrap the regularity of $(b,\psi) \in \T^{s,p}$ to $(b,\psi) \in \T^{1,2}$ since we have an overdetermined elliptic boundary value problem for $(b,\psi)$.  This proves the result.\End

The next result is essentially a linear version of the previous theorem.  Observe that if $(B,\Psi) \in \fM^{s,p}(Y)$ with $s > 3/p$ then $(B,\Psi) \in L^\infty\fC(Y)$ and furthermore, $(B,\Psi) \in C^\infty_{\loc}\fC(Y)$, where $C^\infty_{\loc} = C^\infty_{\loc}(Y)$ is the space of functions which belong to $C^\infty(K)$ for every compact subset $K$ contained in the interior of $Y$.  Indeed, the Seiberg-Witten equations in Coulomb-gauge are elliptic in the interior of $Y$ and so we can bootstrap the regularity of $(B,\Psi)$ to any desired regularity in the interior.  However, since we do not have stronger control of $(B,\Psi)$ at the boundary, this regularity can not in general be boostrapped to all of $Y$.

Thus, for $(B,\Psi) \in \fM^{s,p}(Y)$, we have the map
\begin{align}
  \bd_{\c,t}: \{\xi \in B^{2,2}\Omega^0(Y; i\R): \xi|_\Sigma = 0\} & \to \T^{1,2}_\loc \nonumber \\
  \xi & \mapsto (-d\xi,\xi\Psi), \label{Jloc}
\end{align}
whose image is the formal tangent space to the gauge orbit of $\c$ determined by the gauge group $\G_{\partial}^{2,2}$.  Here, $\T^{1,2}_\loc$ is the closure of the space $\T$ in the $B^{1,2}_\loc(Y)$ topology, where the subscript ``$\loc$" has the same meaning as above. Denote the image of (\ref{Jloc}) by $\J^{1,2,\loc}_{\c,t}$.

\begin{Theorem}\label{ThmUCP2}
  Let $\c \in \fM^{s,p}(Y)$, $s > 3/p$.  Suppose $(b,\psi) \in \T^{1,2}$ satisfies $\H_\c(b,\psi) = 0$ and $r_\Sigma(b,\psi) = 0$.  Then either (i) $(b,\psi) \in \J^{1,2,\loc}_{\c,t}$ or else (ii) $\Psi \equiv 0$, and then $\psi \equiv 0$ and $b \in \ker d$.
\end{Theorem}

\begin{Corollary}\label{CorUCP}
  Let $\c \in \fM^{s,p}(Y)$, $s > 3/p$.  Suppose $(b,\psi) \in \K^{1,2}_\c$ satisfies $\H_\c(b,\psi) = 0$ and $r_\Sigma(b,\psi)=0$.  Then either (i) $(b,\psi) = 0$ or else (ii) $\Psi \equiv 0$ and $b \in H^1(Y,\Sigma; i\R) \cong \{a \in \Omega^1(Y; i\R) : da = d^*a = 0, a|_\Sigma = 0\}$.
\end{Corollary}

\Proof This immediately follows from the previous theorem and $\J^{1,2,\loc}_{\c,t} \cap \K^{1,2}_\c = 0$.\En


\section*{Index of Notation}

\pagestyle{plain}

\setlength{\extrarowheight}{.2cm}
\addcontentsline{toc}{section}{\numberline{}Index of Notation}

\begin{tabularx}{450pt}{l|X}
$\A(X)$ & The space of $\spinc$ connections on the manifold $X$, p. \pageref{p:cA}.\\
$B$ & A typical $\spinc$ connection, usually on a $3$-manifold, p. \pageref{p:B}.\\
$B^{s,p}$ & A Besov space. Used as a prefix, it denotes closure with respect to said topology, p. \pageref{p:funspace}.\\
$\C$ & A Coulomb slice in $\T$, p. \pageref{p:C}. \\
$\fC(X)$ & The smooth configuration space of $\spinc$ connections and spinors on the manifold $X$, p. \pageref{p:fC}.\\
$\fC^{s,p}(X)$ & The $B^{s,p}(X)$ configuration space on $X$, p. \pageref{p:Csp}. \\
$\bd_\bullet$ & The operator associated to the infinitesimal action of the gauge group at $\bullet \in \fC(X)$, p. \pageref{p:bd}.\\
$\bd_\bullet^*$ & The formal adjoint of $\bd_\bullet$, p. \pageref{p:bd^*}.\\
$E_{\partial X}$ & Abbreviation for the restriction of a bundle $E$ on $X$ to $\partial X$.\\
$E_\bullet$ & A chart map for $\bullet \in \fM$, $\M$, or $\L$, p. \pageref{p:Eco}, \pageref{p:Erco}.\\
$E_\bullet^1$ & The nonlinear part of the chart map $E_\bullet$, p. \pageref{p:E1co}, \pageref{p:E1rco}.\\
$F_\c$ & A local straightening map for $\fM$ at $\c$, p. \pageref{p:Fco}.\\
$F_{\Sigma, \c}$ & A local straightening map for $\fL$ at $r_\Sigma\c$, p. \pageref{p:FSigma}.\\
$\G(X)$ & The gauge group of transformations on $X$, p. \pageref{p:G}.\\
$\G_\partial(X)$ & The gauge group of transformations that is the identity on $\partial X$, p. \pageref{p:Gpartial}.\\
$H^{s,p}$ & A Bessel potential space (otherwise known as fractional Sobolev spaces.) Used as a prefix, it denotes closure with respect to said topology, p. \pageref{p:funspace}.\\
$H^s$ & Abbreviation for $H^{s,2}$.\\
$\H_\c$ & The Hessian of a configuration $\c \in \fC(Y)$, p. \pageref{p:Hess}.\\
$\tH_\c$ & The augmented Hessian of a configuration $\c$, p. \pageref{p:aHess}.\\
$\J_\c$, $\J_{\c,t}$ & The subspace of $\T$ given by the infinitesimal action of $\G(Y)$ and $\G_{\partial}(Y)$, respectively, p. \pageref{p:J}. \\
$J_\Sigma$ & The compatible complex structure on $\T_\Sigma$, p. \pageref{p:JSigma}\\
$\tJ_\Sigma$ & The compatible complex structure on $\tT_\Sigma$, p. \pageref{p:tJSigma}.\\
$\K_\c$, $\K_{\c,n}$ & The orthogonal complement of $\J_\c$ and $\J_{\c,t}$ in $\T_{\c}$, respectively, p. \pageref{p:K}.\\
$\K(Y)$ & The bundle over $\fC(Y)$ whose fiber over $\c$ is $\K_\c$, p. \pageref{p:K(Y)}.\\
$\L(Y)$, $\L^{s-1/p,p}(Y)$ & The tangential boundary values of the space of monopoles $\fM$ and
    $\fM^{s,p}$ on $Y$, respectively, p. \pageref{p:L}.  \\
$\fM(Y)$, $\fM^{s,p}(Y)$ & The space of all monopoles in $\fC(Y)$ and $\fC^{s,p}(Y)$, respectively, p. \pageref{p:fM}. \\
$\M(Y)$, $M^{s,p}(Y)$ & The space of all monopoles on $Y$ in $\fC(Y)$ and $\fC^{s,p}(Y)$, respectively, that are in global Coulomb gauge, p. \pageref{p:fM}.\\
$\nu$ & The outward unit normal vector field to $Y$, p. \pageref{p:nu}.\\
\end{tabularx}

\begin{tabularx}{450pt}{l|X}
$\omega$ & The symplectic form on $\T_\Sigma$, p. \pageref{omega}.\\
$\tilde\omega$ & The symplectic form on $\tT_\Sigma$, \pageref{tomega}.\\
$P^+_\c$ & The ``Calderon projection'' of the Hessian $\H_\c$, p. \pageref{p:P}.\\
$\tP^+_\c$ & The Calderon projection of the augmented Hessian $\tH_\c$, p. \pageref{p:tP}.\\
$P_\c$ & The ``Poisson operator" of the Hessian $\H_\c$, p. \pageref{p:P}.\\
$\tP_\c$ & The Poisson operator of the augmented Hessian $\tH_\c$, p. \pageref{p:tP}.\\
$\Pi_{\K_\c}$ & The projection onto $\K_\bullet$ through $\J_{\bullet,t}$, p. \pageref{p:PiK}.\\
$\Psi$ & A spinor, p. \pageref{p:s}.\\
$\rho$ & Clifford multiplication on $Y$, \pageref{p:s}.\\
$\q$ & The quadratic map associated to the Seiberg-Witten map $SW_3$, p. \pageref{p:q}.\\
$r_\Sigma$  & The tangential restriction map, p. \pageref{p:rSigma}, \pageref{p:rSigma2}.\\
$r$ & The full restriction map, p. \pageref{p:r}.\\
$\#$ & Some pointwise bilinear multiplication operation.\\
$*$ & The Hodge star operator on $Y$.\\
$\check{*}$ & The Hodge star operator on $\Sigma = \partial Y$. \\
$\s$ & A $\spinc$ structure, p. \pageref{p:s}\\
$\S$ & The spinor bundle on $Y$, p. \pageref{p:s}\\
$SW_3$ & The Seiberg-Witten map in $3$-dimensions, \pageref{p:SW3}.  \\
$\T_\bullet$ & The tangent space $T_\bullet\fC(X)$ for a configuration $\bullet \in \fC(X)$, p. \pageref{p:Tc}.\\
$\T$, $\T^{s,p}$ & The space $\Omega^1(Y; i\R) \oplus \Gamma(\S)$ and its $B^{s,p}(Y)$ closure, isomorphic to any tangent space of $\fC(Y)$ and $\fC^{s,p}(Y)$, respectively, p. \pageref{p:T}.\\
$\T_\Sigma$, $\T^{s,p}_\Sigma$ & The space $\Omega^1(\Sigma; i\R) \oplus \Gamma(\S_\Sigma)$ and its $B^{s,p}(\Sigma)$ closure, isomorphic to any tangent space of $\fC(\Sigma)$ and $\fC^{s,p}(\Sigma)$, respectively, p. \pageref{p:TSigma}.\\
$\tT$ & The augmented space $\tT \oplus \Omega^0(Y;i\R)$, p. \pageref{p:tT}.\\
$\tT_\Sigma$ & The augmented space $\T_\Sigma \oplus \Omega^0(\Sigma; i\R) \oplus \Omega^0(\Sigma; i\R)$, p. \pageref{p:tTSigma}.\\
$Y$ & A $3$-manifold.  \\
$X$ & A manifold or a Banach space.\\
$\tilde X^{s,p}_\c$, $X^{s,p}_\c$ & Subspaces of $\tT^{s,p}$ and $\T^{s,p}$ on which $\tH_\c$ and $\H_\c$ are invertible, respectively, p. \pageref{p:tX}.
\end{tabularx}


\end{document}